\documentclass[11pt,a4paper]{article}
\usepackage[margin=27mm]{geometry}
\usepackage{amsmath,amssymb,amsthm,mathtools,microtype,mathrsfs}
\usepackage[colorlinks=true,linkcolor=blue,citecolor=blue,urlcolor=blue]{hyperref}
\allowdisplaybreaks
\numberwithin{equation}{section}
\newtheorem{theorem}{Theorem}[section]
\newtheorem{lemma}{Lemma}[section]
\newtheorem{proposition}[lemma]{Proposition}
\theoremstyle{definition}
\newtheorem{definition}{Definition}[section]
\theoremstyle{remark}
\newtheorem{remark}[lemma]{Remark}
\newcommand{\Om}{\Omega}
\newcommand{\dd}{\,\mathrm dx}
\newcommand{\ds}{\,\mathrm ds}
\newcommand{\diver}{\operatorname{div}}
\newcommand{\curl}{\operatorname{curl}}
\newcommand{\Dt}{D_t}
\newcommand{\eps}{\varepsilon}
\title{\bfseries Global strong solutions and large-time behavior\\ for the two-dimensional compressible\\ Navier–Stokes–Poisson system with general bulk viscosity and Navier slip boundary conditions}
\author{Hongyu Wang$^a$, Rong Zhang$^b$\thanks{Corresponding author: Rong Zhang.
Email addresses: \texttt{hyuwang\_a@163.com} (H. Y. Wang);
\texttt{rzhang0921@gmail.com} (R. Zhang).}\\[3mm]
{\small $^a$ School of Mathematics and Computer Science,}\\
{\small Nanchang University, Nanchang 330031, P. R. China}\\[1mm]
{\small $^b$ Institute of Mathematics and Interdisciplinary Sciences,}\\
{\small Nanchang University, Nanchang 330031, P. R. China}}
\date{}
\begin{document}
\maketitle
\begin{abstract}
The global existence of a unique strong solution to the two-dimensional
compressible Navier--Stokes--Poisson system with Navier slip boundary
conditions is established in a smooth bounded domain that may be
multiply connected. Arbitrarily large initial data and vacuum are
allowed, and the positive doping profile may have large spatial
variations. The bulk viscosity satisfies a logarithmic derivative
bound and a growth condition on its primitive, including every
superlinear power law and near-linear laws with a logarithmic gain
larger than one half. Time-uniform density and velocity bounds and
exponential convergence to the nonconstant hydrostatic state are
obtained. For strictly positive initial density, uniform positivity and
exponential convergence to the hydrostatic state in higher Sobolev norms
are proved. An initial vacuum point persists and causes
exponential growth of density gradients.
\end{abstract}

\section{Introduction and main results}
The compressible Navier--Stokes--Poisson system describes the motion of
charged particles under a self-consistent electric field;
see Markowich--Ringhofer--Schmeiser~\cite{MarkowichRinghoferSchmeiser1990}.
We study the
two-dimensional problem in a bounded connected smooth domain
$\Om\subset\mathbb R^2$, which may be multiply connected:
\begin{equation}\label{nsp:model}
\begin{cases}
 \rho_t+\diver(\rho u)=0,\\
 \rho\dot u+\nabla\rho^\gamma
 =\mu\Delta u+\nabla((\mu+\lambda(\rho))\diver u)+\rho\nabla\Phi,\\
 \Delta\Phi=\rho-b(x).
\end{cases}
\end{equation}
The unknowns $\rho\ge0$, $u=(u_1,u_2)^{\mathsf T}$, and $\Phi$ are
the density, velocity, and electric potential, respectively. The
superscript $\mathsf T$ denotes transpose. The prescribed function
$b$ is the doping profile. The pressure is
$P=\rho^\gamma$, where $\gamma>1$, the shear viscosity $\mu>0$ is
constant, and the bulk viscosity $\lambda$ depends on the density.
We write $D_t=\partial_t+u\cdot\nabla$ and $\dot u=D_tu$.

Let $n$ be the unit outward normal to $\Om$ on each component of
$\partial\Om$. We fix the conventions
\[
 J=\begin{pmatrix}0&1\\-1&0\end{pmatrix},\qquad
 \tau=Jn,\qquad
 \omega=\partial_2u_1-\partial_1u_2,\qquad \nabla^\perp=J\nabla.
\]
Here $\tau$ is the unit tangent and $\omega$ is the scalar vorticity.
For a vector field $v$, $v_\tau=v\cdot\tau$ is its scalar tangential
component and $(\nabla v)_{ij}=\partial_jv_i$. We use
$\partial_n=n\cdot\nabla$, $\partial_\tau=\tau\cdot\nabla$, and the
signed curvature $\partial_\tau n=\kappa\tau$. For an integrable
scalar function $g$, let
$\bar g=|\Om|^{-1}\int_\Om g\dd$.

Let $\rho_0$ and $u_0$ denote the prescribed initial density and
velocity. We supplement \eqref{nsp:model} with the initial conditions
\begin{equation}\label{nsp:initial}
	\rho|_{t=0}=\rho_0,\qquad (\rho u)|_{t=0}=\rho_0u_0,
\end{equation}
and the Navier slip boundary conditions
\begin{equation}\label{nsp:boundary}
 u\cdot n=0,\qquad \omega=-\alpha u_\tau,\qquad
 \partial_n\Phi=0,\qquad \bar\Phi=0,\qquad\text{on } \partial{\Omega},
\end{equation}
where $\alpha\in C^\infty(\partial\Om)$ satisfies
$\alpha\ge0$ and $\alpha\not\equiv0$. With the deformation tensor
$D(u)=(\nabla u+\nabla u^{\mathsf T})/2$, the vorticity condition
is equivalent to the classical Navier stress law
\[
 2\mu D(u)n\cdot\tau+\mu(\alpha+2\kappa)u_\tau=0.
\]
The coefficient $\mu(\alpha+2\kappa)$ is the scalar boundary friction
coefficient in this stress law. Indeed, $2D(u)n\cdot\tau=\omega-2\kappa u_\tau$ when $u\cdot n=0$;
see Coron~\cite{Coron1996}. Thus the sign assumption concerns the
curvature-corrected friction. On a domain with holes, the nontrivial
boundary term also controls tangent vector fields with zero divergence
and curl, as expressed in \eqref{uni:coercive}.

The doping profile and initial charge satisfy
\begin{equation}\label{nsp:doping}
\begin{gathered}
 b\in W^{1,q}(\Om),\qquad 2<q<\infty,\qquad
 0<b_-\le b(x)\le b_+,\\
 \int_\Om\rho_0\dd=\int_\Om b\dd=M>0.
\end{gathered}
\end{equation}
Here $b_-$ and $b_+$ are fixed positive constants, and $M$ is the
conserved total mass. The mass equality in \eqref{nsp:doping} is the compatibility condition
for the electric Neumann problem. It is preserved by
\eqref{nsp:model}$_1$ and $u\cdot n=0$. Consequently, $\Phi$ is
determined by $\rho$ at every time and is not an independent initial
datum. No smallness is imposed on the initial energy or on the spatial
variation of $b$.

The global existence of large solutions to the two-dimensional
compressible Navier--Stokes equations depends strongly on the bulk
viscosity. For $\lambda(\rho)=\rho^\beta$,
Vaigant--Kazhikhov~\cite{VaigantKazhikhov1995} established global
strong solutions with positive periodic density when $\beta>3$.
Huang--Li~\cite{HuangLi2016} subsequently obtained global strong
and weak solutions for arbitrarily large periodic data with vacuum
under the weaker restriction $\beta>4/3$. For the Cauchy problem
with far-field vacuum, global strong and classical solutions were
also obtained under the same restriction~\cite{HuangLi2022}.
These proofs rely on estimates for the effective viscous flux to
control the density.

In bounded domains, this control also involves the boundary terms
in the effective viscous flux.
Fan--Li--Li~\cite{FanLiLi2022} proved global existence of strong
and weak solutions for large data with vacuum in simply connected domains,
with Navier slip boundary conditions and $\beta>4/3$.
For multiply connected domains, time-uniform bounds and exponential
convergence to the constant steady state were obtained by
Fan--Jiang--Li~\cite{FanJiangLi2024} under nontrivial boundary dissipation.
More recently, Fan--Jiang--Li--Lv~\cite{NS} established global
strong solutions and exponential convergence to the constant
steady state for every $\beta>1$
in smooth bounded simply connected domains with zero boundary
vorticity. The same authors~\cite{NSPeriodicCauchy} treated general
bulk viscosities in the periodic and Cauchy problems. Near-linear
growth with a logarithmic gain is allowed under their conditions. The density
estimate uses fixed Orlicz moments and the dissipation needed for
successive increases of the viscosity primitive.

For the compressible Navier--Stokes--Poisson system, the Poisson
equation couples the density to a nonlocal force. The local theory
already allows vacuum: Tan--Zhang~\cite{TanZhang2010} proved local
existence of strong solutions for the three-dimensional isentropic system under
a compatibility condition, while Yin--Tan~\cite{YinTan2009}
obtained local strong solutions for the related full heat-conducting
system in a self-gravitating setting, also under compatibility
conditions. Global weak solutions have been constructed under
various pressure and viscosity assumptions. Donatelli~\cite{Donatelli2003}
proved global existence of weak solutions for the Cauchy problem, and
Zhang--Tan~\cite{ZhangTan2007} obtained global finite-energy weak
solutions in two dimensions with logarithmic pressure growth.
With density-dependent viscosities,
Ducomet--Ne\v{c}asov\'a--Vasseur~\cite{DucometNecasovaVasseur2010}
studied global weak solutions for a self-gravitating gas, whose
force is attractive. Ye--Dou~\cite{YeDou2017} established global
weak solutions on the three-dimensional torus with density-dependent
viscosity and nonmonotone pressure. These weak-solution results use
pressure laws or viscosity assumptions different from those considered here.

For the Cauchy problem, global strong solutions near a constant
equilibrium were obtained by Hao--Li~\cite{HaoLi2009} in Besov
spaces in three and higher dimensions. In three-dimensional whole space,
Li--Matsumura--Zhang~\cite{LiMatsumuraZhang2010} established global
solutions and optimal algebraic decay near a constant equilibrium;
optimal decay estimates for higher spatial derivatives were
obtained by Wang~\cite{Wang2012} using energy estimates and
negative Sobolev norms.
For nonconstant doping, the steady state is generally nonconstant.
Tan--Wang--Wang~\cite{TanWangWang2015} proved global classical
solutions near such states in three-dimensional whole space,
allowing large spatial variations of the doping profile.
Global strong solutions with vacuum and large oscillations were
obtained by Liu--Xu--Zhang~\cite{LiuXuZhang2020,LiuXuZhang2021}
under a smallness assumption on the initial energy.
Xu--Zhang--Zhong~\cite{XuZhangZhong2022} obtained global strong
solutions for the full heat-conducting system with small initial
energy, large oscillations and vacuum, again allowing large variations
of the doping profile. In these existence results with nonconstant
doping, the smallness is imposed on the initial perturbation or
energy, rather than on the variation of the fixed doping profile.

The initial-boundary value problem requires control of both the
velocity and the electric potential at the boundary. In two-dimensional
bounded domains, Li--Zhou--Liao~\cite{LiZhouLiao2018} obtained global
regular solutions for well-prepared data and justified the
zero-electron-mass limit with estimates uniform in time and in the
electron mass. In three-dimensional bounded domains,
Liu--Zhong~\cite{LiuZhong2021} proved global existence and
exponential stability of smooth solutions near a positive nonconstant steady state with
no-slip velocity and Neumann electric conditions. Large variations
of the doping profile and steady density are allowed.
For slip boundary conditions, Liu--Si--Yu~\cite{LiuSiYu2023}
established global strong solutions and decay with constant doping
and small initial energy, allowing vacuum.
Global classical solutions with small initial energy, vacuum, and
large variations of both the doping profile and the steady state
were obtained by Chen--Huang--Shi~\cite{ChenHuangShi2024}.
Wu~\cite{Wu2025} proved global strong and smooth solutions near
a nonconstant steady state with large variations of the doping profile.

There are also results for boundary problems in unbounded domains.
Liu--Luo--Zhong~\cite{LiuLuoZhong2020} obtained global solutions
for large radially symmetric data outside a ball; for exterior domains
without the symmetry assumption, their stability result near a constant
state includes an additional damping term.
Du--Liu~\cite{DuLiu2025} treated the problem without damping outside
a ball in three dimensions, with general Navier slip and Neumann
electric conditions. They obtained global strong solutions near a
nonconstant steady state for doping profiles with large variations.
Without damping, Liu--Wu--Zhong~\cite{LiuWuZhong2026} established
global existence of strong solutions near a positive constant state in the
three-dimensional half-space and smooth exterior domains with
no-slip velocity and Neumann electric conditions. Additional
low-frequency or integrability assumptions on the initial perturbation
yield algebraic decay.

For arbitrarily large initial data, global strong solutions have
been obtained under further assumptions on the doping, viscosity or
symmetry. In the three-dimensional Cauchy problem,
Yu~\cite{Yu2025} proved global existence and decay of strong solutions with
density-dependent viscosities when the nonconstant doping profile
is sufficiently large. On solid balls,
Fan--Huang--Lei~\cite{FanHuangLei2026} established global classical
solutions for arbitrarily large spherically symmetric data bounded
away from vacuum, with density-dependent shear and bulk viscosities.
For two-dimensional bounded simply connected $C^{2,1}$ domains,
Du--Zhao~\cite[Theorems~1.2 and~1.3]{DuZhao} obtained global strong
and weak solutions with arbitrarily large data and vacuum for
$\lambda(\rho)=\rho^\beta$, $\beta>4/3$. Their results allow
nonconstant positive doping, Navier slip for the velocity and the
Neumann condition for the electric potential. The proof retains
the electric contribution and the harmonic boundary term in the
effective viscous flux. Concerning large-time behavior,
Xu~\cite{Xu2026} proved exponential convergence to a nonconstant
hydrostatic state in three-dimensional bounded domains, allowing
large nonconstant doping and vacuum, under an essential upper
density bound.

In this paper, we establish global existence and large-time behavior
of strong solutions to the planar Navier--Stokes--Poisson system with the
general bulk viscosities specified below. The initial data may be
arbitrarily large and may contain vacuum; the positive doping profile
may have large spatial variations, and the domain may have holes.
In addition to global existence, we obtain density and velocity
bounds independent of time and exponential convergence to the
nonconstant hydrostatic state. In particular, the density upper bound
is proved from the initial data. Compared with the strong-solution
result of Du--Zhao~\cite{DuZhao}, we allow a wider class of bulk
viscosities and multiply connected domains, while assuming
$\alpha\not\equiv0$ in addition to $\alpha\ge0$.
Their weak-solution conclusion is not extended here. The zero
boundary vorticity condition in~\cite{NS} gives a different boundary
coercivity from the condition considered in this paper.

We assume
\begin{equation}\label{nsp:viscosity}
\begin{gathered}
 \lambda\in C^1([0,\infty))\cap C^2((0,\infty)),\qquad
 \lambda\ge0,\qquad |r\lambda'(r)|\le C_\lambda(1+\lambda(r))
       \quad(r>0),\\
 \theta(r)=2\mu\log r+\int_1^r\frac{\lambda(z)}{z}\,\mathrm dz,
 \qquad
 \theta(r)\ge c_\lambda r[\log(e+r)]^a
       \quad(r\ge r_\lambda),
\end{gathered}
\end{equation}
where $a\in(1/2,1]$, $C_\lambda,c_\lambda>0$, and
$r_\lambda\ge2$ are fixed constants.
These are the general viscosity conditions used for the periodic and
Cauchy problems in~\cite{NSPeriodicCauchy}. Their logarithmic
derivative bound gives
\begin{equation}\label{nsp:polynomial}
 2\mu+\lambda(r)\le C(1+r^\ell),\qquad
 \lambda(r)\le C(1+\theta(r)_+),\qquad
 \ell=\max\{1,C_\lambda\}.
\end{equation}
Thus the upper power in \eqref{nsp:polynomial} is a consequence of
\eqref{nsp:viscosity}, not an additional growth assumption. We use $b$
only for the doping profile and $\ell$ for this upper exponent.

The limiting hydrostatic pair $(\rho_s,\Phi_s)$ satisfies
\begin{equation}\label{nsp:steady-model}
 \nabla\rho_s^\gamma=\rho_s\nabla\Phi_s,\quad
 \Delta\Phi_s=\rho_s-b,\quad
 \partial_n\Phi_s=0,\quad\bar\Phi_s=0.
\end{equation}
Lemma~\ref{nsp:stationary} gives its existence and uniqueness, together
with $b_-\le\rho_s\le b_+$ and the prescribed total mass. We write
$P_s=\rho_s^\gamma$. The variation of this fixed state cannot be
treated as a small coefficient in the estimates.

We use the relative effective viscous flux
\begin{equation}\label{nsp:intro-relative-flux}
 G=(2\mu+\lambda(\rho))\diver u-(P-P_s).
\end{equation}

Before stating the main results, we specify the function-space
notation and the solution class.
For $1\le p\le\infty$ and a nonnegative integer $k$, we use
\[
 L^p=L^p(\Om),\qquad W^{k,p}=W^{k,p}(\Om),\qquad
 H^k=W^{k,2}(\Om).
\]
The same notation is used componentwise for vector and matrix fields.
Unless another domain is displayed, $\|\cdot\|_{L^p}$,
$\|\cdot\|_{W^{k,p}}$, and $\|\cdot\|_{H^k}$ are the corresponding
norms on $\Om$. We use the usual Bochner spaces in time and write
$z_+=\max\{z,0\}$. The function $\theta$ is defined for $r>0$ and
tends to $-\infty$ as $r\to0^+$. Accordingly, $\theta(\rho)_+$ is assigned the value zero at vacuum
by continuous extension. The symbols $C$ and $c$, with or without
subscripts, denote positive constants whose dependence is stated with
the relevant result. Generic bounds denoted by $C$ may change from
line to line; rates in a stated exponential bound are fixed.
The structural constants in \eqref{nsp:viscosity} retain their values.

\begin{samepage}
\begin{definition}[Strong solution]\label{nsp:def-strong}
Let $T>0$ and let $q$ be as in \eqref{nsp:doping}. A triple
$(\rho,u,\Phi)$ is a strong solution of
\eqref{nsp:model}--\eqref{nsp:initial} on $[0,T]$ if $\rho\ge0$,
the equations hold almost everywhere in $\Om\times(0,T)$, the
boundary conditions hold in the Sobolev trace sense for almost every
$t\in(0,T)$, and, for every $1<r<q/(q-1)$,
\begin{equation}\label{nsp:strong-class}
\left\{
\begin{aligned}
 &\rho\in C([0,T];W^{1,q}),\qquad
       \rho_t\in L^\infty(0,T;L^2),\\
 &u\in L^\infty(0,T;H^1)\cap L^2(0,T;H^2)
                   \cap L^r(0,T;W^{2,q}),\\
 &\sqrt t\,u\in L^2(0,T;W^{2,q}),\qquad
       \sqrt t\,u_t\in L^2(0,T;H^1),\\
 &\sqrt\rho\,u_t\in L^2(0,T;L^2),\qquad
       \rho u,\sqrt\rho\,u\in C([0,T];L^2),\\
 &\nabla\Phi\in C([0,T];W^{2,q}),\qquad
       \nabla\Phi_t\in L^\infty(0,T;H^1).
\end{aligned}
\right.
\end{equation}
The initial conditions in \eqref{nsp:initial} are attained through
the continuous density and momentum traces; the continuous weighted
velocity has initial value $\sqrt{\rho_0}u_0$. A global strong solution
has these properties on every finite interval $[0,T]$.
\end{definition}
\end{samepage}

Our first result gives global existence and estimates independent of
time and of a positive lower density bound in the approximations.

\begin{theorem}[Strong solution and uniform bounds]\label{nsp:thm-global}
Let $\Om\subset\mathbb R^2$ be a bounded connected smooth domain,
possibly with holes, and let $\alpha\in C^\infty(\partial\Om)$ satisfy
$\alpha\ge0$, $\alpha\not\equiv0$. Assume
\eqref{nsp:doping} and \eqref{nsp:viscosity},
$0\le\rho_0\in W^{1,q}$, and $u_0\in H^1$ with $u_0\cdot n=0$.
Then \eqref{nsp:model}--\eqref{nsp:initial} admits a unique global
strong solution in the sense of Definition~\ref{nsp:def-strong}.
No compatibility condition involving the initial acceleration is
required. With $\sigma(t)=\min\{1,t\}$,
\begin{equation}\label{nsp:uniform-conclusion}
\begin{aligned}
 &\operatorname*{ess\,sup}_{t>0}\bigl(\|\rho(t)\|_{L^{\infty}}+\|u(t)\|_{H^1}^2\bigr)
       +\int_0^\infty\|\sqrt\rho\,\dot u\|_{L^{2}}^2\,\mathrm dt\\
 &+\operatorname*{ess\,sup}_{t>0}\sigma(t)\|\sqrt\rho\,\dot u(t)\|_{L^{2}}^2
       +\int_0^\infty\sigma(t)\|\nabla\dot u\|_{L^{2}}^2\,\mathrm dt\le C.
\end{aligned}
\end{equation}
The constant depends only on $\Om,\alpha,\mu,\gamma,q,M$,
$\|\rho_0\|_{L^{\infty}}$, $\|u_0\|_{H^1}$, the constants in
\eqref{nsp:viscosity}, a bound for $\lambda$ on $[0,2]$, and the fixed
doping profile through $b_-,b_+$ and $\|b\|_{W^{1,q}}$.
It is independent of a positive lower density bound in the
approximations and of higher initial density norms.
\end{theorem}

The second result describes convergence to $(\rho_s,\Phi_s)$
and distinguishes positive initial density from initial vacuum.

\begin{theorem}[Large-time behavior and vacuum]\label{nsp:thm-longtime}
Under the hypotheses of Theorem~\ref{nsp:thm-global}, let
$(\rho_s,\Phi_s)$ satisfy \eqref{nsp:steady-model} and let $G$ be given
by \eqref{nsp:intro-relative-flux}. For each $1\le p<\infty$,
\begin{equation}\label{nsp:decay-conclusion}
 \|\rho-\rho_s\|_{L^{p}}+\|u\|_{W^{1,p}}
       +\|\nabla(\Phi-\Phi_s)\|_{W^{1,p}}
       \le C_pe^{-c_pt},\qquad t\ge1.
\end{equation}
For almost every $t\ge1$,
\begin{equation}\label{nsp:accel-decay-conclusion}
 \|\sqrt\rho\,\dot u\|_{L^{2}}+\|\sqrt\rho\,u_t\|_{L^{2}}
       +\|G\|_{H^1}+\|\omega\|_{H^1}
       +\left(\int_t^\infty\|\dot u\|_{H^1}^2\,\mathrm ds\right)^{1/2}
       \le Ce^{-ct}.
\end{equation}
If $\min_{\bar\Om}\rho_0>0$, then
$\inf_{(t,x)\in[0,\infty)\times\bar\Om}\rho(t,x)>0$ and
\begin{equation}\label{nsp:positive-conclusion}
\begin{aligned}
 &\|\rho-\rho_s\|_{W^{1,q}}+\|\rho_t\|_{L^{q}}
       +\|u\|_{H^2}+\|u_t\|_{L^{2}}
       +\|\nabla(\Phi-\Phi_s)\|_{W^{2,q}}\le Ce^{-ct},\\
 &\int_t^\infty\bigl(\|u\|_{W^{2,q}}^2+\|u_t\|_{H^1}^2
                         +\|\nabla u\|_{L^{\infty}}\bigr)\,\mathrm ds
       \le Ce^{-ct},\qquad t\ge1.
\end{aligned}
\end{equation}
In \eqref{nsp:positive-conclusion}, pointwise velocity and time-derivative
assertions hold for almost every $t\ge1$; the density and tail assertions
hold for every $t\ge1$. These higher-order constants may also depend on
$\|\rho_0\|_{W^{1,q}}$ and $(\min\rho_0)^{-1}$.
If $\rho_0(x_0)=0$ at some $x_0\in\bar\Om$, vacuum persists and
$\|\rho(t)-\rho_s\|_{L^{\infty}}\ge\min\rho_s$ for every $t\ge0$.
Moreover, for $2<r\le q$ and $0\le\varepsilon<1$,
\begin{equation}\label{nsp:vacuum-conclusion}
\begin{aligned}
 &\|\nabla\rho(t)\|_{L^{r}}\ge C_r^{-1}e^{c_rt},\\
 &|\{x\in\Om:\rho(t,x)\le\varepsilon\rho_s(x)\}|
       \le C_\varepsilon e^{-ct},\qquad t\ge0.
\end{aligned}
\end{equation}
\end{theorem}

\begin{remark}[Admissible bulk viscosities]\label{rem:viscosity}
The assumptions \eqref{nsp:viscosity} include
$\lambda(r)=cr^\beta$ for every $c>0$ and $\beta>1$, as well as
$\lambda(r)=r[\log(e+r)]^\eta$ for every $\eta>1/2$.
The latter grows more slowly than every superlinear power.
Monotonicity of $\lambda$ and a uniform bound for $\lambda''$ near
vacuum are not required. Neither the linear law nor the law with
logarithmic exponent $\eta=1/2$ is covered by these assumptions.
The range is inherited from the general viscosity class
in~\cite{NSPeriodicCauchy}; the electric and boundary estimates
described below are needed to retain it in the present problem.
\end{remark}

We next explain the proof of Theorems~\ref{nsp:thm-global}
and~\ref{nsp:thm-longtime}. The main step is to obtain an upper
bound for the density independent of time and of a positive lower
density bound in the approximations. We begin with fixed density
moments under \eqref{nsp:viscosity}. For a vector field
$f\in L^2(\Om;\mathbb R^2)$, let $Qf\in H^1(\Om)$ be the
mean-zero Neumann potential defined by
\[
 \int_\Om\nabla Qf\cdot\nabla\zeta\dd
       =-\int_\Om f\cdot\nabla\zeta\dd,
 \qquad \overline{Qf}=0\qquad(\zeta\in H^1(\Om)).
\]
In \eqref{mom:transport}, the electric contribution
$Q(\rho\nabla\Phi)$ has no velocity dissipation factor. Bounding
it by a constant would leave a source that is not integrable in
time. We instead retain the pressure term and use it to absorb
this contribution. Write
\[
 A_3^2=\int_\Om\bigl((2\mu+\lambda(\rho))(\diver u)^2
                      +\mu\omega^2\bigr)\dd
                 +\mu\int_{\partial\Om}\alpha u_\tau^2\ds.
\]
Fix a sufficiently large integer $k$. For each $\varepsilon>0$, choose
a fixed $L$ sufficiently large as in Lemma~\ref{nsp:lem-electric-moment}
and put $\psi=Q(\rho u)$ and
$\varphi=(\theta(\rho)-\psi-L)_+$, with the value zero at vacuum
by continuous extension. Then \eqref{nsp:eq-electric-absorption} gives
\[
 \int_\Om\rho\varphi^{k-1}|Q(\rho\nabla\Phi)|\dd
 \le\varepsilon\int_\Om\rho^{\gamma+1}\varphi^{k-1}\dd
    +C_{k,L,\varepsilon}A_3^2
                     \left(1+\int_\Om\rho\varphi^k\dd\right).
\]
The pressure absorbs the first term on the right, and the remaining
term contains $A_3^2$, whose integral over all time is bounded by
the basic energy. To prove \eqref{nsp:eq-electric-absorption}, we use the monotonicity of a
density truncation in the mass measure $\rho\dd/M$ and estimate
the errors introduced by $\psi$. No density maximum is used.
The resulting fixed moment, together with the planar Trudinger
inequality, gives the required control of $\rho u$.
The logarithmic powers require $ak\ge(k+1)/2$, or $(2a-1)k\ge1$;
this is where the restriction $a>1/2$ enters.

To proceed from these moments to a uniform density bound, we need
a relative pressure estimate that is integrable in time. In general,
the hydrostatic state is nonconstant, so
$P_s-\overline{P_s}$ need not vanish, and therefore
$\|P-\bar P\|_{L^2}^2$ need not be integrable over all time. Testing by
$\nabla(\Phi-\Phi_s)/\rho_s$ in \eqref{nsp:weighted-pressure-test}
produces a positive relative pressure term and the remainder
\[
 -\int_\Om
 \frac{P-P_s-\gamma\rho_s^{\gamma-1}(\rho-\rho_s)}{\rho_s^2}
       \nabla\rho_s\cdot\nabla(\Phi-\Phi_s)\dd.
\]
The sign of this remainder is not fixed, and $\nabla\rho_s$ may
be large. We therefore first establish the uniform entry time in
Lemma~\ref{nsp:entry-time}, before closing the density upper bound.
The fixed moments give compactness of the electric field. By taking
unit intervals with small dissipation and using uniqueness of the
nonnegative hydrostatic density, we show that
$\|\nabla(\Phi-\Phi_s)\|_{L^\infty}$ becomes arbitrarily small after
a time independent of the unknown density maximum. This smallness
allows us to absorb the remainder and obtain
\eqref{nsp:relative-pressure}. The viscosity terms at high density
are controlled by the time-integrable tails already obtained in
\eqref{nsp:uniform-tails}.

With this relative pressure bound, we estimate $G$. The energy in
\eqref{nsp:uni-corrected-energy} includes the cross term
$-2\int_\Om(\rho-\rho_s)\nabla\Phi_s\cdot u\dd$, which accounts for
the static part of the electric force. We retain $\dot u$ in the
acceleration estimates, since $\dot u-\nabla\Phi$ need not vanish at
a nonconstant hydrostatic state. The boundary condition gives
$\dot u\cdot n=-\kappa u_\tau^2$, so we use the tangent correction
in \eqref{acc:qbounds}. Quantitative elliptic estimates and
logarithmic interpolation then give \eqref{nsp:uni-local-flux}.
Its factor is proportional to $R\log(e+R)$, where $R\ge2$ is a
temporary upper bound for the density on the interval under
consideration. The large interpolation exponents are chosen in terms
of $R$, and their constants must be tracked explicitly.
Proposition~\ref{uni:apriori-proposition} uses this bound to control
successive doublings of $\theta(\rho)$ along particle paths by the
relative pressure and velocity dissipation. Since the total
dissipation is finite, this gives the desired time-uniform density bound.

Once this density bound is available, the relative energy gives
exponential convergence. For positive initial density, we estimate
$(2\mu+\lambda(\rho))\nabla\rho
 -(2\mu+\lambda(\rho_s))\nabla\rho_s$, so that the stationary gradient
does not leave a constant source. If an initial vacuum point is
present, it persists along the flow. Its coexistence with finite-$L^p$
convergence implies growth of density gradients by interpolation.

The rest of the paper is organized as follows.
Section~\ref{sec:preliminaries} collects the quantitative elliptic,
transport, and interpolation estimates.
Section~\ref{sec:finite-estimates} proves the finite-time bounds.
Section~\ref{sec:uniform-moments} establishes uniform density moments,
the hydrostatic state, and the uniform entry time.
Section~\ref{sec:uniform-estimates} obtains the relative-pressure
dissipation and closes the time-uniform density and velocity bounds.
Section~\ref{sec:construction} constructs the global strong solution,
proves uniqueness, and establishes exponential convergence.
Section~\ref{sec:higher-relaxation} proves convergence in higher
Sobolev norms for positive density and the conclusions concerning
persistent vacuum.

\section{Elliptic and transport estimates}\label{sec:preliminaries}
\subsection{Elliptic estimates}
\label{sec:elliptic}

For a vector field $f$, define the mean-zero scalar $Qf$ by
\begin{equation}\label{mom:Q}
 \int_\Om\nabla Qf\cdot\nabla\phi\dd
       =-\int_\Om f\cdot\nabla\phi\dd,\qquad \overline{Qf}=0.
\end{equation}
Thus $Q\nabla g=-g+\bar g$.

We recall the quantitative elliptic estimates of
Fan, Jiang, Li, and Lv~\cite[Lemma~2.1, (2.6)--(2.9)]{NS}.

\begin{lemma}[Quantitative elliptic estimates]\label{mom:potential}
On the fixed smooth domain, for every $p\ge2$,
\begin{equation}\label{mom:elliptic}
 \|\nabla Qf\|_{L^{p}}\le Cp\|f\|_{L^{p}}.
\end{equation}
The same estimate holds for the gradient of the zero-Dirichlet
solution with divergence data. For the scalar homogeneous Dirichlet
or mean-zero Neumann Poisson problem $-\Delta w=g$,
\begin{equation}\label{ell:Hessian}
 \|D^2w\|_{L^{p}}\le Cp^6\|g\|_{L^{p}} .
\end{equation}
The Neumann equation has the usual mean-zero compatibility condition.
Consequently, tangent vector fields satisfy
\begin{equation}\label{mom:Hodge}
 \|\nabla v\|_{L^{p}}\le Cp^6(\|\diver v\|_{L^{p}}+\|\curl v\|_{L^{p}})
                          +C\|v\|_{L^{2}}.
\end{equation}
For each fixed $1<q<2$, one also has
\begin{equation}\label{mom:potential-low}
 \|Qf\|_{L^{2q/(2-q)}}\le C_q\|f\|_{L^{q}}.
\end{equation}
\end{lemma}

If $\rho\ge0$, $\int_\Om\rho\dd=M>0$, and
$\|\rho\|_{L^{\gamma}}\le C_0$ for $\gamma>1$, then
\begin{equation}\label{ell:weightedHodge}
 \|v\|_{H^1}\le C\bigl(\|\diver v\|_{L^{2}}+\|\curl v\|_{L^{2}}
                                    +\|\sqrt\rho\,v\|_{L^{2}}\bigr),
 \qquad v\cdot n=0.
\end{equation}
The constant depends only on $\Om,M,C_0,\gamma$.

If $\alpha\ge0$ is smooth and not identically zero, then
\begin{equation}\label{uni:coercive}
 \|v\|_{H^1}^2\le C\left(\|\diver v\|_{L^{2}}^2+\|\curl v\|_{L^{2}}^2
                 +\int_{\partial\Om}\alpha v_\tau^2\ds\right),
 \qquad v\cdot n=0,
\end{equation}
with $C$ depending only on $\Om,\alpha$.
On multiply connected domains, \eqref{ell:weightedHodge} and
\eqref{uni:coercive} follow from
Fan, Jiang, and Li~\cite[Lemma~3.4(2) and~(1)]{FanJiangLi2024},
respectively, together with Poincar\'e's inequality for tangent vector
fields. In that lemma take $p=2$, $r=\gamma$, $M_1=M$,
$M_2=C_0^\gamma$, and $K=\alpha$ as applicable.
On simply connected domains, the stronger div--curl coercivity in
Fan, Li, and Li~\cite[Lemma~2.2]{FanLiLi2022} gives both estimates.

\begin{lemma}[Flux and vorticity gradients]\label{mom:fluxelliptic-lemma}
Suppose $\nabla F+\mu\nabla^\perp\omega=f$ and
$\omega=-\alpha u_\tau$ on the boundary. For every $p\ge2$,
\begin{equation}\label{mom:fluxelliptic}
 \|\nabla F\|_{L^{p}}+\|\nabla\omega\|_{L^{p}}
 \le Cp\bigl(\|f\|_{L^{p}}+\|u\|_{W^{1,p}}\bigr).
\end{equation}
\end{lemma}
After a fixed smooth extension of $\alpha\tau$ into $\Om$, apply the
quantitative Dirichlet estimate in Lemma~\ref{mom:potential} to
$\omega+\alpha u\cdot\tau$, which has zero trace, and use
$\nabla F=f-\mu\nabla^\perp\omega$;
see also~\cite[(5.17)--(5.18)]{FanJiangLi2024}.

\begin{lemma}[Transport commutator]
\label{mom:commutator}
Define
\[
 \mathfrak C(u,f)=u\cdot\nabla Qf-Q\diver(f\otimes u),
 \qquad (\diver(f\otimes u))_j=\partial_i(f_j u_i).
\]
For a tangent $u$ and every smooth test function,
\begin{equation}\label{ell:commidentity}
\begin{aligned}
 &\int_\Om\nabla\mathfrak C(u,f)\cdot\nabla\phi\dd\\
 &=\int_\Om\left[(\nabla u+\nabla u^{\mathsf T}-(\diver u)I)\nabla Qf
                            +(\nabla u)f\right]\cdot\nabla\phi\dd.
\end{aligned}
\end{equation}
For each fixed $m>2$,
\begin{equation}\label{mom:commutator-bound}
 \|\mathfrak C(u,f)\|_{L^{m}}\le C_m\|\nabla u\|_{L^{2}}\|f\|_{L^{m}}.
\end{equation}
\end{lemma}
This is the Neumann commutator estimate of
Fan, Jiang, Li, and Lv~\cite[Lemma~2.3, (2.12)--(2.13)]{NS}.

\subsection{The harmonic boundary contribution}

\begin{lemma}[The harmonic boundary contribution]\label{mom:boundary}
Let $B=Q(\nabla^\perp\omega)$. Then, for every finite $p\ge1$,
\begin{equation}\label{mom:boundary-bound}
 \|B\|_{H^1}+\|B\|_{L^{p}}\le C_{\Om,\alpha,p}\|u\|_{H^1}.
\end{equation}
\end{lemma}
This is the harmonic Neumann estimate in
Du and Zhao~\cite[proof of Lemma~3.2, (3.31)--(3.33)]{DuZhao},
with their $h=-\mu B$ and $A=\alpha$.
The weak formulation applies equally to several closed boundary
components.

\subsection{Weighted Sobolev and logarithmic interpolation estimates}

\begin{lemma}[Mass-weighted Poincar\'e and Sobolev inequalities]
\label{ell:weightedPoincare}
If $\rho\ge0$, $\int_\Om\rho\dd=M>0$, and $\|\rho\|_{L^{\gamma}}\le C_0$
for some $\gamma>1$, then every $v\in H^1(\Om)$ satisfies
\begin{equation}\label{ell:weightedSobolev}
 \|v\|_{H^1}\le C(\|\sqrt\rho\,v\|_{L^{2}}+\|\nabla v\|_{L^{2}}),\qquad
 \|v\|_{L^{p}}\le C\sqrt p(\|\sqrt\rho\,v\|_{L^{2}}+\|\nabla v\|_{L^{2}}),
 \quad p\ge2.
\end{equation}
The constant depends only on $\Om,M,C_0,\gamma$.
\end{lemma}
These inequalities are recalled in~\cite[Section~2.1, following
(2.10)]{NS}; see also~\cite[Lemma~2.3]{FanLiLi2022}.

For a zero-mean scalar function the quantitative
Gagliardo--Nirenberg inequality is
\begin{equation}\label{eng:GN}
 \|f\|_{L^{p}}\le C\sqrt p\,\|f\|_{L^{2}}^{2/p}
                      \|\nabla f\|_{L^{2}}^{1-2/p},\qquad p\ge2.
\end{equation}
This is the quantitative planar interpolation inequality in
Fan, Li, and Li~\cite[Lemma~2.3, (2.5)]{FanLiLi2022}, followed by
extension and the zero-mean Poincar\'e inequality. The corresponding
Trudinger bound is
\begin{equation}\label{ell:Trudinger}
 \int_\Om\exp\left(\frac{c|v|^2}{\|v\|_{H^1}^2}\right)\dd\le C
 \quad(v\ne0).
\end{equation}
It follows from the same $C\sqrt p$ Sobolev estimate by expanding the
exponential; no boundary condition on $v$ is required.

\begin{lemma}[Logarithmic interpolation]\label{ell:logInterpolation}
For $\bar f=0$ and $p,m\ge8$,
\begin{equation}\label{cl:interpolation}
 \|f\|_{L^{\infty}}\le C\sqrt p\,
 \|f\|_{L^{2}}^\eta\|\nabla f\|_{L^{2}}^{1-\eta-\sigma}
                           \|\nabla f\|_{L^{m}}^\sigma,
 \qquad
 \eta=\frac{2(m-2)}{p(m-2)+2m},\quad
 \sigma=\frac{2m}{p(m-2)+2m}.
\end{equation}
The constant is independent of $p,m$.
\end{lemma}
This estimate, with $C$ independent of $p,m$, is recalled
from~\cite[Section~3, (3.1)]{NS}.

\section{Estimates on finite time intervals}\label{sec:finite-estimates}
\subsection{Basic energy and the finite-time a priori bound}
\label{nsp:sec-finite-time-transfer}

Consider a smooth positive solution of \eqref{nsp:model} and
\eqref{nsp:boundary} on $[0,S]$, $0<S\le T$.  For the finite-time
argument it is enough to assume $\alpha\ge0$; nontrivial boundary
damping is needed only in the later time-uniform estimates. Put
\begin{equation}\label{eq:flux}
 \rho(\dot u-\nabla\Phi)=\nabla F+\mu\nabla^\perp\omega,
 \qquad F=(2\mu+\lambda(\rho))\diver u-P,\qquad P=\rho^\gamma.
\end{equation}
Define
\begin{equation}\label{eq:energy}
\begin{aligned}
 h&=2\mu+\lambda(\rho),\qquad
 A_3^2=\int_\Om(h(\diver u)^2+\mu\omega^2)\dd
                 +\mu\int_{\partial\Om}\alpha u_\tau^2\ds,\\
 y=A_2^2&=\int_\Om\rho|\dot u|^2\dd,\qquad
 E=e+\int_\Om\left(\mu\omega^2+\frac{F^2}h\right)\dd
                 +\mu\int_{\partial\Om}\alpha u_\tau^2\ds,\\
 R_t&=\max\{2,\sup_{0\le s\le t}\|\rho(s)\|_{L^{\infty}}\}.
\end{aligned}
\end{equation}

\begin{proposition}[Finite-time estimates]\label{thm:main}
Fix $T>0$ and $C_0\ge1$.  Assume \eqref{nsp:viscosity},
$\|\rho_0\|_{L^{\infty}}+\|u_0\|_{H^1}\le C_0$, and
$\int\rho_0=\int b=M>0$.  Then
\begin{equation}\label{eq:main}
\begin{aligned}
 &\sup_{0\le t\le S}
 (\|\rho(t)\|_{L^{\infty}}+\|u(t)\|_{H^1}^2+E(t)+A_3^2(t))
          +\int_0^S A_2^2\,\mathrm dt\\
 &\quad+\sup_{0<t\le S}tA_2^2(t)
          +\int_0^S t\|\nabla\dot u\|_{L^{2}}^2\,\mathrm dt\le C_T.
\end{aligned}
\end{equation}
The constant depends only on $\Om,\alpha,\mu,\gamma,M,C_0,T$,
$\|b\|_{L^{\infty}}$, the constants in \eqref{nsp:viscosity}, and a
bound for $\lambda$ on $[0,2]$. It is independent of $S$, a positive
lower density bound, and higher initial norms.
\end{proposition}

All constants $C_T$ below have this dependence. The derivative bound
on $\lambda$ implies \eqref{nsp:polynomial}: integration gives its
polynomial upper bound, while integration on $[r/2,r]$ uses
$1+\lambda(z)\ge2^{-C_\lambda}(1+\lambda(r))$ to give
$\lambda(r)\le C(1+\theta(r)_+)$. No monotonicity of $\lambda$ is
used.

Testing momentum by $u$ gives
\begin{equation}\label{eq:basic}
 \frac{\mathrm d}{\mathrm dt}\int_\Om
 \left(\frac12\rho|u|^2+\frac{\rho^\gamma}{\gamma-1}
                    +\frac12|\nabla\Phi|^2\right)\dd+A_3^2=0,
 \qquad \|u\|_{H^1}\le C(1+A_3).
\end{equation}
Indeed, $\Delta\Phi_t=-\diver(\rho u)$ and the impermeability and
electric Neumann conditions give
\[
 \int_\Om\rho u\cdot\nabla\Phi\dd
 =\int_\Om\Phi\Delta\Phi_t\dd
 =-\tfrac12\frac{\mathrm d}{\mathrm dt}\|\nabla\Phi\|_{L^{2}}^2.
\]
The viscous boundary term is
$-\mu\int_{\partial\Om}\alpha u_\tau^2\ds$.
The last inequality in \eqref{eq:basic} follows from
\eqref{ell:weightedHodge}, the positive mass, and the kinetic energy.

\subsection{A fixed Orlicz density moment}\label{mom:section}

Fix an integer $k\ge2$ such that
\begin{equation}\label{mom:k}
 k+1\ge4\gamma,\qquad (2a-1)k\ge1,
\end{equation}
where $a>1/2$ is the exponent in \eqref{nsp:viscosity}.

\begin{lemma}[A fixed Orlicz density moment]\label{mom:moment}
For the fixed $k$ in \eqref{mom:k},
\begin{equation}\label{mom:moment-bound}
 \sup_{0\le t\le S}\int_\Om\rho\theta(\rho)_+^k\dd
 +\sup_{0\le t\le S}\int_\Om
       \rho^{k+1}[\log(e+\rho)]^{ak}\dd\le C_T.
\end{equation}
The constant has the dependence in Proposition~\ref{thm:main};
it is independent of the density lower bound and of $R$.
\end{lemma}
\begin{proof}
Put $\psi=Q(\rho u)$, $K=Q(\rho\nabla\Phi)$, and
$B=Q(\nabla^\perp\omega)$.
Applying $Q$ to conservative momentum gives
\[
 \Dt\psi=-F+\bar F+\mu B+\mathfrak C(u,\rho u)+K.
\]
The continuity equation gives $\Dt\theta(\rho)+P=-F$. Hence
\begin{equation}\label{mom:transport}
 \Dt(\theta(\rho)-\psi)+P
    =-\bar F-\mu B-\mathfrak C(u,\rho u)-K.
\end{equation}
Set
\[
 \mathcal M_k(t)=\int_\Om\rho(\theta(\rho)-\psi)_+^k\dd,
 \qquad m=k+1.
\]
The growth condition implies
\[
 \int_\Om\rho^m[\log(e+\rho)]^{ak}\dd
       \le C\left(1+\int_\Om\rho\theta(\rho)_+^k\dd\right).
\]
For $q=2m/(m+1)$, the potential estimate and the basic kinetic energy
bound yield
\[
 \|\psi\|_{L^{2m}}\le C_m\|\rho u\|_{L^{q}}
       \le C_m\|\rho\|_{L^{m}}^{1/2}.
\]
Interpolation with the fixed mass gives
$\|\rho\|_{L^{q}}\le C_M\|\rho\|_{L^{m}}^{1/2}$. Consequently,
\[
 \int_\Om\rho|\psi|^k\dd
  \le\|\rho\|_{L^{q}}\|\psi\|_{L^{2m}}^k
  \le C_m\|\rho\|_{L^{m}}^{m/2}
  \le C_m\left(1+\int_\Om\rho\theta(\rho)_+^k\dd\right)^{1/2}.
\]
Using $\theta_+\le(\theta-\psi)_++|\psi|$ and absorbing this square
root gives
\begin{equation}\label{mom:comparison}
 \int_\Om\rho\theta(\rho)_+^k\dd
 +\int_\Om\rho^m[\log(e+\rho)]^{ak}\dd
       \le C_m(1+\mathcal M_k).
\end{equation}

By \eqref{eq:basic} and the Trudinger inequality,
\[
 \int_\Om\exp\!\left(\frac{c|u|^2}{1+A_3^2}\right)\dd\le C.
\]
Apply the Orlicz--Young inequality
\[
 xy\le C_mx[\log(e+x)]^{m/2}
               +C_m\exp(c_my^{2/m})
\]
with $x=\rho^m$, $y=(c_*|u|/(1+A_3))^m$, and $c_*$ small.
Since $ak\ge m/2$ by \eqref{mom:k}, \eqref{mom:comparison} gives
\begin{equation}\label{mom:momentum}
 \|\rho u\|_{L^{m}}\le C_m(1+A_3)(1+\mathcal M_k)^{1/m}.
\end{equation}
Lemmas~\ref{mom:commutator} and \ref{mom:boundary} give
\[
 \|\mathfrak C(u,\rho u)\|_{L^{m}}
    \le C_m (1+A_3)^2(1+\mathcal M_k)^{1/m},\qquad
 \|B\|_{L^{m}}\le C_m(1+A_3).
\]

For the electric term use the Maxwell identity
\[
 \rho\nabla\Phi=\diver\left(\nabla\Phi\otimes\nabla\Phi
                     -\tfrac12|\nabla\Phi|^2I\right)+b\nabla\Phi.
\]
The stress traction is $-\tfrac12|\nabla\Phi|^2n$ because
$\partial_n\Phi=0$.  It therefore gives no boundary contribution
against the gradient of a homogeneous Neumann dual solution.
Fixed-exponent Neumann duality gives
\[
 \|K\|_{L^{m}}\le C_m(\|\nabla\Phi\|_{L^{2m}}^2+\|b\nabla\Phi\|_{L^{m}}),\qquad
 \|\nabla\Phi\|_{L^{2m}}\le C_m(1+\|\rho\|_{L^{2m/(m+1)}})
          \le C_m(1+\|\rho\|_{L^{m}}^{1/2}).
\]
Consequently
\begin{equation}\label{nsp:ft-K-moment-bound}
 \|K\|_{L^{m}}\le C_m(1+\|\rho\|_{L^{m}})
              \le C_m(1+\mathcal M_k)^{1/m}.
\end{equation}

Test \eqref{mom:transport} by $\rho(\theta-\psi)_+^{k-1}$.
Since $u\cdot n=0$ and $P\ge0$, discard the pressure term and use
\[
 \int_\Om\rho(\theta-\psi)_+^{k-1}|Z|\dd
       \le\mathcal M_k^{1-1/k}\|\rho\|_{L^{m}}^{1/k}\|Z\|_{L^{m}}.
\]
Since $1-1/k+1/[k(k+1)]+1/(k+1)=1$, the commutator contribution
is bounded by $C_m(1+A_3^2)(1+\mathcal M_k)$.
The $B$ term has exponent $k/(k+1)<1$ and satisfies the same bound.
The electric term has the same linear bound by
\eqref{nsp:ft-K-moment-bound} and the displayed exponent identity.

For the mean, since $\int_\Om\diver u\dd=0$,
\[
 \bar F=|\Om|^{-1}\int_\Om\lambda\diver u\dd-\bar P.
\]
By \eqref{nsp:polynomial}, $\theta_+=0$ for $\rho\le1$, and
\eqref{mom:comparison},
\[
 \int_\Om\lambda(\rho)\dd
   \le C\left(1+\left[\int_\Om\rho\theta(\rho)_+^k\dd\right]^{1/k}\right),
 \qquad
 |\bar F|\le C+C_mA_3(1+\mathcal M_k)^{1/(2k)}.
\]
By H\"older's inequality in the mass measure, the mean contribution
has exponent $1-1/(2k)<1$. Hence
\[
 \mathcal M_k'\le C_m(1+A_3^2)(1+\mathcal M_k).
\]
The initial bounds, the potential estimate at a fixed $p>2$, and
Sobolev embedding give
$\|\psi_0\|_{L^{\infty}}+\|\theta(\rho_0)_+\|_{L^{\infty}}\le C$,
so $\mathcal M_k(0)\le C_m$.
Gronwall's inequality, \eqref{eq:basic}, and
\eqref{mom:comparison} prove \eqref{mom:moment-bound}.
\end{proof}

Since $k+1\ge4\gamma$, the preceding lemma gives
\begin{equation}\label{mom:pressure-viscosity}
 \sup_{0\le t\le S}
 \left(\|P(t)\|_{L^{4}}+\int_\Om\lambda(\rho(t))\dd\right)\le C_T.
\end{equation}

Since $m=k+1\ge4\gamma>4$, fixed-exponent Neumann regularity gives
\begin{equation}\label{nsp:ft-electric-fixed-bounds}
 \|\nabla\Phi\|_{L^{\infty}}+\|D^2\Phi\|_{L^{4}}
       +\|\rho\nabla\Phi\|_{L^{2}}+\|\sqrt\rho\,\nabla\Phi\|_{L^{2}}\le C_T.
\end{equation}

Since $F=(2\mu+\lambda)\diver u-P$ and $2\mu+\lambda\ge2\mu$,
\eqref{mom:pressure-viscosity} and \eqref{eq:basic} give
\begin{equation}\label{mom:Ecomparison}
 C_T^{-1}(1+A_3^2)\le E\le C_T(1+A_3^2),\qquad
                  \int_0^S E(t)\,\mathrm dt\le C_T.
\end{equation}
Testing momentum by $u$, including its slip boundary term, gives
\[
 A_3^2=\int_\Om P\diver u\dd-\int_\Om\rho\dot u\cdot u\dd
       +\int_\Om\rho\nabla\Phi\cdot u\dd
       \le C_T A_3+C_T\sqrt y+C_T.
\]
Young's inequality and \eqref{mom:Ecomparison} therefore imply
\begin{equation}\label{mom:Esquare}
                       E^2\le C_T(1+y).
\end{equation}

\subsection{Logarithmic energy estimate}\label{eng:section}

On $[0,S]$, fix
\[
 R=\max\left\{2,\sup_{0\le t\le S}\|\rho(t)\|_{L^{\infty}}\right\}.
\]
Then $h=2\mu+\lambda(\rho)\le CR^\ell$, where $\ell\ge1$ is fixed by
\eqref{nsp:polynomial}.

\begin{lemma}[Low-exponent elliptic estimate]\label{eng:elliptic-lemma}
For every $t\in[0,S]$,
\begin{equation}\label{eng:ellipticlow}
 \|\nabla F\|_{L^{2}}+\|\nabla\omega\|_{L^{2}}
 \le C_T\left(\sqrt{Ry}+\sqrt E\right).
\end{equation}
Moreover, $|\bar F|\le C_T\sqrt E$ and
$\|F\|_{L^{2}}\le CR^{\ell/2}\sqrt E$.
\end{lemma}
\begin{proof}
Apply \eqref{mom:fluxelliptic} at $p=2$ with
$f=\rho(\dot u-\nabla\Phi)$. The bounds
$\|\rho\dot u\|_{L^2}\le\sqrt{Ry}$,
\eqref{mom:Ecomparison}, and \eqref{nsp:ft-electric-fixed-bounds}
give \eqref{eng:ellipticlow}. Since $\int_\Om\diver u\dd=0$,
\[
 |\bar F|\le C\|\lambda(\rho)\|_{L^1}^{1/2}A_3+C\|P\|_{L^1}
             \le C_T\sqrt E,
 \qquad \|F\|_{L^2}^2\le CR^\ell E,
\]
by \eqref{mom:pressure-viscosity} and $h\le CR^\ell$.
\end{proof}

\begin{proposition}[Logarithmic energy]\label{eng:proposition}
There is a positive constant $c$, independent of $R$ and the positive
lower density bound, such that
\begin{equation}\label{eng:differential}
 E'+cy\le C_T R\log(e+R)E^2.
\end{equation}
Consequently,
\begin{equation}\label{eng:logenergy}
 \sup_{0\le t\le S}\log E(t)
 +\int_0^S\frac{e+y}{E}\,\mathrm dt
 \le C_T R\log(e+R).
\end{equation}
\end{proposition}
\begin{proof}
Testing momentum by $2\dot u$ and using
\[
\begin{aligned}
 \diver\dot u&=\Dt(\diver u)+(\diver u)^2-2\det\nabla u,
 &\curl\dot u&=\Dt\omega+\omega\diver u,\\
 \Dt\lambda&=-\rho\lambda'\diver u,
 &\Dt P&=-\gamma P\diver u
\end{aligned}
\]
gives
\begin{equation}\label{eng:identity}
\begin{aligned}
 E'+2y={}&
 \int_\Om\left(1+\frac{\rho\lambda'}h\right)\frac{F^2}h\diver u\dd\\
 &-2\int_\Om\left(1+\frac{\rho\lambda'}h\right)F(\diver u)^2\dd
 +2\gamma\int_\Om\frac{FP}h\diver u\dd\\
 &+4\int_\Om F\det\nabla u\dd
 -\mu\int_\Om\omega^2\diver u\dd\\
 &-2\int_{\partial\Om}\kappa F u_\tau^2\ds
 +\frac{2\mu}{3}\int_{\partial\Om}(\partial_\tau\alpha)u_\tau^3\ds\\
 &+2\int_\Om\rho\nabla\Phi\cdot\dot u\dd .
\end{aligned}
\end{equation}
The boundary terms follow from $\dot u\cdot n=-\kappa u_\tau^2$,
$\dot u\cdot\tau=\partial_tu_\tau+u_\tau\partial_\tau u_\tau$,
and $\omega=-\alpha u_\tau$:
\[
 2\mu\int_{\partial\Om}\omega(\dot u\cdot\tau)\ds
 =-\mu\frac{\mathrm d}{\mathrm dt}
       \int_{\partial\Om}\alpha u_\tau^2\ds
   +\frac{2\mu}{3}\int_{\partial\Om}(\partial_\tau\alpha)u_\tau^3\ds.
\]
The electric work is bounded by
\[
 2\left|\int_\Om\rho\nabla\Phi\cdot\dot u\dd\right|
 \le\varepsilon y+C_\varepsilon\|\sqrt\rho\,\nabla\Phi\|_{L^{2}}^2
 \le\varepsilon y+C_{T,\varepsilon}.
\]

Since $|\rho\lambda'|/h\le C$ and $\diver u=(F+P)/h$, the absolute sum of
the first three terms is at most
\begin{equation}\label{eng:coefficient-terms}
 C\|\diver u\|_{L^{2}}
 \left(\left\|\frac{F^2}h\right\|_{L^{2}}
              +\left\|\frac{FP}h\right\|_{L^{2}}\right).
\end{equation}
For $0<\delta\le1/4$, H\"older's inequality and $h\ge2\mu$ give
\[
 \left\|\frac{F^2}h\right\|_{L^{2}}
 \le C\left(\int_\Om\frac{F^2}h\dd\right)^{(1-\delta)/2}
             \|F\|_{L^{2(1+\delta)/\delta}}^{1+\delta}.
\]
Applying \eqref{eng:GN} to $F-\bar F$ gives
\[
 \|F-\bar F\|_{L^{2(1+\delta)/\delta}}^{1+\delta}
 \le \frac C{\sqrt\delta}
       \|F-\bar F\|_{L^{2}}^\delta\|\nabla F\|_{L^{2}}.
\]
The constant is uniform because
$\delta^{-(1+\delta)/2}\le C\delta^{-1/2}$. Hence
\begin{equation}\label{eng:weighted-quadratic}
 \left\|\frac{F^2}h\right\|_{L^{2}}
 \le \frac{C_T}{\sqrt\delta}R^{\ell\delta/2}
             \sqrt E\,\|\nabla F\|_{L^{2}}+C_TE.
\end{equation}
For the pressure term, use the fixed bound
$\|P\|_{L^{2+\delta}}\le C_T$ and
$q=2(2+\delta)/\delta$ in \eqref{eng:GN}. Weighted arithmetic--geometric
mean then gives
\begin{equation}\label{eng:pressure-quadratic}
 \left\|\frac{FP}h\right\|_{L^{2}}
 \le \frac{C_T}{\sqrt\delta}
       R^{\ell\delta/[2(2+\delta)]}
       \left(\sqrt E+\|\nabla F\|_{L^{2}}\right).
\end{equation}
Using $\|\diver u\|_{L^{2}}\le C\sqrt E$ and \eqref{eng:ellipticlow} in
\eqref{eng:coefficient-terms}--\eqref{eng:pressure-quadratic}, followed
by Young's inequality, bounds these terms by
\begin{equation}\label{eng:cubic-bound}
 \eps y+\frac{C_{T,\eps}}\delta R^{1+\ell\delta}E^2 .
\end{equation}

For the Jacobian term, let $p\ge8$. Lemma~\ref{eng:elliptic-lemma}
and \eqref{eng:GN} give
\[
 \|F\|_{L^{p}}+\|\omega\|_{L^{p}}
 \le C_T\sqrt p\,R^{\ell/p}E^{1/p}
                         (\sqrt E+\sqrt{Ry})^{1-2/p}.
\]
Using \eqref{mom:Hodge}, $\|P\|_{L^{p}}\le CR^\gamma$, and
$\|u\|_{L^{2}}\le C_T\sqrt E$, we obtain
\[
 \|\nabla u\|_{L^{p}}
 \le C_Tp^{13/2}R^{\ell/p}E^{1/p}(\sqrt E+\sqrt{Ry})^{1-2/p}
                                      +Cp^6R^\gamma.
\]
Interpolating $\|\nabla u\|_{L^{2}}\le C_T\sqrt E$ with this bound at
exponent $2p/(p-1)$, followed by H\"older's inequality, gives
\[
\begin{aligned}
 \left|\int_\Om F\det\nabla u\dd\right|
 &\le C_T\sqrt p
      \left(R^{\ell/(p-2)}+R^{\ell/p+2\gamma/(p-2)}\right)\\
 &\qquad\times E(\sqrt E+\sqrt{Ry}).
\end{aligned}
\]
The two powers of $E$ are $1$ and $1-2/[p(p-2)]$, and those of
$\sqrt E+\sqrt{Ry}$ are $1$ and $1-2/p$; both bases are at least one.
Since $p^{13/(p-2)}+p^{12/(p-2)}\le C$ for $p\ge8$, choose
\[
 p=8+(\ell+2\gamma+1)\log(e+R).
\]
The displayed powers of $R$ are now bounded, so
\begin{equation}\label{eng:Jacobian}
\begin{aligned}
 \left|\int_\Om F\det\nabla u\dd\right|
 &\le C_T\sqrt{\log(e+R)}\,E
                    \left(\sqrt{Ry}+\sqrt E\right)\\
 &\le\eps y+C_{T,\eps}R\log(e+R)E^2 .
\end{aligned}
\end{equation}
For the vorticity term, the bounded-domain inequality
$\|\omega\|_{L^{4}}^2\le
C\|\omega\|_{L^{2}}(\|\nabla\omega\|_{L^{2}}+\|\omega\|_{L^{2}})$ gives
\[
 \left|\int_\Om\omega^2\diver u\dd\right|
 \le C_TE\left(\sqrt{Ry}+\sqrt E\right)
 \le\eps y+C_{T,\eps}RE^2.
\]
The trace embedding $H^1(\Om)\to L^3(\partial\Om)$ and the controlled
mean of $F$ similarly imply
\[
\begin{aligned}
 \left|\int_{\partial\Om}\kappa F u_\tau^2\ds\right|
 &\le C\|F\|_{H^1}\|u\|_{H^1}^2
 \le C_TE\left(\sqrt{Ry}+\sqrt E\right)
 \le\eps y+C_{T,\eps}RE^2,\\
 \left|\int_{\partial\Om}(\partial_\tau\alpha)u_\tau^3\ds\right|
 &\le C\|u\|_{H^1}^3\le CE^{3/2}.
\end{aligned}
\]
Take $\delta=[4\log(e+R)]^{-1}$ in \eqref{eng:cubic-bound}.
Then $R^{\ell\delta}\le e^{\ell/4}$. Substitution of all these estimates into
\eqref{eng:identity} and absorption of $\eps y$ prove
\eqref{eng:differential}.

Dividing that inequality by $E$ gives
\[
 (\log E)'+c\frac yE
 \le C_TR\log(e+R)E.
\]
Since $E(0)\le C$ and $\int_0^S E\,\mathrm dt\le C_T$ by
\eqref{mom:Ecomparison}, integration proves \eqref{eng:logenergy};
also $e/E\le1$.
\end{proof}

\begin{lemma}[High-exponent flux gradient]\label{eng:highgradient-lemma}
For every $m\ge8$,
\begin{equation}\label{eng:highgradient}
 \|\nabla F\|_{L^{m}}
 \le C_Tm^{15/2}R^\gamma
       \left(\sqrt{e+y}+\|\nabla\dot u\|_{L^{2}}\right).
\end{equation}
The constant is independent of $m$ and $R$.
\end{lemma}
\begin{proof}
The weighted Sobolev estimate \eqref{ell:weightedSobolev},
\eqref{eng:ellipticlow}, \eqref{mom:Esquare}, and the quantitative
Hodge estimate \eqref{mom:Hodge} give
\[
 \|\rho\dot u\|_{L^m}
   \le C_TR\sqrt m\bigl(\sqrt y+\|\nabla\dot u\|_{L^2}\bigr),
 \qquad
 \|u\|_{W^{1,m}}\le C_Tm^{13/2}R^\gamma\sqrt{e+y}.
\]
By \eqref{nsp:ft-electric-fixed-bounds},
$\|\rho\nabla\Phi\|_{L^m}\le C_TR$, uniformly for $m\ge8$.
Substituting these bounds into \eqref{mom:fluxelliptic} with
$f=\rho(\dot u-\nabla\Phi)$ proves \eqref{eng:highgradient}.
\end{proof}

\subsection{Acceleration estimate}\label{acc:section}

We correct the normal acceleration to obtain a tangent test field.

\begin{lemma}[Positive-time acceleration]\label{acc:lemma}
There is a fixed $N\ge1$, depending only on $\gamma$, the coefficient
and domain constants, such that, for $0<\tau<S\le T$,
\begin{equation}\label{acc:logacceleration}
 \sup_{\tau\le t\le S}\log(e+A_2^2(t))
 +\int_\tau^S
 \left(1+\frac{\|\nabla\dot u\|_{L^{2}}}{\sqrt{e+A_2^2}}\right)^2\,\mathrm dt
 \le C_{\tau,T}R^N.
\end{equation}
The constant is independent of the positive density lower bound and of
$A_2(0)$.
\end{lemma}

\begin{proof}
Keep $R$ fixed on $[0,S]$. The fixed exponent
$N$ may be enlarged in the estimates below.

By \eqref{mom:Hodge}, \eqref{eng:GN}, \eqref{eng:ellipticlow}, and
\eqref{mom:pressure-viscosity},
\begin{equation}\label{acc:gradientfour}
 \|\nabla u\|_{L^{4}}^4
 \le C_T R^N E(y+E)\le C_T R^N E(1+y),
\end{equation}
where the last inequality follows from \eqref{mom:Esquare}.

Extend $n$ constantly along normal segments in a boundary collar and
fix a bounded $C^{1,1}$ extension to $\overline\Om$. Set
\[
 Q_n=\bigl((u\cdot\nabla)n\cdot u\bigr)n,\qquad
 v=\dot u+Q_n,\qquad z=\int_\Om\rho|v|^2\dd.
\]
On the boundary,
\[
 Q_n=\kappa u_\tau^2n,\qquad
 \dot u\cdot n=-\kappa u_\tau^2,\qquad
 v\cdot n=0,\qquad v_\tau=\dot u\cdot\tau.
\]
The extension is fixed. The density moment bound and
$\|u\|_{L^8}\le C_T\sqrt E$ give
$\int_\Om\rho|Q_n|^2\dd\le C\|\rho\|_{L^2}\|u\|_{L^8}^4\le C_TE^2$.
Consequently,
\begin{equation}\label{acc:qbounds}
\begin{aligned}
 |Q_n|&\le C|u|^2,&
 |\nabla Q_n|&\le C\bigl(|u|^2+|u||\nabla u|\bigr),&
 |\Dt Q_n|&\le C\bigl(|u||\dot u|+|u|^3\bigr),\\
 \int_\Om\rho|Q_n|^2\dd&\le C_T E^2,&
 y&\le2z+C_T E^2,&
 z&\le2y+C_T E^2.
\end{aligned}
\end{equation}
By \eqref{acc:gradientfour}, \eqref{mom:Esquare}, and Young's inequality,
\begin{equation}\label{acc:qgradient}
\begin{aligned}
 \|\nabla Q_n\|_{L^{2}}^2
 &\le C_T\bigl(E^2+E\|\nabla u\|_{L^{4}}^2\bigr)\\
 &\le C_T R^N\bigl(E^2+E^{3/2}(1+y)^{1/2}\bigr)
 \le C_T R^N E(1+y).
\end{aligned}
\end{equation}
To differentiate momentum, use $(\nabla u)_{ij}=\partial_j u_i$ and set
\[
 \mathsf S=F I+\mu\omega J,\qquad
 \Sigma=\Dt\mathsf S+(\diver u)\mathsf S
                         -\mathsf S(\nabla u)^{\mathsf T}.
\]
Since $\diver\mathsf S=\rho(\dot u-\nabla\Phi)$, the continuity
equation and the identity
$\diver\Sigma=(\Dt+\diver u)\diver\mathsf S$ give
\begin{equation}\label{acc:commutation}
 \rho\Dt\dot u_i=\partial_j\Sigma_{ij}+\rho\Dt\partial_i\Phi.
\end{equation}
Since $F=h\diver u-P$,
\[
\begin{aligned}
 \Dt F&=h\diver\dot u-\rho\lambda'(\diver u)^2
          -h\operatorname{tr}((\nabla u)^2)+\gamma P\diver u,\\
 \Dt\omega&=\curl\dot u-\omega\diver u.
\end{aligned}
\]
Consequently,
\begin{equation}\label{acc:stressremainder}
\begin{aligned}
 \Sigma&=h(\diver\dot u)I+\mu(\curl\dot u)J+\mathcal Q,\\
 \mathcal Q&=\bigl((h-\rho\lambda')(\diver u)^2
          -h\operatorname{tr}((\nabla u)^2)
          +(\gamma-1)P\diver u\bigr)I\\
  &\qquad-F(\nabla u)^{\mathsf T}-\mu\omega J(\nabla u)^{\mathsf T},\\
 |\mathcal Q|&\le C\bigl(h|\nabla u|^2+P|\nabla u|\bigr).
\end{aligned}
\end{equation}
Define
\[
 H=1+z+E^2.
\]
Equations~\eqref{mom:Esquare} and \eqref{acc:qbounds} imply
\begin{equation}\label{acc:Hcomparison}
 C_T^{-1}(1+y)\le H\le C_T(1+y).
\end{equation}
Test $\rho\Dt v=\diver\Sigma+\rho\Dt Q_n+\rho\Dt\nabla\Phi$ by $v$ and set
$\mathcal D_v=\int_\Om(h(\diver v)^2+\mu(\curl v)^2)\dd$.
Since $v\cdot n=0$, \eqref{ell:weightedHodge} gives
\[
 \|v\|_{H^1}^2\le C(\mathcal D_v+z),\qquad
 \mathcal D_v\ge c\|\nabla v\|_{L^{2}}^2-Cz.
\]
Thus $\eps\|\nabla v\|_{L^{2}}^2\le C\eps\mathcal D_v+C\eps z$.
The lower order terms satisfy
\[
 \sqrt{Ry}\sqrt E\sqrt z\le CR E(y+z),\qquad
 \sqrt z(E+\sqrt E\|\nabla u\|_{L^{4}})
 \le z+C(E^2+E\|\nabla u\|_{L^{4}}^2).
\]
These are bounded by $CR^NEH$ since $z\le H$ and $E\ge e$.
The cross terms containing $Q_n$ are bounded by a small multiple of the
principal dissipation plus $C R^\ell\|\nabla Q_n\|_{L^{2}}^2$. Also,
\[
\begin{aligned}
 \left|\int_\Om \mathcal Q:\nabla v\dd\right|
 &\le\eps\|\nabla v\|_{L^{2}}^2
   +C_\eps R^{2\ell}\|\nabla u\|_{L^{4}}^4
   +C_\eps R^{2\gamma}\|\nabla u\|_{L^{2}}^2\\
 &\le\eps\mathcal D_v+C_{\eps,T}R^N E H.
\end{aligned}
\]
The correction in the material derivative satisfies
\[
\begin{aligned}
 \int_\Om\rho|u||\dot u||v|\dd
 &\le\sqrt{Ry}\,\|u\|_{L^{4}}\|v\|_{L^{4}}
 \le\eps\mathcal D_v+C_{\eps,T}R E(y+z),\\
 \int_\Om\rho|u|^3|v|\dd
 &\le\sqrt{Rz}\,\|u\|_{L^{6}}^3
 \le C_T\sqrt{Rz}\,E^{3/2}
 \le C_T E(z+R E^2).
\end{aligned}
\]

For the new term, differentiating the Poisson equation and using
$\partial_n\Phi_t=0$ gives
\[
 \|\nabla\Phi_t\|_{L^{2}}\le\|\rho u\|_{L^{2}}
       \le\|\rho\|_{L^{4}}\|u\|_{L^{4}}\le C_T\sqrt E,\qquad
 \|u\cdot D^2\Phi\|_{L^{2}}\le C_T\sqrt E.
\]
Consequently
\begin{equation}\label{nsp:ft-electric-acceleration-error}
 \left|\int_\Om\rho\Dt\nabla\Phi\cdot v\dd\right|
 \le\|\rho\|_{L^4}\|\Dt\nabla\Phi\|_{L^2}\|v\|_{L^4}
 \le C_T\sqrt E\,(\mathcal D_v+z)^{1/2}
 \le\varepsilon\mathcal D_v+C_{T,\varepsilon}EH.
\end{equation}
This produces no new boundary term.

For the boundary traction, use $\omega=-\alpha u_\tau$ and
$\Dt u_\tau=v_\tau$ on $\partial\Om$:
\[
\begin{aligned}
 (\nabla u)^{\mathsf T}n
  &=(\partial_nu\cdot n)n-\kappa u_\tau\tau,\\
 \diver u&=\partial_nu\cdot n+\partial_\tau u_\tau,\\
 \Dt\omega&=-\alpha v_\tau-(\partial_\tau\alpha)u_\tau^2.
\end{aligned}
\]
Using $Jn=\tau$ and $v\cdot n=0$ gives
\begin{equation}\label{acc:boundarytraction}
\begin{aligned}
 (\Sigma n)\cdot v
 &=\mu\bigl(\Dt\omega+\omega\partial_\tau u_\tau\bigr)v_\tau
                                      +\kappa F u_\tau v_\tau\\
 &=-\mu\alpha v_\tau^2
   -\mu(\partial_\tau\alpha)u_\tau^2v_\tau
   -\mu\alpha u_\tau(\partial_\tau u_\tau)v_\tau
   +\kappa F u_\tau v_\tau.
\end{aligned}
\end{equation}
The first term is nonpositive, and the trace inequalities give
\[
\begin{aligned}
 \left|\int_{\partial\Om}\kappa F u_\tau v_\tau\ds\right|
 &\le C\|F\|_{H^1}\|u\|_{H^1}\|v\|_{H^1}\\
 &\le\eps\mathcal D_v+C_{\eps,T}R^N E H,\\
 \left|\int_{\partial\Om}(\partial_\tau\alpha)u_\tau^2v_\tau\ds\right|
 &\le C_T E(\|\nabla v\|_{L^{2}}+\sqrt z)
 \le\eps\mathcal D_v+C_{\eps,T}R^N E H.
\end{aligned}
\]
For the remaining term, use fixed smooth extensions of $\alpha,\tau$
and write $u_\tau=u\cdot\tau$, $v_\tau=v\cdot\tau$ in $\Om$.
Stokes' formula gives
\[
\begin{aligned}
 \left|\int_{\partial\Om}\alpha u_\tau(\partial_\tau u_\tau)v_\tau\ds\right|
 &=\frac12\left|\int_{\partial\Om}\alpha v_\tau
                                      \partial_\tau(u_\tau^2)\ds\right|\\
 &\le C\int_\Om|\nabla(\alpha v_\tau)|\,|\nabla(u_\tau^2)|\dd\\
 &\le C_T(\|\nabla v\|_{L^{2}}+\sqrt z)
                  \bigl(E+\sqrt E\,\|\nabla u\|_{L^{4}}\bigr)\\
 &\le\eps\mathcal D_v
       +C_{\eps,T}\bigl(z+E^2+E\|\nabla u\|_{L^{4}}^2\bigr)\\
 &\le\eps\mathcal D_v+C_{\eps,T}R^N E H.
\end{aligned}
\]
Combining these estimates and absorbing $\eps\mathcal D_v$ gives
\begin{equation}\label{acc:zdifferential}
 z'+c\|\nabla v\|_{L^{2}}^2\le C_T R^N E H.
\end{equation}
By \eqref{eng:differential},
\[
 E'+cy\le C_T R^N E^2.
\]
Multiply by $2E$ and add \eqref{acc:zdifferential}.
Since $E^3\le EH$, dropping the nonnegative $Ey$ term gives
\begin{equation}\label{acc:Hdifferential}
 H'+c\|\nabla v\|_{L^{2}}^2\le C_T R^N E H.
\end{equation}
Here $\int_0^S E\,\mathrm dt\le C_T$ by \eqref{mom:Ecomparison}.

Finally, \eqref{eng:logenergy} implies, after enlarging $N$ if needed,
\[
 \int_0^S y\,\mathrm dt\le C_T\exp(C_T R^N).
\]
There exists $s_\tau\in[\tau/2,\tau]$ such that
\[
 y(s_\tau)\le\frac2\tau\int_0^\tau y\,\mathrm dt,
 \qquad \log H(s_\tau)\le C_{\tau,T}R^N,
\]
where the second bound uses \eqref{acc:Hcomparison}.
Divide \eqref{acc:Hdifferential} by $H$ and integrate from $s_\tau$.
Using again \eqref{mom:Ecomparison}, we find
\[
 \sup_{\tau\le t\le S}\log H(t)
       +\int_\tau^S\frac{\|\nabla v\|_{L^{2}}^2}{H}\,\mathrm dt
 \le C_{\tau,T}R^N.
\]
Since $\dot u=v-Q_n$, \eqref{acc:qgradient} and
\eqref{acc:Hcomparison} give
\[
 \frac{\|\nabla\dot u\|_{L^{2}}^2}{e+y}
 \le C_T\frac{\|\nabla v\|_{L^{2}}^2}{H}+C_T R^N E.
\]
Integration, the first inequality in \eqref{acc:Hcomparison}, and
$(1+x)^2\le2+2x^2$ prove \eqref{acc:logacceleration}.
\end{proof}

\subsection{Density upper bound}

\begin{lemma}[Weighted and local flux estimates]\label{cl:flux}
For $0<\tau<S\le T$,
\begin{equation}\label{cl:weightedflux}
 \int_\tau^S\frac{\|F\|_{L^{\infty}}^2}{E}\,\mathrm dt
       \le C_{\tau,T}R_S^2[\log(e+R_S)]^2.
\end{equation}
For $\tau\le s<t\le S$,
\begin{equation}\label{cl:localflux}
 \int_s^t\|F\|_{L^{\infty}}\,\mathrm d\xi
 \le C_{\tau,T}R_t\log(e+R_t)
       \left(\int_s^t(1+A_3^2)\,\mathrm d\xi\right)^{1/2}.
\end{equation}
\end{lemma}
\begin{proof}
Freeze $R=R_S$ and write
\[
 w=1+\frac{e+y}{E},\qquad
 Z=1+\frac{\|\nabla\dot u\|_{L^{2}}}{\sqrt{e+y}}.
\]
From \eqref{eng:ellipticlow},
$\|F\|_{L^{2}}\le CR^{\ell/2}\sqrt E$,
\eqref{eng:highgradient}, and \eqref{cl:interpolation},
\[
 \frac{\|F-\bar F\|_{L^{\infty}}}{\sqrt E}
 \le C_T\sqrt p\,
 R^{1/2+(\ell-1)\eta/2+(\gamma-1/2)\sigma}
          m^{15\sigma/2}w^{(1-\eta)/2}Z^\sigma.
\]
Choose
\[
 p=8+(N+\ell+2\gamma+2)\log(e+R),\qquad m=e^p,
\]
where $N$ is enlarged to cover \eqref{acc:logacceleration}.
Then
\[
 p\eta\le2,\quad p\sigma\le\frac83,\quad
 \sigma-\eta=\frac4{p(m-2)+2m}.
\]
The extra powers of $R$ and $m^{15\sigma/2}$ are bounded. Moreover,
\[
 \sup_{[\tau,S]}w^{(\sigma-\eta)/2}
 \le\exp\left(\frac{C_{\tau,T}R^N}{pm}\right)\le C_{\tau,T},
\]
by \eqref{acc:logacceleration} and $E\ge e$.
Consequently
\[
 \frac{\|F-\bar F\|_{L^{\infty}}}{\sqrt E}
       \le C_{\tau,T}\sqrt{pR}\,
                 w^{(1-\sigma)/2}Z^\sigma.
\]
H\"older in time, \eqref{eng:logenergy}, and
\eqref{acc:logacceleration} give
\[
 \left(\int_\tau^S\frac{\|F-\bar F\|_{L^{\infty}}^2}{E}\,\mathrm dt\right)^{1/2}
 \le C_{\tau,T}\sqrt p\,
 R^{1+(N-1)\sigma/2}[\log(e+R)]^{(1-\sigma)/2}
 \le C_{\tau,T}R\log(e+R).
\]
Since $|\bar F|\le C_T\sqrt E$ by Lemma~\ref{eng:elliptic-lemma},
\eqref{cl:weightedflux} follows. Apply it on $[0,t]$; Cauchy--Schwarz
on $[s,t]$ and \eqref{mom:Ecomparison} give \eqref{cl:localflux}.
\end{proof}

\begin{lemma}[A uniform initial interval]\label{cl:initial}
There are $\tau_0>0$ and $M_0$, depending only on the data in
Proposition~\ref{thm:main}, such that
\[
 \sup_{0\le t\le\min\{S,\tau_0\}}\|\rho(t)\|_{L^{\infty}}\le M_0.
\]
\end{lemma}
\begin{proof}
Work before first contact with $R_*=2(1+C_0)$, with $t\le1$.
By \eqref{eng:logenergy} and \eqref{acc:Hcomparison},
\[
 \sup_{s\le t}E(s)+\int_0^t(y+H)\,\mathrm ds\le C_*.
\]
Multiplying \eqref{acc:Hdifferential} by time and using
\eqref{acc:qgradient} gives
\[
 \sup_{0<s\le t}sy(s)
       +\int_0^t s\|\nabla\dot u\|_{L^{2}}^2\,\mathrm ds\le C_*.
\]
The weighted initial term vanishes for each smooth solution.
By \eqref{eng:highgradient}
with $m=8$,
$\|\nabla F\|_{L^{4}}\le C_*(\sqrt{e+y}+\|\nabla\dot u\|_{L^{2}})$.
Since $\|F\|_{L^{2}}\le C_*$, the GN inequality gives
\[
 \|F\|_{L^{\infty}}
 \le C_*\bigl(1+y^{1/3}+\|\nabla\dot u\|_{L^{2}}^{2/3}\bigr).
\]
H\"older's inequality therefore yields
\[
\begin{aligned}
 \int_0^t\|F\|_{L^{\infty}}\,\mathrm ds
 &\le C_*\left[t+t^{2/3}\left(\int_0^t y\,\mathrm ds\right)^{1/3}
       +t^{1/3}\left(\int_0^t
                   s\|\nabla\dot u\|_{L^{2}}^2\,\mathrm ds\right)^{1/3}\right]\\
 &\le C_*t^{1/3}.
\end{aligned}
\]
Since $u\cdot n=0$, the flow preserves $\overline\Om$, and
\begin{equation}\label{cl:thetaflow}
 \Dt\theta(\rho)+P=-F.
\end{equation}
Thus $\theta(\rho(t))\le\theta(C_0)+C_*t^{1/3}$.
Choose $\tau_0$ with
$\theta(C_0)+C_*\tau_0^{1/3}<\theta(R_*)$ to exclude first contact.
\end{proof}

\begin{proof}[Proof of Proposition~\ref{thm:main}]
If $S\le\tau_0$, use Lemma~\ref{cl:initial}. Otherwise fix
$\tau=\tau_0/2$ and
$M_*>\max\{2,r_\lambda,M_0\}$ with $\theta(M_*)>0$.
Define
\[
 r_j=\theta^{-1}(2^j\theta(M_*)),\qquad t_0=\tau,
\]
and let $t_j$, $j\ge1$, be the first time the running maximum reaches
$r_j$, considering only attained levels. Write
\[
 \mathcal E_j=\int_{t_{j-1}}^{t_j}(1+A_3^2)\,\mathrm dt.
\]
For a trajectory ending at density $r_j$, integration of
\eqref{cl:thetaflow} and \eqref{cl:localflux} give
\[
 \frac{\theta(r_j)}2
 \le C_{\tau,T}r_j\log(e+r_j)\sqrt{\mathcal E_j}.
\]
Since $\theta(r_j)\ge cr_j[\log(e+r_j)]^a$,
$r_j\le C2^j\theta(M_*)$, and $a\le1$, every finite number $L$
of attained levels satisfies
\[
\begin{aligned}
 C_T\ge\sum_{j=1}^{L}\mathcal E_j
 &\ge c_T\sum_{j=1}^{L}
       \left[\frac{\theta(r_j)}{r_j\log(e+r_j)}\right]^2\\
 &\ge c_T\sum_{j=1}^{L}(j+1)^{2a-2}.
\end{aligned}
\]
The last sum diverges for $a>1/2$, so $R_S\le C_T$.

Now \eqref{eng:logenergy} gives $\sup E+\int_0^S y\,\mathrm dt\le C_T$,
and \eqref{mom:Ecomparison} controls $A_3$ and $u$.
Equation \eqref{acc:Hdifferential} has a bounded coefficient and,
by \eqref{acc:Hcomparison}, $\int_0^S H\,\mathrm dt\le C_T$.
Multiplying by time and integrating yields
$\sup tH+\int_0^S t\|\nabla v\|_{L^{2}}^2\,\mathrm dt\le C_T$.
Equations \eqref{acc:Hcomparison} and \eqref{acc:qgradient}
give the remaining terms of \eqref{eq:main}.
\end{proof}

\section{Uniform moments and the hydrostatic state}\label{sec:uniform-moments}
\subsection{Absorption of the electric term in a uniform density moment}
\label{nsp:sec-electric-moment}

Let $\Omega\subset\mathbb R^2$ be a fixed smooth bounded connected
domain, possibly multiply connected.  The electric potential is normalized by
\[
 \Delta\Phi=\rho-b,\qquad \partial_n\Phi=0,\qquad \overline\Phi=0,
 \qquad \int_\Omega\rho\,\mathrm dx=\int_\Omega b\,\mathrm dx=M>0,
\]
where $b\in L^\infty(\Omega)$ is fixed. We use the Neumann operator $Q$ defined in \eqref{mom:Q}.
Thus $Q\nabla g=-g+\bar g$.  In particular, with
\[
 \psi=Q(\rho u),\qquad K=Q(\rho\nabla\Phi),
\]
one has
\begin{equation}\label{nsp:eq-K-neumann}
 \Delta K=-\operatorname{div}(\rho\nabla\Phi),\qquad
 \partial_nK=-\rho\partial_n\Phi=0,\qquad \bar K=0.
\end{equation}
The boundary condition in \eqref{nsp:eq-K-neumann} follows from the
electric Neumann condition; no additional boundary assumption is imposed.

Write
\[
 \theta(r)=2\mu\log r+\int_1^r\frac{\lambda(z)}z\,\mathrm dz,
 \qquad \mu>0,\quad \lambda\ge0.
\]
Assume that $\theta$ is strictly increasing and
\begin{equation}\label{nsp:eq-theta-electric-growth}
 \theta(r)\ge c_\lambda r[\log(e+r)]^a
 \quad(r\ge r_\lambda),\qquad \frac12<a\le1.
\end{equation}
The argument below uses neither a power-law formula for $\lambda$ nor
an upper power bound for $\theta$.

\begin{lemma}[Uniform absorption of the electric moment term]
\label{nsp:lem-electric-moment}
Let $\gamma>1$.  Suppose that $\rho,u,\Phi$ are smooth on a time interval,
$\rho>0$, $u\cdot n=0$, and
\begin{equation}\label{nsp:eq-electric-moment-inputs}
 \int_\Omega\rho|u|^2\,\mathrm dx\le C_0,
 \qquad \|u\|_{H^1}\le C_0 A_3.
\end{equation}
Choose an integer $k$ such that
\begin{equation}\label{nsp:eq-electric-k}
 k\ge7,\qquad k\ge2\gamma+1,\qquad (2a-1)k\ge1,
 \qquad m=k+1.
\end{equation}
For every $\varepsilon>0$ there is $L_\varepsilon>0$ such that, for
each fixed $L\ge L_\varepsilon$, with
\[
 \varphi=(\theta(\rho)-\psi-L)_+,\qquad
 \mathcal M_k=\int_\Omega\rho\varphi^k\,\mathrm dx,
\]
one has
\begin{equation}\label{nsp:eq-electric-absorption}
 \int_\Omega\rho\varphi^{k-1}|K|\,\mathrm dx
 \le\varepsilon\int_\Omega\rho^{\gamma+1}\varphi^{k-1}\,\mathrm dx
       +C_{k,L,\varepsilon}A_3^2(1+\mathcal M_k).
\end{equation}
The constants depend on the fixed domain, $M,C_0,\gamma$, the growth
constants in \eqref{nsp:eq-theta-electric-growth}, and $\|b\|_{L^{\infty}}$.
The choice of $L_\varepsilon$ may additionally depend on the value of
$\theta$ at a fixed finite density determined by these data and
$k,\varepsilon$.  They are independent of the length of the time
interval, a positive lower density bound, and an upper density bound.
\end{lemma}

\begin{proof}
We first record the density and momentum bounds needed for the errors
caused by $\psi$.  Put $X=1+\mathcal M_k$ and
$Y=\int_\Omega\rho\theta(\rho)_+^k\,\mathrm dx$.
For $q=2m/(m+1)$, the fixed-exponent Neumann potential estimate and
the kinetic bound give
\[
 \|\psi\|_{L^{2m}}\le C_m\|\rho u\|_{L^{q}}
 \le C_m\|\rho\|_{L^{m}}^{1/2}
             \|\sqrt\rho\,u\|_{L^{2}}
 \le C_m\|\rho\|_{L^{m}}^{1/2}.
\]
Interpolation with the mass gives $\|\rho\|_{L^{q}}\le
M^{1/2}\|\rho\|_{L^{m}}^{1/2}$, and hence
\[
 \int_\Omega\rho|\psi|^k\,\mathrm dx
 \le\|\rho\|_{L^{q}}\|\psi\|_{L^{2m}}^k
 \le C_m\|\rho\|_{L^{m}}^{m/2}.
\]
By \eqref{nsp:eq-theta-electric-growth},
\[
 \int_\Omega\rho^m[\log(e+\rho)]^{ak}\,\mathrm dx
 \le C(1+Y).
\]
Since $\theta_+\le\varphi+|\psi|+L$, the preceding estimates imply
$Y\le C_{k,L}X+C_k(1+Y)^{1/2}$.  Absorbing the square root yields
\begin{equation}\label{nsp:eq-electric-moment-comparison}
\begin{aligned}
 &\int_\Omega\bigl(\rho\theta(\rho)_+^k
       +\rho^m[\log(e+\rho)]^{ak}\bigr)\,\mathrm dx\le C_{k,L}X,\\
 &\|\rho\|_{L^{m}}\le C_{k,L}X^{1/m},\qquad
 \|\psi\|_{L^{2m}}\le C_{k,L}X^{1/(2m)},\qquad
 \int_\Omega\rho|\psi|^k\,\mathrm dx\le C_{k,L}X^{1/2}.
\end{aligned}
\end{equation}
Neumann regularity, $W^{1,8/3}\hookrightarrow L^\infty$, and
\eqref{nsp:eq-electric-moment-inputs} give
\begin{equation}\label{nsp:eq-psi-infty-electric}
 \|\psi\|_{L^{\infty}}
 \le C\|\rho u\|_{L^{8/3}}
 \le C\|\rho\|_{L^{4}}\|u\|_{L^{8}}
 \le C_{k,L}A_3X^{1/m}.
\end{equation}

The part without $\psi$ can be absorbed entirely by the pressure.
Since $\gamma+1>2$, the Poisson equation and \eqref{mom:Q} give
\begin{equation}\label{nsp:eq-K-infty-electric}
 \|K\|_{L^\infty}
 \le C_\gamma\|\rho\|_{L^{\gamma+1}}\|\nabla\Phi\|_{L^\infty}
 \le C_\gamma(1+\|\rho\|_{L^{\gamma+1}}^2).
\end{equation}
For $L>0$ define
\[
 g(\rho)=(\theta(\rho)-L)_+^{k-1},\qquad
 W_0=\int_\Omega\rho g(\rho)\,\mathrm dx,\qquad
 D_0=\int_\Omega\rho^{\gamma+1}g(\rho)\,\mathrm dx.
\]
If $W_0=0$, then $\|K\|_{L^{\infty}} W_0=0$.  Otherwise use the probability
measure $\mathrm d\nu=\rho\,\mathrm dx/M$.  Since both $s^\gamma$ and
$g(s)$ are increasing,
\begin{equation}\label{nsp:eq-electric-covariance}
\begin{aligned}
 &\int\rho^\gamma g(\rho)\,\mathrm d\nu
       -\left(\int\rho^\gamma\,\mathrm d\nu\right)
        \left(\int g(\rho)\,\mathrm d\nu\right)\\
 &=\frac12\iint
    (\rho(x)^\gamma-\rho(y)^\gamma)
    (g(\rho(x))-g(\rho(y)))\,\mathrm d\nu(x)\,\mathrm d\nu(y)
 \ge0.
\end{aligned}
\end{equation}
In particular,
\[
 \int\rho^\gamma\,\mathrm d\nu\le\frac{D_0}{W_0}.
\]
Writing $A=D_0/W_0$, we therefore have
\[
 \|\rho\|_{L^{\gamma+1}}^2
 =M^{2/(\gamma+1)}
       \left(\int\rho^\gamma\,\mathrm d\nu\right)^{2/(\gamma+1)}
 \le M^{2/(\gamma+1)}A^{2/(\gamma+1)}.
\]
On the support of $g(\rho)$, $\rho>\theta^{-1}(L)$. Hence
\[
 A\ge[\theta^{-1}(L)]^\gamma,\qquad
 \frac{\|K\|_{L^\infty}W_0}{D_0}
 \le C_\gamma\left(A^{-1}+A^{-(\gamma-1)/(\gamma+1)}\right).
\]
Both powers tend to zero as $A\to\infty$. Choose $R_\delta$ so that
the last expression is at most $\delta$ whenever $A\ge R_\delta^\gamma$.
Then every $L\ge\max\{1,\theta(R_\delta)\}$ gives
\begin{equation}\label{nsp:eq-electric-unshifted-absorption}
 \|K\|_{L^{\infty}} W_0\le\delta D_0.
\end{equation}
The choice of $R_\delta$ does not depend on $L$ or the density.

It remains to estimate the errors introduced by $\psi$.  First,
\eqref{nsp:eq-K-infty-electric} and
\eqref{nsp:eq-electric-moment-comparison}, with $\gamma+1\le m$, give
$\|K\|_{L^{\infty}}\le C_{k,L}X^{2/m}$.  By H\"older's inequality in the
mass measure,
\[
 \int_\Omega\rho|\psi|^{k-3}\,\mathrm dx
 \le M^{3/k}
       \left(\int_\Omega\rho|\psi|^k\,\mathrm dx\right)^{(k-3)/k}
 \le C_{k,L}X^{1/2}.
\]
It follows from \eqref{nsp:eq-psi-infty-electric} that
\begin{equation}\label{nsp:eq-electric-psi-error}
\begin{aligned}
 \|K\|_{L^{\infty}}\int_\Omega\rho|\psi|^{k-1}\,\mathrm dx
 &\le\|K\|_{L^{\infty}}\|\psi\|_{L^{\infty}}^2
       \int_\Omega\rho|\psi|^{k-3}\,\mathrm dx\\
 &\le C_{k,L}A_3^2X^{1/2+4/m}
 \le C_{k,L}A_3^2X,
\end{aligned}
\end{equation}
where the last inequality uses $m\ge8$.  Likewise, set
$s=2m(\gamma+1)/(k+5)$.  The choice of $k$ ensures $1<s\le m$ and
$\frac{\gamma+1}s+\frac{k-3}{2m}=1$.  Hence
\begin{equation}\label{nsp:eq-electric-pressure-psi-error}
\begin{aligned}
 \int_\Omega\rho^{\gamma+1}|\psi|^{k-1}\,\mathrm dx
 &\le\|\psi\|_{L^{\infty}}^2\|\rho\|_{L^{s}}^{\gamma+1}
                  \|\psi\|_{L^{2m}}^{k-3}\\
 &\le C_{k,L}A_3^2X^{(k+2\gamma+3)/(2m)}
 \le C_{k,L}A_3^2X.
\end{aligned}
\end{equation}
Here $(k+2\gamma+3)/(2m)\le1$ because $k\ge2\gamma+1$.

Finally, with $c_k=2^{k-2}$,
\[
 \varphi^{k-1}\le c_k\bigl(g(\rho)+|\psi|^{k-1}\bigr),\qquad
 g(\rho)\le c_k\bigl(\varphi^{k-1}+|\psi|^{k-1}\bigr).
\]
Apply \eqref{nsp:eq-electric-unshifted-absorption} with
$\delta=\varepsilon/c_k^2$.  The preceding two inequalities yield
\[
\begin{aligned}
 \int_\Omega\rho\varphi^{k-1}|K|\,\mathrm dx
 &\le c_k\delta D_0
       +c_k\|K\|_{L^{\infty}}\int_\Omega\rho|\psi|^{k-1}\,\mathrm dx\\
 &\le\varepsilon\int_\Omega\rho^{\gamma+1}\varphi^{k-1}\,\mathrm dx
       +\varepsilon\int_\Omega\rho^{\gamma+1}|\psi|^{k-1}\,\mathrm dx
       +c_k\|K\|_{L^{\infty}}\int_\Omega\rho|\psi|^{k-1}\,\mathrm dx.
\end{aligned}
\]
Now \eqref{nsp:eq-electric-psi-error} and
\eqref{nsp:eq-electric-pressure-psi-error} prove
\eqref{nsp:eq-electric-absorption}.
\end{proof}

With $K=Q(\rho\nabla\Phi)$, the transport identity
\eqref{mom:transport} reads
\begin{equation}\label{nsp:eq-electric-theta-transport}
 D_t(\theta(\rho)-\psi)+P
 =-\bar F-\mu B-\mathfrak C(u,\rho u)-K.
\end{equation}
\subsection{Uniform density moments and high-density tails}

We first obtain time-uniform moments without estimating the density
maximum. This order is needed for the entry argument in the next
subsection. Throughout, $\alpha\ge0$ and $\alpha\not\equiv0$.

\begin{lemma}\label{uni:moment-lemma}
Choose an integer $k$ with
\[
 k\ge7,\qquad k\ge2\gamma+1,\qquad (2a-1)k\ge1.
\]
Then, uniformly on every interval of smooth positive existence,
\begin{equation}\label{nsp:uniform-moments}
 \sup_{0\le t\le S}\int_\Om
 \left(\rho\theta(\rho)_+^k+
       \rho^{k+1}[\log(e+\rho)]^{ak}\right)\dd\le C_k.
\end{equation}
Consequently, every fixed finite density moment is uniformly bounded,
and
\begin{equation}\label{uni:momentum}
 \|\rho u\|_{L^{p}}\le C_pA_3,\qquad 2\le p<\infty.
\end{equation}
For each fixed $p\ge1$ there are $L_p,C_p$ such that
\begin{equation}\label{nsp:uniform-tails}
 \int_0^S\int_\Om(\rho-L_p)_+^p\dd\,\mathrm dt\le C_p.
\end{equation}
The constants are independent of $S$, the density maximum, and the
positive lower density bound.
\end{lemma}

\begin{proof}
Set $\psi=Q(\rho u)$, $B=Q(\nabla^\perp\omega)$,
$K=Q(\rho\nabla\Phi)$, and
\[
 \varphi=(\theta(\rho)-\psi-L)_+,\quad
 \mathcal M_k=\int_\Om\rho\varphi^k\dd,\quad
 m=k+1,\quad X=1+\mathcal M_k .
\]
The shift $L$ is a sufficiently large fixed constant, chosen below.
We use the comparison \eqref{nsp:eq-electric-moment-comparison}.
Coercivity and Trudinger's inequality yield
\[
 \int_\Om\exp\left(\frac{c|u|^2}{A_3^2}\right)\dd\le C
       \quad(A_3>0),\qquad u=0\quad(A_3=0).
\]
Apply the Orlicz--Young inequality
\[
 xy\le C_mx[\log(e+x)]^{m/2}+C_m\exp(c_my^{2/m})
\]
with $x=\rho^m$ and $y=(c_*|u|/A_3)^m$, where $c_*>0$ is small.
Since $ak\ge m/2$,
\[
 \|\rho u\|_{L^{m}}\le C_{k,L}A_3X^{1/m}.
\]
Lemmas~\ref{mom:commutator} and \ref{mom:boundary} therefore give
\[
 \|\mathfrak C(u,\rho u)\|_{L^{m}}\le C_{k,L}A_3^2X^{1/m},
 \qquad \|B\|_{L^{2m}}\le C_kA_3 .
\]
For every scalar $Z$,
\[
 \int_\Om\rho\varphi^{k-1}|Z|\dd
 \le\mathcal M_k^{1-1/k}\|\rho\|_{L^{m}}^{1/k}\|Z\|_{L^{m}}.
\]
The commutator term is bounded by $C_{k,L}A_3^2X$, because
$1-1/k+1/(km)+1/m=1$.

Write $V=\overline{\lambda(\rho)\diver u}$.
The basic energy bounds $\bar P$, while
\[
 |V|^2\le CA_3^2\int_\Om\lambda(\rho)\dd
       \le C_{k,L}A_3^2X^{1/k}.
\]
The last inequality follows from
$\lambda(r)\le C(1+\theta(r)_+)$ and
\eqref{nsp:eq-electric-moment-comparison}. Moreover,
$\int_\Om\rho\varphi^{k-1}\dd\le M^{1/k}\mathcal M_k^{1-1/k}$.
Using $|V|\le V^2+1$ gives
\[
 \int_\Om\rho\varphi^{k-1}|V|\dd
 \le C_{k,L}A_3^2X+C\int_\Om\rho\varphi^{k-1}\dd.
\]
Similarly, $\mu|B|\le\mu^2 B^2+1$ and the bound for $B$ give
the same estimate for its contribution. The constants multiplying
the last integral can be chosen independently of $L$.

Test \eqref{nsp:eq-electric-theta-transport} in the transported
mass measure. Lemma~\ref{nsp:lem-electric-moment}, with a sufficiently
small absorption constant, yields
\[
 \mathcal M_k'+\frac{3k}{4}
       \int_\Om\rho^{\gamma+1}\varphi^{k-1}\dd
 \le C_{k,L}A_3^2X+Ck\int_\Om\rho\varphi^{k-1}\dd.
\]
Increase $L$ so that
$L>1+\max\{\theta((8C)^{1/\gamma}),0\}$ and so that the electric
absorption lemma applies. On $\{\rho^\gamma\ge8C\}$ the pressure
absorbs the last term. On its complement $\varphi\le|\psi|$, and
\eqref{nsp:eq-psi-infty-electric} implies
\[
\begin{aligned}
 \int_{\{\rho^\gamma<8C\}}\rho\varphi^{k-1}\dd
 &\le\|\psi\|_{L^{\infty}}^2\int_\Om\rho|\psi|^{k-3}\dd\\
 &\le C_{k,L}A_3^2X^{2/m+1/2}
 \le C_{k,L}A_3^2X .
\end{aligned}
\]
Thus
\begin{equation}\label{nsp:uniform-moment-differential}
 \mathcal M_k'+c_k
       \int_\Om\rho^{\gamma+1}\varphi^{k-1}\dd
       \le C_{k,L}A_3^2(1+\mathcal M_k).
\end{equation}
The initial density and velocity bounds give $\mathcal M_k(0)\le C_k$.
The basic energy bounds $\int_0^S A_3^2$, so Gronwall yields
\begin{equation}\label{uni:moment-damping}
 \sup_{t\le S}\mathcal M_k(t)
       +\int_0^S\int_\Om\rho^{\gamma+1}\varphi^{k-1}
                         \dd\,\mathrm dt\le C_k .
\end{equation}
This proves \eqref{nsp:uniform-moments}. Taking $k$ arbitrarily large
gives every fixed density moment; H\"older's inequality and
$\|u\|_{H^1}\le CA_3$ give \eqref{uni:momentum}.

For $p\ge3$, choose an admissible $k\ge p$ and increase $L_p$ so that
$\theta(r)\ge2L$ and $\theta(r)\ge cr$ for $r\ge L_p$.
On this set either $\varphi\ge\theta(\rho)/4$ or
$\psi\ge\theta(\rho)/4\ge c\rho$. Consequently
\[
\begin{aligned}
 \int_{\{\rho\ge L_p\}}\rho^p\dd
 &\le C_p\int_\Om\rho^{\gamma+1}\varphi^{k-1}\dd
       +C_p\|\psi\|_{L^{\infty}}^2
                           \int_\Om\rho|\psi|^{p-3}\dd\\
 &\le C_p\int_\Om\rho^{\gamma+1}\varphi^{k-1}\dd+C_pA_3^2.
\end{aligned}
\]
Here the now established fixed moments imply
$\|\psi\|_{L^{\infty}}\le C_pA_3$ and a uniform bound for the last
weighted potential moment. Integration proves
\eqref{nsp:uniform-tails} for $p\ge3$.
For $1\le p<3$ increase the threshold by one and use
$(\rho-L_3-1)_+^p\le(\rho-L_3)_+^3$.
\end{proof}
\subsection{The stationary state and a uniform entry time}

Assume that $b\in W^{1,q}(\Om)$ for some $q>2$ and
$0<b_-\le b\le b_+$. The potential is normalized by $\bar\Phi=0$,
and the charge compatibility condition is
$\int_\Om\rho\dd=\int_\Om b\dd$.
The dissipation is
\[
 A_3^2=\int_\Om\bigl((2\mu+\lambda(\rho))(\diver u)^2
                         +\mu\omega^2\bigr)\dd
              +\mu\int_{\partial\Om}\alpha u_\tau^2\ds.
\]
We use $\|u\|_{H^1}\le CA_3$, which follows from
$\alpha\ge0$ and $\alpha\not\equiv0$.

\begin{lemma}[The stationary state]\label{nsp:stationary}
There is a unique strictly positive stationary pair
$(\rho_s,\Phi_s)$ satisfying
\[
 \nabla\rho_s^\gamma=\rho_s\nabla\Phi_s,\qquad
 \Delta\Phi_s=\rho_s-b,\qquad
 \partial_n\Phi_s=0,\qquad \bar\Phi_s=0.
\]
It satisfies
\[
 b_-\le\rho_s\le b_+,\qquad
 \int_\Om\rho_s\dd=\int_\Om b\dd.
\]
Moreover, every nonnegative hydrostatic density with the same mass
and sufficiently high finite Lebesgue integrability equals $\rho_s$;
positivity need not be assumed for the competing density.
\end{lemma}

\begin{proof}
Define $\mathfrak h(r)=\gamma r^{\gamma-1}/(\gamma-1)$ for $r\ge0$.
The constants $\mathfrak h(b_-)$ and $\mathfrak h(b_+)$ are a
subsolution and a supersolution of
\[
 -\Delta w+\mathfrak h^{-1}(w)=b,\qquad \partial_nw=0.
\]
The standard monotone elliptic construction gives a solution between
these constants. Testing the difference of two such solutions by
their difference proves uniqueness. Set
$\rho_s=\mathfrak h^{-1}(w)$ and $\Phi_s=w-\bar w$.
Integration of the equation gives the mass identity, and
$\nabla\rho_s^\gamma=\rho_s\nabla w$ gives the stationary momentum
equation. Since $b\in W^{1,q}$ and $w$ is bounded away from zero,
elliptic regularity gives $w,\rho_s,\Phi_s\in W^{3,q}$.
In particular, their second derivatives are bounded.

For completeness, let $r\ge0$ and $\Psi$ satisfy
\[
 \nabla r^\gamma=r\nabla\Psi,\qquad
 \Delta\Psi=r-b,\qquad \partial_n\Psi=0,\qquad
 \int_\Om r\dd=\int_\Om b\dd,
\]
with $r\in L^p$ for a sufficiently large fixed $p$.
Then $\nabla\Psi\in L^\infty$ and $r^\gamma\in W^{1,p}$,
so $r$ is continuous. On $\{r>0\}$,
$\nabla\mathfrak h(r)=\nabla\Psi$.
The truncations $\max\{\mathfrak h(r),\varepsilon\}$ have gradients
bounded by $\|\nabla\Psi\|_{L^{\infty}}$, independently of $\varepsilon$.
Consequently $z=\mathfrak h(r)$ belongs to $W^{1,\infty}$ and
\[
 \nabla z=\boldsymbol 1_{\{r>0\}}\nabla\Psi
 \quad\hbox{almost everywhere}.
\]
Let $H_\varepsilon$ be a smooth nondecreasing approximation to
$\boldsymbol 1_{(0,\infty)}$, vanishing on $(-\infty,0]$.
For a smooth $\zeta\ge0$, testing the Neumann equation for $\Psi$
by $H_\varepsilon(z)\zeta$ gives
\[
\begin{aligned}
 &\int_\Om H_\varepsilon(z)\nabla\Psi\cdot\nabla\zeta\dd
  +\int_\Om H_\varepsilon'(z)|\nabla\Psi|^2\zeta\dd\\
 &\hspace{35mm}=\int_\Om(b-r)H_\varepsilon(z)\zeta\dd.
\end{aligned}
\]
The second term is nonnegative. Passing to the limit yields
\[
 \int_\Om\nabla z\cdot\nabla\zeta\dd+\int_\Om r\zeta\dd
 \le\int_{\{r>0\}}b\zeta\dd\le\int_\Om b\zeta\dd.
\]
Thus $z$ is a weak Neumann subsolution of
$-\Delta z+\mathfrak h^{-1}(z)=b$.
Subtract the equation for $w$ and test by $(z-w)_+$.
Monotonicity gives $z\le w$, hence $r\le\rho_s$.
Equality of the masses implies $r=\rho_s$ almost everywhere.
\end{proof}

Set $\phi=\Phi-\Phi_s$ and define the relative energy
\begin{equation}\label{nsp:relative-energy}
 \mathcal E_s(t)=\int_\Om\left[
  \frac12\rho|u|^2+
  \frac{\rho^\gamma-\rho_s^\gamma
         -\gamma\rho_s^{\gamma-1}(\rho-\rho_s)}{\gamma-1}
  +\frac12|\nabla\phi|^2\right]\dd.
\end{equation}
The energy calculation, the continuity equation, and
$\mathfrak h(\rho_s)-\Phi_s=\bar w$ give
\begin{equation}\label{nsp:relative-energy-identity}
 \mathcal E_s'(t)+A_3^2(t)=0.
\end{equation}
Indeed, the relative energy minus the uncentered energy equals,
up to a constant,
$-\int_\Om[\mathfrak h(\rho_s)-\Phi_s]\rho\dd$,
which is constant in time. No smallness of $b$ is used.

\begin{lemma}[A uniform entry time]\label{nsp:entry-time}
Suppose that the previously established density moment bounds hold
uniformly on the interval of existence:
\[
 \sup_t\|\rho(t)\|_{L^{p}}\le C_p
 \quad\hbox{for every fixed finite }p.
\]
Assume also the energy bound $\mathcal E_s(0)\le C$ and the
polynomial upper bound $2\mu+\lambda(r)\le C(1+r^\ell)$.
For each $\varepsilon>0$ there is $T_\varepsilon<\infty$,
depending only on these bounds and the fixed coefficients, such that
\[
 \mathcal E_s(t)\le\varepsilon\qquad(t\ge T_\varepsilon)
\]
whenever the solution exists up to that time.
The constants are independent of an upper or a positive lower
pointwise bound for the density. Only finitely many of the moment
bounds are needed for each fixed $\varepsilon$.
\end{lemma}

\begin{proof}
We first prove that sufficiently small dissipation on a unit interval
forces small relative energy at its right endpoint, uniformly for
solutions with the stated bounds. Otherwise there are solutions
$(\rho_j,u_j,\Phi_j)$ on $[0,1]$ such that
\[
 \int_0^1 A_{3,j}^2\,\mathrm dt\longrightarrow0,
 \qquad \mathcal E_{s,j}(1)\ge\varepsilon.
\]
Fix a sufficiently large moment exponent
$p_0>\max\{8,4\gamma,2\ell+4\}$, increasing it if necessary in the
fixed-exponent H\"older estimates below. Coercivity and the moment
bounds imply
\[
 u_j\longrightarrow0\quad\hbox{in }L^2(0,1;H^1),\qquad
 \rho_j u_j\longrightarrow0\quad\hbox{in }L^2(0,1;L^p)
\]
for each of the fixed finite exponents $p$ used below. The continuity
equation gives
$\partial_t\rho_j\to0$ in $L^2(0,1;W^{-1,p})$.
After extraction, $\rho_j\rightharpoonup r$ in the required finite
Lebesgue spaces, with $r$ independent of time.

Neumann regularity gives uniform bounds for $\Phi_j$ in
$L^\infty(0,1;W^{2,p})$, while differentiating its weak equation gives
\[
 \|\nabla\partial_t\Phi_j\|_{L^{p}}\le C_p\|\rho_j u_j\|_{L^{p}}.
\]
In particular, $\nabla\Phi_j$ is bounded in $L^\infty W^{1,p_0}$
and its time derivative tends to zero in $L^2L^4$.
The compact embedding $W^{1,p_0}\Subset C(\overline\Om)$ and
interpolation with this time estimate yield, after extraction,
\[
 \nabla\Phi_j\longrightarrow\nabla\Psi
       \quad\hbox{in }C([0,1];L^\infty),\qquad
 \Delta\Psi=r-b,\quad\partial_n\Psi=0,\quad\bar\Psi=0.
\]
The limit is independent of time.

We next identify the pressure without assuming strong compactness
of the density. Let
\[
 \varphi_j=\Phi_j-\Psi,\qquad
 \Delta\varphi_j=\rho_j-r,\quad\partial_n\varphi_j=0.
\]
Test conservative momentum by $\chi(t)\nabla\varphi_j$, where
$0\le\chi\in C_c^\infty(0,1)$. Its normal trace vanishes.
The time terms tend to zero because $\rho_j u_j\to0$ in $L^2L^2$,
$\nabla\varphi_j\to0$ uniformly, and
$\|\nabla\partial_t\varphi_j\|_{L^{2}}\le C\|\rho_j u_j\|_{L^{2}}$.
The convective term tends to zero by the fixed density moments,
the Neumann estimate for $D^2\varphi_j$, and $u_j\to0$ in $L^2H^1$.
For the bulk viscous term, Cauchy--Schwarz gives
\[
\begin{aligned}
 &\left|\int_0^1\!\!\int_\Om
   \chi(2\mu+\lambda(\rho_j))\diver u_j\,(\rho_j-r)\dd\,\mathrm dt
                                                      \right|\\
 &\quad\le C\left(\int_0^1 A_{3,j}^2\,\mathrm dt\right)^{1/2}
       \left(\int_0^1\!\!\int_\Om
           (1+\rho_j^\ell)|\rho_j-r|^2\dd\,\mathrm dt\right)^{1/2}
       \longrightarrow0.
\end{aligned}
\]
The vorticity interior term vanishes because
$\curl\nabla\varphi_j=0$; its Navier boundary contribution tends to
zero by the trace estimate and $u_j\to0$ in $L^2H^1$.
Finally, the electric term tends to zero since
$\nabla\varphi_j\to0$ uniformly, $\nabla\Phi_j$ is uniformly bounded,
and $\rho_j$ has fixed mass. Consequently
\[
 \int_0^1\!\!\int_\Om
     \chi\rho_j^\gamma(\rho_j-r)\dd\,\mathrm dt\longrightarrow0.
\]
Subtracting the weakly vanishing term containing
$\chi r^\gamma(\rho_j-r)$ and using
\[
 (a^\gamma-c^\gamma)(a-c)\ge|a-c|^{\gamma+1},\qquad a,c\ge0,
\]
proves $\rho_j\to r$ strongly in
$L^{\gamma+1}_{\mathrm{loc}}((0,1)\times\Om)$.
Interpolation with the higher fixed moments gives every additional
finite strong convergence required here.

The inertial and viscous terms vanish in distributions, so the
momentum limit is $\nabla r^\gamma=r\nabla\Psi$.
The mass is preserved. Lemma~\ref{nsp:stationary} therefore gives
$r=\rho_s$ and $\Psi=\Phi_s$.
The strong density convergence, the potential convergence, and the
vanishing kinetic energy imply
\[
 \int_{1/3}^{2/3}\mathcal E_{s,j}(t)\,\mathrm dt\longrightarrow0.
\]
By monotonicity, $\mathcal E_{s,j}(1)$ tends to zero, a contradiction.

Thus for each $\varepsilon>0$ there is $\delta_\varepsilon>0$,
depending only on the common bounds, such that
\[
 \int_a^{a+1}A_3^2\,\mathrm dt\le\delta_\varepsilon
       \quad\Longrightarrow\quad
       \mathcal E_s(a+1)\le\varepsilon.
\]
Choose an integer $N$ with $N\delta_\varepsilon>\mathcal E_s(0)$.
By \eqref{nsp:relative-energy-identity}, at least one of the first
$N$ unit intervals has dissipation at most $\delta_\varepsilon$.
Its endpoint energy is at most $\varepsilon$, and energy
monotonicity then gives the result with $T_\varepsilon=N$.
If the interval of existence ends before $N$, no entry conclusion
is needed on that interval.
\end{proof}

In particular, the moment bounds, the strict positivity of $\rho_s$,
and \eqref{nsp:relative-energy} imply that, for every fixed $p<\infty$,
$\|\rho-\rho_s\|_{L^{p}}$ becomes arbitrarily small after a data-dependent
time. The Neumann estimate then gives
\begin{equation}\label{nsp:small-electric-perturbation}
 \|\nabla(\Phi-\Phi_s)\|_{L^{\infty}}\le\eta
 \qquad(t\ge T_\eta),
\end{equation}
where $T_\eta$ is independent of a pointwise density bound.
For example, relative internal energy controls the density difference
in $L^{\min\{2,\gamma\}}$; interpolation with a fixed higher moment
gives its smallness in a chosen $L^p$, $p>2$.


\section{Relative pressure and uniform estimates}\label{sec:uniform-estimates}
\subsection{Relative pressure dissipation and conditional decay}

Let $(s,\Phi_s)$ be the stationary solution, normalized by
$\overline{\Phi_s}=0$, so that
\[
 \nabla s^\gamma=s\nabla\Phi_s,\qquad
 \Delta\Phi_s=s-b,\qquad \partial_n\Phi_s=0,\qquad
 \int_\Om s\dd=\int_\Om b\dd.
\]
We use $0<c_s\le s\le C_s$ and $s\in W^{2,q}$, $q>2$.
Put $P_s=s^\gamma$, $\psi=\Phi-\Phi_s$, and
$h=2\mu+\lambda(\rho)$. Thus
\[
 \Delta\psi=\rho-s,\qquad \partial_n\psi=0,
 \qquad\bar\psi=0.
\]
In this subsection set
\[
 A_3^2=\int_\Om\bigl(h(\diver u)^2+\mu\omega^2\bigr)\dd
                +\mu\int_{\partial\Om}\alpha u_\tau^2\ds.
\]
The boundary coercivity gives $\|u\|_{H^1}\le CA_3$.
Fix an interval $[0,S]$ of smooth positive existence. All constants below
are independent of $S$, the density maximum, and the positive lower bound
in the approximation.

We use the previously established bounds
\begin{equation}\label{nsp:pressure-inputs}
\begin{aligned}
 &\sup_{0\le t\le S}\|\rho(t)\|_{L^{p}}\le C_p,\qquad
       \int_0^S A_3^2\,\mathrm dt\le C,\\
 &\int_0^S\int_\Om(\rho-L_p)_+^p\dd\,\mathrm dt\le C_p,
       \qquad p\ge1,
\end{aligned}
\end{equation}
where each $L_p$ is fixed and may be increased.
The polynomial upper bound $h\le C(1+\rho^\ell)$ has a fixed $\ell\ge1$.
We also use the uniform entry lemma: there is a fixed $t_e$ such that
\begin{equation}\label{nsp:small-electric-entry}
 \|\nabla\psi(t)\|_{L^{\infty}}\le\varepsilon_s,
                  \qquad t_e\le t\le S,
\end{equation}
where $\varepsilon_s>0$ is sufficiently small in terms of $s$ and
$\gamma$. Both $t_e$ and this assertion are uniform over the
approximations. No density maximum is used in these inputs.

We first prove
\begin{equation}\label{nsp:relative-pressure}
\begin{aligned}
 &\int_0^S\int_\Om
          (\rho+s)^{\gamma-1}(\rho-s)^2\dd\,\mathrm dt\\
 &\qquad+\int_0^S
       \bigl(\|P-P_s\|_{L^{2}}^2+\|P-P_s\|_{L^{4}}^4\bigr)\,\mathrm dt\le C.
\end{aligned}
\end{equation}
For this purpose use the tangent field $v=\nabla\psi/s$ and put
$K=\int_\Om\rho u\cdot v\dd$. Testing conservative momentum by $v$
gives the exact identity
\begin{equation}\label{nsp:weighted-pressure-test}
\begin{aligned}
 K'={}&\int_\Om\frac{(P-P_s)(\rho-s)}s\dd
       +\int_\Om\frac\rho s|\nabla\psi|^2\dd\\
 &-\int_\Om
       \frac{P-P_s-\gamma s^{\gamma-1}(\rho-s)}{s^2}
                       \nabla s\cdot\nabla\psi\dd\\
 &+\int_\Om\rho u\cdot\frac{\nabla\Phi_t}s\dd
       +\int_\Om\rho u\otimes u:\nabla v\dd\\
 &-\int_\Om h\diver u\diver v\dd
       -\mu\int_\Om\omega\curl v\dd
       -\mu\int_{\partial\Om}\alpha u_\tau v_\tau\ds.
\end{aligned}
\end{equation}
Indeed,
\[
 \diver v=\frac{\rho-s}s-\frac{\nabla s\cdot\nabla\psi}{s^2},
 \qquad \nabla P_s=\gamma s^{\gamma-1}\nabla s,
\]
and the non-gradient part of the relative electric force is
$(\rho-s)\nabla\Phi_s+\rho\nabla\psi$.
The last term in \eqref{nsp:weighted-pressure-test} is precisely the
Navier boundary contribution. The curl of $v$ need not vanish.

For every $z\ge0$,
\[
 0\le z^\gamma-1-\gamma(z-1)
       \le C_\gamma(z^\gamma-1)(z-1).
\]
Apply this with $z=\rho/s$. By decreasing $\varepsilon_s$ in
\eqref{nsp:small-electric-entry}, the absolute value of the third
term in \eqref{nsp:weighted-pressure-test} is at most one quarter
of its first term. Moreover,
\[
 \frac{(P-P_s)(\rho-s)}s
       \ge c(\rho+s)^{\gamma-1}(\rho-s)^2
       \ge c|\rho-s|^2.
\]

We next estimate all cross terms without the density maximum.
The weak Neumann equation
$\Delta\Phi_t=-\diver(\rho u)$ gives
$\|\nabla\Phi_t\|_{L^{2}}\le\|\rho u\|_{L^{2}}$.
Consequently, \eqref{nsp:pressure-inputs}, the fixed-exponent Sobolev
inequalities, and Neumann regularity imply
\[
\begin{aligned}
 \left|\int_\Om\rho u\cdot\frac{\nabla\Phi_t}s\dd\right|
       &\le C\|\rho u\|_{L^{2}}^2\le CA_3^2,\\
 \left|\int_\Om\rho u\otimes u:\nabla v\dd\right|
       &\le\|\rho\|_{L^{4}}\|u\|_{L^{8}}^2\|\nabla v\|_{L^{2}}\le CA_3^2,\\
 \|v\|_{H^1}+\|\curl v\|_{L^{2}}&\le C\|\rho-s\|_{L^{2}}.
\end{aligned}
\]
Here $\|\rho-s\|_{L^{2}}$ is uniformly bounded; it is not assumed small.
The vorticity and boundary terms therefore contribute at most
$\varepsilon\|\rho-s\|_{L^{2}}^2+C_\varepsilon A_3^2$.

For the bulk-viscosity term write
\[
 \diver v=\frac{\rho-s}s-B,\qquad
 B=\frac{\nabla s\cdot\nabla\psi}{s^2}.
\]
The uniform viscosity moment and the fixed Neumann estimate give
\[
 \int_\Om hB^2\dd
       \le\|h\|_{L^{2}}\|B\|_{L^{4}}^2\le C\|\rho-s\|_{L^{2}}^2,
\]
so the term with $B$ is bounded by
$\varepsilon\|\rho-s\|_{L^{2}}^2+C_\varepsilon A_3^2$.
Choose a fixed $L\ge L_{\ell+2}$ with $L\ge2\|s\|_{L^{\infty}}+1$.
On $\{\rho\le2L\}$, $h$ is bounded, while on
$\{\rho>2L\}$ one has $h|\rho-s|^2\le C\rho^{\ell+2}$.
Weighted Cauchy--Schwarz and Young's inequality yield
\[
 \left|\int_\Om h\diver u\frac{\rho-s}s\dd\right|
 \le\varepsilon\|\rho-s\|_{L^{2}}^2+C_\varepsilon A_3^2
             +C\int_\Om(\rho-L)_+^{\ell+2}\dd.
\]
Substituting these estimates into
\eqref{nsp:weighted-pressure-test} and choosing $\varepsilon$ small gives
\begin{equation}\label{nsp:weighted-pressure-lower}
 K'\ge c\int_\Om(\rho+s)^{\gamma-1}(\rho-s)^2\dd
                 -CA_3^2-C\int_\Om(\rho-L)_+^{\ell+2}\dd,
                         \qquad t\ge t_e.
\end{equation}
The basic kinetic energy, \eqref{nsp:pressure-inputs}, and Neumann
regularity also give
\[
 |K|\le C\|\sqrt\rho u\|_{L^{2}}
               \|\rho\|_{L^{2}}^{1/2}\|\nabla\psi\|_{L^{4}}\le C.
\]
If $S>t_e$, integrate \eqref{nsp:weighted-pressure-lower} on $[t_e,S]$.
A fixed density moment controls $[0,\min\{S,t_e\}]$.
This proves the first term of \eqref{nsp:relative-pressure}, including
the case $S\le t_e$, with a constant independent of $S$.
For $p=2,4$, increasing a fixed $L_{p\gamma}$ if necessary gives
\[
 |P-P_s|^p\le C_p|\rho-s|^2
                  +C_p(\rho-L_{p\gamma})_+^{p\gamma}.
\]
This proves all of \eqref{nsp:relative-pressure}.

We now establish the conditional decay needed for the flux
interpolation. Set
\[
 R=\max\{2,2\|s\|_{L^{\infty}},
                     \sup_{0\le t\le S}\|\rho(t)\|_{L^{\infty}}\}.
\]
There is a fixed $N\ge1$ such that
\begin{equation}\label{nsp:conditional-relative-decay}
 \|\rho(t)-s\|_{L^{2}}^2\le CR^N\exp(-cR^{-N}t),
                          \qquad 0\le t\le S,
\end{equation}
and, uniformly for $m\ge8$,
\begin{equation}\label{nsp:relative-pressure-high}
 \int_0^S\|P-P_s\|_{L^{m}}^2\,\mathrm dt\le CmR^N.
\end{equation}
The constants and $N$ are independent of $m$, $S$, and $R$.

Use the relative energy $\mathcal E_s$ from
\eqref{nsp:relative-energy} and its dissipation identity
\eqref{nsp:relative-energy-identity}.
For a fixed sufficiently large $N_0>\ell+\gamma+4$, the relative internal
energy, the uniform $L^2$ density moment, and Neumann regularity give
\[
\begin{aligned}
 cR^{-N_0}\|\rho-s\|_{L^{2}}^2&\le\mathcal E_s
             \le CR^{N_0}(A_3^2+\|\rho-s\|_{L^{2}}^2),\\
 |K|&\le C\|\sqrt\rho u\|_{L^{2}}\|\rho-s\|_{L^{2}}
             \le CR^{N_0}\mathcal E_s.
\end{aligned}
\]
For the second inequality in the last line, the bound preceding it
uses $\|\rho\|_{L^{2}}\le C$ and
$\|\nabla\psi\|_{L^{4}}\le C\|\rho-s\|_{L^{2}}$.
When $\rho\le R$, estimate the bulk-viscosity cross term directly,
without the tail split:
\[
 \left|\int_\Om h\diver u\frac{\rho-s}s\dd\right|
       \le\varepsilon\|\rho-s\|_{L^{2}}^2
                      +C_\varepsilon R^{\ell}A_3^2.
\]
The preceding estimates in the weighted identity thus imply
\[
 K'\ge c\|\rho-s\|_{L^{2}}^2-CR^{\ell}A_3^2,
                                  \qquad t_e\le t\le S.
\]
Choose $\delta=c_0R^{-N_0}$, with a fixed $c_0>0$ small enough that
$C\delta R^\ell\le1/2$ and $\delta|K|\le\mathcal E_s/2$. Then
\[
 \mathcal E_s-\delta K\asymp\mathcal E_s,\qquad
 (\mathcal E_s-\delta K)'
       +c\bigl(A_3^2+\delta\|\rho-s\|_{L^2}^2\bigr)\le0.
\]
Since $A_3^2+\delta\|\rho-s\|_{L^2}^2
\ge cR^{-2N_0}(\mathcal E_s-\delta K)$, we obtain
\[
 (\mathcal E_s-\delta K)'
       +cR^{-2N_0}(\mathcal E_s-\delta K)\le0
                              \quad\hbox{on }[t_e,S].
\]
The bound $\mathcal E_s(t_e)\le\mathcal E_s(0)\le C$ and the lower bound
for $\mathcal E_s$ prove \eqref{nsp:conditional-relative-decay} with
$N=2N_0$ at this stage. The same estimate on $[0,t_e]$ follows from the
uniform $L^2$ moment after increasing the data-dependent prefactor,
since $\exp(cR^{-2N_0}t_e)\le\exp(ct_e)$. This also covers $S<t_e$.

Finally, the mean-value theorem and interpolation give, for $m\ge8$,
\[
\begin{aligned}
 \|P-P_s\|_{L^m}^2
 &\le CR^{2\gamma-2}\|\rho-s\|_{L^m}^2
       \le CR^{2\gamma}\|\rho-s\|_{L^2}^{4/m}\\
 &\le CR^{2\gamma+N_0/2}
                         \exp(-2cR^{-2N_0}t/m).
\end{aligned}
\]
The constant is uniform in $m\ge8$: the prefactor in the conditional
decay estimate is raised to the power $2/m$, and $4N_0/m\le N_0/2$.
Integration therefore gives
\[
 \int_0^S\|P-P_s\|_{L^m}^2\,\mathrm dt
        \le CmR^{2\gamma+5N_0/2}.
\]
Fix $N\ge4N_0+2\gamma+1$. Both
\eqref{nsp:conditional-relative-decay} and
\eqref{nsp:relative-pressure-high} then hold with this same exponent,
independently of $m$, $S$, and $R$.

\paragraph{Electric terms in the subsequent estimates.}
The relative flux is
$G=h\diver u-(P-P_s)$. Its energy identity has the additional terms
\[
 2\int_\Om\frac Gh u\cdot\nabla P_s\dd
       +2\int_\Om\bigl((\rho-s)\nabla\Phi_s
                              +\rho\nabla\psi\bigr)\cdot\dot u\dd.
\]
The static part is an exact time derivative up to velocity-quadratic
terms:
\begin{equation}\label{nsp:static-electric-cross}
\begin{aligned}
 \int_\Om(\rho-s)\nabla\Phi_s\cdot\dot u\dd
 ={}&\frac{\mathrm d}{\mathrm dt}
       \int_\Om(\rho-s)\nabla\Phi_s\cdot u\dd\\
 &-\int_\Om\rho u^{\mathsf T}\nabla^2\Phi_s u\dd
       -\int_\Om\nabla P_s\cdot(u\cdot\nabla u)\dd.
\end{aligned}
\end{equation}
The last two integrals are bounded by $CA_3^2$, while
\[
 \left|\int_\Om\rho\nabla\psi\cdot\dot u\dd\right|
       \le\varepsilon\|\sqrt\rho\dot u\|_{L^{2}}^2
                          +C_\varepsilon\|P-P_s\|_{L^{2}}^2.
\]
One therefore subtracts twice the cross integral in
\eqref{nsp:static-electric-cross} from the relative-flux energy,
adding a fixed sufficiently large constant to preserve positivity.
\subsection{Relative logarithmic energy and acceleration}
\label{nsp:uni-energy-section}

Write $s=\rho_s$, $P_s=s^\gamma$, $q=P-P_s$, and
$\phi=\Phi-\Phi_s$.  Thus $\Delta\phi=\rho-s$ and
$\partial_n\phi=0$.  The uniform moments, relative pressure estimate
\eqref{nsp:relative-pressure}, and boundary coercivity give
\begin{equation}\label{nsp:uni-budget}
 D=A_3^2+\|q\|_{L^{2}}^2+\|q\|_{L^{4}}^4,\qquad
 \int_0^S D\,\mathrm dt\le C,\qquad
 \sup_{t\le S}\bigl(\|q\|_{L^{2}}+\|q\|_{L^{4}}\bigr)\le C,\qquad
 \|u\|_{H^1}\le CA_3.
\end{equation}
All constants in this subsection are independent of $S$, the density
maximum, and the positive lower bound in a smooth approximation.
Fix
\[
 R=\max\{2,2\|s\|_{L^{\infty}},\sup_{0\le t\le S}\|\rho(t)\|_{L^{\infty}}\},
 \qquad h=2\mu+\lambda(\rho),\qquad
 y=\int_\Om\rho|\dot u|^2\dd,
\]
and set
\begin{equation}\label{nsp:uni-relative-flux-energy}
 G=h\diver u-q,\qquad
 A_1^2=\int_\Om\left(\mu\omega^2+\frac{G^2}h\right)\dd
       +\mu\int_{\partial\Om}\alpha u_\tau^2\ds.
\end{equation}
The corrected energy used below is
\begin{equation}\label{nsp:uni-corrected-energy}
 C_s(t)=\int_\Om(\rho-s)\nabla\Phi_s\cdot u\dd,
 \qquad \mathscr{E}=e_0+A_1^2-2C_s(t).
\end{equation}
Here $e_0$ is a fixed sufficiently large constant.  Indeed,
\begin{equation}\label{nsp:uni-energy-comparison}
\begin{gathered}
 A_1^2\le C(A_3^2+\|q\|_{L^{2}}^2)\le CD,\qquad
 A_3^2\le C(A_1^2+\|q\|_{L^{2}}^2),\qquad |C_s|\le CA_3,\\
 \mathscr{E}\asymp1+A_1^2,\qquad \mathscr{E}\ge1,\qquad
 A_3^2\le C\mathscr{E},\qquad \mathscr{E}\le C(1+D).
\end{gathered}
\end{equation}
For $|C_s|\le CA_3$, use the uniform $L^2$ density moment and the
bounded fixed field $\nabla\Phi_s$.  The other comparisons follow
from $h\diver u=G+q$ and $h\ge2\mu$.

The relative momentum equation is
\begin{equation}\label{nsp:uni-relative-momentum}
 \rho\dot u=\nabla G+\mu\nabla^\perp\omega
                  +(\rho-s)\nabla\Phi_s+\rho\nabla\phi.
\end{equation}
Since $|\rho-s|\le C|q|$, fixed-exponent Neumann regularity and the
uniform density moments imply
\[
 \|(\rho-s)\nabla\Phi_s+\rho\nabla\phi\|_{L^{2}}\le C\|q\|_{L^{2}}.
\]
The flux--vorticity estimate \eqref{mom:fluxelliptic} therefore gives
\begin{equation}\label{nsp:uni-low-gradient}
\begin{aligned}
 \|\nabla G\|_{L^{2}}+\|\nabla\omega\|_{L^{2}}
   &\le C(\sqrt{Ry}+A_3+\|q\|_{L^{2}}),\\
 |\bar G|&\le CA_1,\qquad
 \|G\|_{L^{2}}\le CR^{\ell/2}A_1.
\end{aligned}
\end{equation}
The uniform viscosity moment gives $\int_\Om h\dd\le C$. Hence
weighted Cauchy--Schwarz directly yields
\[
 |\bar G|\le\frac1{|\Om|}
       \left(\int_\Om\frac{G^2}h\dd\right)^{1/2}
       \left(\int_\Om h\dd\right)^{1/2}\le CA_1.
\]
This does not require $\bar q=0$; in general,
$\bar G=\overline{\lambda\diver u}-\bar q$.

We claim that
\begin{equation}\label{nsp:uni-log-differential}
 \mathscr{E}'+cy\le CR\log(e+R)\mathscr{E}D,
\end{equation}
and consequently
\begin{equation}\label{nsp:uni-log-energy}
 \sup_{t\le S}\log \mathscr{E}(t)+\int_0^S\frac{y}{\mathscr{E}}\,\mathrm dt
          \le CR\log(e+R).
\end{equation}
To obtain the exact identity, write $d_u=\diver u$ and use
\[
\begin{aligned}
 D_tG&=h\diver\dot u-\rho\lambda'd_u^2
       -h\operatorname{tr}((\nabla u)^2)+\gamma P d_u
                                      +u\cdot\nabla P_s,\\
 D_t\omega&=\curl\dot u-\omega d_u,\qquad
 \diver\dot u=D_td_u+d_u^2-2\det\nabla u.
\end{aligned}
\]
On the boundary $\dot u\cdot n=-\kappa u_\tau^2$ and
$\dot u\cdot\tau=\partial_tu_\tau+u_\tau\partial_\tau u_\tau$.
Thus the tangential boundary term satisfies
\[
 2\mu\int_{\partial\Om}\omega\dot u\cdot\tau\ds
 =-\mu\frac{\mathrm d}{\mathrm dt}
           \int_{\partial\Om}\alpha u_\tau^2\ds
       +\frac{2\mu}3\int_{\partial\Om}
                              (\partial_\tau\alpha)u_\tau^3\ds.
\]
Testing \eqref{nsp:uni-relative-momentum} by $2\dot u$, substituting
these identities, and using
$\operatorname{tr}((\nabla u)^2)=d_u^2-2\det\nabla u$ gives
\begin{equation}\label{nsp:uni-exact-energy}
\begin{aligned}
 \mathscr{E}'+2y={}&
 \int_\Om\left(1+\frac{\rho\lambda'}h\right)
                       \frac{G^2}h\diver u\dd
 -2\int_\Om\left(1+\frac{\rho\lambda'}h\right)
                       G(\diver u)^2\dd\\
 &+2\gamma\int_\Om\frac{GP}h\diver u\dd
       +4\int_\Om G\det\nabla u\dd
       -\mu\int_\Om\omega^2\diver u\dd\\
 &-2\int_{\partial\Om}\kappa G u_\tau^2\ds
       +\frac{2\mu}3\int_{\partial\Om}
                    (\partial_\tau\alpha)u_\tau^3\ds\\
 &+2\int_\Om\frac Gh u\cdot\nabla P_s\dd
       +2\int_\Om\rho\nabla\phi\cdot\dot u\dd\\
 &-2\int_\Om\rho u^{\mathsf T}\nabla^2\Phi_s u\dd
       -2\int_\Om\nabla P_s\cdot(u\cdot\nabla u)\dd.
\end{aligned}
\end{equation}
The static electric contribution has been removed using
\eqref{nsp:static-electric-cross}.  In particular,
\eqref{nsp:uni-log-differential} is an inequality for the corrected
$\mathscr{E}$, not for the derivative of $A_1^2$ alone.

We estimate the terms in this identity. Put $d=A_3+\|q\|_{L^{2}}$.
The uniform pressure moments and \eqref{nsp:uni-energy-comparison}
imply
\begin{equation}\label{nsp:uni-small-energy-factors}
 A_1^2+A_3^2+d^2\le CD,\qquad
 D\le C\mathscr E,\qquad \mathscr E\ge1.
\end{equation}
For $0<\delta\le1/4$, weighted H\"older gives
\[
 \left\|\frac{G^2}h\right\|_{L^{2}}
 \le C A_1^{1-\delta}
           \|G\|_{L^{2(1+\delta)/\delta}}^{1+\delta}.
\]
Apply \eqref{eng:GN} to $G-\bar G$. Its relevant power is
\[
 \|G-\bar G\|_{L^{2(1+\delta)/\delta}}^{1+\delta}
 \le\frac C{\sqrt\delta}
          \|G-\bar G\|_{L^{2}}^\delta\|\nabla G\|_{L^{2}}.
\]
The constant is uniform in $\delta$, since
$\delta^{-(1+\delta)/2}\le C\delta^{-1/2}$.
Using \eqref{nsp:uni-low-gradient} and $|\bar G|\le CA_1$ yields
\begin{equation}\label{nsp:uni-weighted-quadratic}
 \left\|\frac{G^2}h\right\|_{L^{2}}
 \le\frac C{\sqrt\delta}R^{\ell\delta/2}
                      A_1(\sqrt{Ry}+d)
                         +CA_1^2.
\end{equation}
Since $|\rho\lambda'|/h\le C$ and $\diver u=(G+q)/h$,
the absolute sum of the first three integrals in
\eqref{nsp:uni-exact-energy} is at most
\[
 CA_3\left(\left\|\frac{G^2}h\right\|_{L^{2}}
                +\left\|\frac{G(P+P_s)}h\right\|_{L^{2}}\right).
\]
The uniform $L^4$ pressure moment and $h\ge2\mu$ give
\[
 \left\|\frac{G(P+P_s)}h\right\|_{L^{2}}
       \le C\left\|\frac{G^2}h\right\|_{L^{2}}^{1/2}.
\]
For example, the part containing $y$ in this square root is bounded by
\[
 C\delta^{-1/4}R^{(1+\ell\delta)/4}
                   A_3A_1^{1/2}y^{1/4}
 \le\varepsilon y+
 C_\varepsilon\delta^{-1/3}R^{(1+\ell\delta)/3}
                   A_3^{4/3}A_1^{2/3}.
\]
The last energy factor is at most $CD$.
In the quadratic expression itself, the $y$ term gives
$\varepsilon y+C_\varepsilon\delta^{-1}R^{1+\ell\delta}A_3^2A_1^2$.
The mean term contributes $CA_3A_1^2\le C\mathscr ED$, and its
square root contributes $CA_3A_1\le CD$.
The remaining terms are bounded by
$C(1+\delta^{-1/2}R^{\ell\delta/2})\mathscr ED$ using
\eqref{nsp:uni-small-energy-factors}; their square-root contributions
are bounded by the same quantity. Choose
$\delta=[4\log(e+R)]^{-1}$, so $R^{\ell\delta}\le e^{\ell/4}$.
We have proved
\begin{equation}\label{nsp:uni-coefficient-bound}
 \text{absolute sum of the first three terms}
       \le\varepsilon y+C_\varepsilon R\log(e+R)\mathscr ED.
\end{equation}

For the Jacobian, let $p\ge8$. Equations \eqref{eng:GN} and
\eqref{nsp:uni-low-gradient} give
\[
 \|G\|_{L^{p}}+\|\omega\|_{L^{p}}
 \le C\sqrt p\,R^{\ell/p}\mathscr E^{1/p}
                              (\sqrt{Ry}+d)^{1-2/p}.
\]
Here the mean contribution is included because $d^2\le C\mathscr E$.
The quantitative Hodge estimate \eqref{mom:Hodge},
$\|q\|_{L^{p}}\le CR^\gamma$, and $\|u\|_{L^{2}}\le CA_3$ imply
\[
 \|\nabla u\|_{L^{p}}
 \le Cp^{13/2}R^{\ell/p}\mathscr E^{1/p}
                  (\sqrt{Ry}+d)^{1-2/p}+Cp^6R^\gamma.
\]
Interpolate this bound with $\|\nabla u\|_{L^{2}}\le CA_3$ at
$2p/(p-1)$, whose interpolation weight on $L^p$ is $1/(p-2)$.
H\"older's inequality gives
\begin{equation}\label{nsp:uni-Jacobian-interpolation}
\begin{aligned}
 \left|\int_\Om G\det\nabla u\dd\right|
 \le C\sqrt p\,A_3^{2(p-3)/(p-2)}
 \bigl[&R^{\ell/(p-2)}\mathscr E^{1/(p-2)}(\sqrt{Ry}+d)\\
 &+R^{\ell/p+2\gamma/(p-2)}\mathscr E^{1/p}
                                  (\sqrt{Ry}+d)^{1-2/p}\bigr].
\end{aligned}
\end{equation}
The factors $p^{13/(p-2)}$ and $p^{12/(p-2)}$ have been absorbed
into $C$. Choose $p=8+(\ell+2\gamma+2)\log(e+R)$; both displayed
extra powers of $R$ are then bounded. Young's inequality absorbs
the contributions of $\sqrt y$. Their energy factors satisfy
\[
\begin{aligned}
 \mathscr E^{2/(p-2)}A_3^{4(p-3)/(p-2)}
       &\le C\mathscr EA_3^2,\\
 \mathscr E^{2/(p+2)}
          A_3^{4p(p-3)/((p-2)(p+2))}
       &\le C\mathscr EA_3^2,
\end{aligned}
\]
by $A_3^2\le C\mathscr E$ and $\mathscr E\ge1$.
The contributions containing $d$ are bounded by
$C\sqrt p\,\mathscr ED$, using
\eqref{nsp:uni-small-energy-factors}. Hence
\begin{equation}\label{nsp:uni-Jacobian-bound}
 \left|\int_\Om G\det\nabla u\dd\right|
       \le\varepsilon y+C_\varepsilon R\log(e+R)\mathscr ED.
\end{equation}
For the vorticity term, the fixed-exponent interpolation inequality
gives
\[
 \left|\int_\Om\omega^2\diver u\dd\right|
 \le C\|\omega\|_{L^{2}}(\|\nabla\omega\|_{L^{2}}+\|\omega\|_{L^{2}})
                                  \|\diver u\|_{L^{2}}
 \le CA_3^2(\sqrt{Ry}+d)
 \le\varepsilon y+C_\varepsilon R\mathscr ED.
\]
The trace embedding $H^1(\Om)\to L^3(\partial\Om)$ and Poincar\'e's
inequality give
\[
\begin{aligned}
 \left|\int_{\partial\Om}\kappa G u_\tau^2\ds\right|
       &\le C\|G\|_{H^1}\|u\|_{H^1}^2
        \le CA_3^2(\sqrt{Ry}+d)
        \le\varepsilon y+C_\varepsilon R\mathscr ED,\\
 \left|\int_{\partial\Om}(\partial_\tau\alpha)u_\tau^3\ds\right|
       &\le C\|u\|_{H^1}^3\le C\mathscr ED.
\end{aligned}
\]
The new terms in the last two lines of
\eqref{nsp:uni-exact-energy} satisfy
\[
\begin{aligned}
 \left|\int_\Om\frac Gh u\cdot\nabla P_s\dd\right|
       &\le CA_1A_3\le CD,\\
 \left|\int_\Om\rho\nabla\phi\cdot\dot u\dd\right|
       &\le\varepsilon y+C_\varepsilon\|q\|_{L^{2}}^2,\\
 \left|\int_\Om\rho u^{\mathsf T}\nabla^2\Phi_s u\dd\right|
 +\left|\int_\Om\nabla P_s\cdot(u\cdot\nabla u)\dd\right|
       &\le CA_3^2.
\end{aligned}
\]
This proves \eqref{nsp:uni-log-differential}; division by $\mathscr{E}$ and
\eqref{nsp:uni-budget} prove \eqref{nsp:uni-log-energy}.

We next prove the positive-time acceleration estimate:
for some fixed $N\ge1$,
\begin{equation}\label{nsp:uni-log-acceleration}
 \sup_{\tau\le t\le S}\log(e+y(t))
 +\int_\tau^S\frac{\|\nabla\dot u\|_{L^{2}}^2}{e+y}\,\mathrm dt
           \le C_\tau R^N,\qquad 0<\tau<S.
\end{equation}
First, testing \eqref{nsp:uni-relative-momentum} by $u$ and using
the basic kinetic energy gives
\[
 A_3^2\le CA_3\|q\|_{L^{2}}+C\sqrt y
       +\|(\rho-s)\nabla\Phi_s+\rho\nabla\phi\|_{L^{2}}\|u\|_{L^{2}}
       \le CA_3+C\sqrt y.
\]
Thus
\begin{equation}\label{nsp:uni-energy-square}
 \mathscr{E}^2\le C(1+y),\qquad
 \|\nabla u\|_{L^{4}}^4\le CR^ND(1+y).
\end{equation}
For the second estimate use the fixed-exponent Hodge and
Gagliardo--Nirenberg inequalities with
\eqref{nsp:uni-low-gradient}; the terms are bounded by
$CR^\ell A_1^2(Ry+d^2)+CA_3^4+C\|q\|_{L^{4}}^4+Cd^4$, which is at
most the displayed right-hand side.

Use the fixed normal extension and correction $Q_n$ from
\eqref{acc:qbounds}, and set
\[
 v=\dot u+Q_n,\qquad z=\int_\Om\rho|v|^2\dd,
\]
so that $v\cdot n=0$. Write
\[
 D_v=\int_\Om(h(\diver v)^2+\mu(\curl v)^2)\dd
                    +\mu\int_{\partial\Om}\alpha v_\tau^2\ds,
 \qquad H=1+z+\mathscr E^2.
\]
The bounds \eqref{acc:qbounds} and \eqref{nsp:uni-energy-square} give
\begin{equation}\label{nsp:uni-correction-bounds}
\begin{gathered}
 \int_\Om\rho|Q_n|^2\dd\le CA_3^4,\qquad
 \|Q_n\|_{H^1}^2
       \le C(A_3^4+A_3^2\|\nabla u\|_{L^{4}}^2)
       \le CR^ND(1+y),\\
 \|v\|_{H^1}^2\le CD_v,\qquad
 c(1+y)\le H\le C(1+y).
\end{gathered}
\end{equation}
The uniform $L^2$ density moment and boundary coercivity give
\[
 \int_\Om\rho|Q_n|^2\dd
 \le C\|\rho\|_{L^{2}}\|u\|_{L^{8}}^4\le CA_3^4.
\]
For the $H^1$ bound use
$A_3^2\|\nabla u\|_{L^{4}}^2\le(A_3^4+\|\nabla u\|_{L^{4}}^4)/2$.
The estimate for $v$ is \eqref{uni:coercive}. Moreover,
$y\le2z+CA_3^4$ and $z\le2y+CA_3^4$; these, together with
$A_3^4\le C\mathscr E^2$ and
\eqref{nsp:uni-energy-square}, prove the last comparison.

Use the stress and commutation identities
\eqref{acc:commutation}--\eqref{acc:stressremainder}, with the physical
flux $F=G-P_s$. Thus
\[
 \rho D_tv=\diver\Sigma+\rho D_tQ_n+\rho D_t\nabla\Phi.
\]
We estimate these terms using the uniform budget $D$.
After testing by $v$, the principal interior terms give
the interior part of $D_v$. Replacing $\dot u$ by $v-Q_n$ costs at
most $\varepsilon D_v+CR^\ell\|\nabla Q_n\|_{L^{2}}^2$.
The stress remainder satisfies
\[
\begin{aligned}
 \left|\int_\Om\mathcal Q:\nabla v\dd\right|
 &\le\varepsilon D_v+CR^{2\ell}\|\nabla u\|_{L^{4}}^4
                         +C\int_\Om P^2|\nabla u|^2\dd\\
 &\le\varepsilon D_v+CR^NDH.
\end{aligned}
\]
For the last step use $P=q+P_s$ and
\[
 \int_\Om P^2|\nabla u|^2\dd
 \le C\bigl(A_3^2+\|q\|_{L^{4}}^4+\|\nabla u\|_{L^{4}}^4\bigr).
\]
The uniform $L^2$ density moment, H\"older's inequality, and boundary
coercivity give
\[
\begin{aligned}
 \int_\Om\rho|D_tQ_n||v|\dd
 &\le C(\sqrt y\,A_3+A_3^3)D_v^{1/2}\\
 &\le\varepsilon D_v+C_\varepsilon(A_3^2y+A_3^6)
       \le\varepsilon D_v+C_\varepsilon DH.
\end{aligned}
\]
The fixed moments and the Neumann equation for $\Phi_t$ imply
\[
 \|D_t\nabla\Phi\|_{L^{2}}
 \le C\|\rho u\|_{L^{2}}+\|u\|_{L^{4}}\|D^2\Phi\|_{L^{4}}\le CA_3.
\]
Triple H\"older and boundary coercivity give the required estimate
without a density-maximum factor:
\begin{equation}\label{nsp:uni-electric-acceleration}
 \left|\int_\Om\rho D_t\nabla\Phi\cdot v\dd\right|
 \le\|\rho\|_{L^{4}}\|D_t\nabla\Phi\|_{L^{2}}\|v\|_{L^{4}}
 \le\varepsilon D_v+C_\varepsilon A_3^2.
\end{equation}
For the boundary terms, apply \eqref{acc:boundarytraction} with
$F=G-P_s$. Its first term is the boundary part of $D_v$.
The trace inequalities and \eqref{nsp:uni-low-gradient} give
\[
\begin{aligned}
 &\left|\int_{\partial\Om}
       [\kappa(G-P_s)u_\tau
                 -\mu(\partial_\tau\alpha)u_\tau^2]v_\tau\ds\right|\\
 &\quad\le C\bigl[(\sqrt{Ry}+d+1)A_3+A_3^2\bigr]D_v^{1/2}
       \le\varepsilon D_v+CR^NDH.
\end{aligned}
\]
For the remaining term extend $\alpha,\tau$ smoothly to the domain.
Stokes' formula for a scalar tangential derivative gives
\[
\begin{aligned}
 \left|\int_{\partial\Om}
                \alpha u_\tau(\partial_\tau u_\tau)v_\tau\ds\right|
 &=\frac12\left|\int_{\partial\Om}
                         \alpha v_\tau\partial_\tau(u_\tau^2)\ds\right|\\
 &\le C\int_\Om
              |\nabla(\alpha v_\tau)|\,|\nabla(u_\tau^2)|\dd\\
 &\le CD_v^{1/2}(A_3^2+A_3\|\nabla u\|_{L^{4}})\\
 &\le\varepsilon D_v
              +C(A_3^4+A_3^2\|\nabla u\|_{L^{4}}^2)
       \le\varepsilon D_v+CR^NDH.
\end{aligned}
\]
All boundary components are included in Stokes' formula with their
induced orientations. Combining the interior and boundary estimates
and absorbing $\varepsilon D_v$ proves
\begin{equation}\label{nsp:uni-z-differential}
 z'+cD_v\le CR^NDH.
\end{equation}
Multiplying \eqref{nsp:uni-log-differential} by $2\mathscr{E}$ and adding
gives, after enlarging $N$,
\begin{equation}\label{nsp:uni-H-differential}
 H'+cD_v\le CR^NDH.
\end{equation}
The logarithmic estimate supplies a point
$s_\tau\in[\tau/2,\tau]$ with $\log H(s_\tau)\le C_\tau R^N$.
Indeed \eqref{nsp:uni-log-energy} bounds
$\int_0^S y\le CR\log(e+R)\exp(CR\log(e+R))$; choose
$s_\tau$ by averaging $y$ over $[\tau/2,\tau]$ and use
$H\le C(1+y)$. No initial acceleration norm is required.
Divide \eqref{nsp:uni-H-differential} by $H$, integrate from this
point, and use \eqref{nsp:uni-budget}.  Finally,
\[
 \frac{\|\nabla\dot u\|_{L^{2}}^2}{e+y}
       \le C\frac{D_v}H+CR^ND
\]
by \eqref{nsp:uni-correction-bounds}.  This proves
\eqref{nsp:uni-log-acceleration}.

\subsection{The uniform flux integral and density maximum}
\label{nsp:uni-density-section}

We prove the two estimates
\begin{equation}\label{nsp:uni-weighted-flux}
 \int_\tau^S\frac{\|G\|_{L^{\infty}}^2}{\mathscr{E}}\,\mathrm dt
          \le C_\tau R^2[\log(e+R)]^2,
\end{equation}
and, on every prefix interval and $\tau\le s_0<t\le S$,
\begin{equation}\label{nsp:uni-local-flux}
 \int_{s_0}^t\|G\|_{L^{\infty}}\,\mathrm d\xi
 \le C_\tau R_t\log(e+R_t)
       \left(t-s_0+\int_{s_0}^tD\,\mathrm d\xi\right)^{1/2}.
\end{equation}
The fixed-exponent bound
$\|\nabla\phi\|_{L^{\infty}}\le C\|\rho-s\|_{L^{4}}\le C\|q\|_{L^{m}}$ for
$m\ge8$ yields, uniformly in $m$,
\[
 \|(\rho-s)\nabla\Phi_s+\rho\nabla\phi\|_{L^{m}}\le CR\|q\|_{L^{m}}.
\]
The weighted Poincar\'e and Sobolev inequalities give
\[
 \|\dot u\|_{L^{m}}\le C\sqrt m
                        (\sqrt y+\|\nabla\dot u\|_{L^{2}}),\qquad
 \|u\|_{L^{m}}\le C\sqrt m A_3.
\]
By \eqref{nsp:uni-low-gradient} and Poincar\'e's inequality,
\[
 \|G\|_{L^{m}}+\|\omega\|_{L^{m}}
       \le C\sqrt m(\sqrt{Ry}+A_3+\|q\|_{L^{2}}).
\]
Apply \eqref{mom:Hodge} to $\diver u=(G+q)/h$, then apply
\eqref{mom:fluxelliptic} to \eqref{nsp:uni-relative-momentum}.
The resulting estimates are
\[
\begin{aligned}
 \|\nabla u\|_{L^{m}}
 &\le Cm^{13/2}(\sqrt{Ry}+A_3+\|q\|_{L^{m}}),\\
 \|\nabla G\|_{L^{m}}
 &\le Cm\bigl(\|\rho\dot u\|_{L^{m}}+\|u\|_{W^{1,m}}
               +\|(\rho-s)\nabla\Phi_s+\rho\nabla\phi\|_{L^{m}}\bigr).
\end{aligned}
\]
Here $\|q\|_{L^{2}}\le C\|q\|_{L^{m}}$ uniformly for $m\ge8$.
Since $\|\rho\dot u\|_{L^{m}}\le R\|\dot u\|_{L^{m}}$, these imply
\begin{equation}\label{nsp:uni-high-gradient}
 \|\nabla G\|_{L^{m}}\le Cm^{15/2}R
       (\sqrt y+\|\nabla\dot u\|_{L^{2}}+A_3+\|q\|_{L^{m}}),\qquad m\ge8.
\end{equation}
The pressure norm is retained, since
\eqref{nsp:relative-pressure-high} controls its integral with only
a polynomial density-maximum cost.

For this interpolation put
\[
 a_0=\frac{A_1}{\sqrt{\mathscr{E}}},\qquad
 L_0=\frac{\sqrt y+A_3+\|q\|_{L^{2}}}{\sqrt{\mathscr{E}}},\qquad
 Z=\frac{\|\nabla\dot u\|_{L^{2}}}{\sqrt{e+y}},\qquad
 \Pi_m=\frac{\|q\|_{L^{m}}}{\sqrt{\mathscr{E}}}.
\]
The estimates already proved imply, after fixing a sufficiently
large $N$,
\begin{equation}\label{nsp:uni-time-factors}
\begin{aligned}
 &a_0\le C,\qquad \int_0^S a_0^2\,\mathrm dt\le C,\qquad
       \int_0^S L_0^2\,\mathrm dt\le CR\log(e+R),\\
 &\int_\tau^S Z^2\,\mathrm dt\le C_\tau R^N,\qquad
       \int_0^S\Pi_m^2\,\mathrm dt\le CmR^N,\qquad
       \sup_{[\tau,S]}\log(e+L_0^2)\le C_\tau R^N.
\end{aligned}
\end{equation}
Choose
\[
 p=8+(N+\ell+4)\log(e+R),\qquad m=e^p,
 \qquad
 \eta=\frac{2(m-2)}{p(m-2)+2m},\quad
 \sigma=\frac{2m}{p(m-2)+2m}.
\]
Apply the two-dimensional interpolation inequality
\eqref{cl:interpolation} to $G-\bar G$.
Equations \eqref{nsp:uni-low-gradient} and
\eqref{nsp:uni-high-gradient} yield
\begin{align*}
 \frac{\|G-\bar G\|_{L^{\infty}}}{\sqrt{\mathscr{E}}}
 &\le C\sqrt p\,
 R^{1/2+(\ell-1)\eta/2+\sigma/2}m^{15\sigma/2}
 a_0^\eta L_0^{1-\eta-\sigma}
       (L_0+\mathscr{E}^{-1/2}Z+L_0Z+\Pi_m)^\sigma\\
 &\le C\sqrt{pR}\,a_0^\eta L_0^{1-\eta-\sigma}
       (L_0^\sigma+Z^\sigma+L_0^\sigma Z^\sigma+\Pi_m^\sigma).
\end{align*}
All the extra powers of $R$ and $m$ are bounded. To spell out the
time integration, let $B=(CR\log(e+R))^{1/2}$, enlarged so $B\ge1$.
The $L^2(\tau,S)$ norms of the first, second, and fourth summands
inside the last product, without its factor $C\sqrt{pR}$, are at most
\[
\begin{aligned}
 &\|a_0\|_{L^{2}}^\eta\|L_0\|_{L^{2}}^{1-\eta}\le CB,\\
 &\|a_0\|_{L^{2}}^\eta\|L_0\|_{L^{2}}^{1-\eta-\sigma}
                                  \|Z\|_{L^{2}}^\sigma\le C_\tau B R^{N\sigma/2},\\
 &\|a_0\|_{L^{2}}^\eta\|L_0\|_{L^{2}}^{1-\eta-\sigma}
                           \|\Pi_m\|_{L^{2}}^\sigma
                             \le CB m^{\sigma/2}R^{N\sigma/2}.
\end{aligned}
\]
Each norm here is over $[\tau,S]$. For the third summand, note that
\[
 \sigma-\eta=\frac4{p(m-2)+2m},\qquad
 \sup_{[\tau,S]}L_0^{\sigma-\eta}
 \le\exp\left(\frac{C_\tau R^N}{pm}\right)\le C_\tau,
\]
discard the bounded factor $a_0^\eta$, and apply H\"older with
weights $(1-\sigma,\sigma)$ to obtain
\[
 C_\tau\|L_0\|_{L^{2}}^{1-\sigma}\|Z\|_{L^{2}}^\sigma
 \le C_\tau B R^{N\sigma/2}.
\]
Because $m^\sigma$ and $R^{N\sigma}$ are bounded,
\eqref{nsp:uni-time-factors} gives
\[
 \left(\int_\tau^S\frac{\|G-\bar G\|_{L^{\infty}}^2}{\mathscr{E}}\,\mathrm dt\right)^{1/2}
 \le C_\tau\sqrt{pR}\sqrt{R\log(e+R)}
 \le C_\tau R\log(e+R).
\]
The mean is controlled by
$\int|\bar G|^2/\mathscr{E}\le C\int D\le C$.
This proves \eqref{nsp:uni-weighted-flux}.
Apply that estimate on $[0,t]$ and use
\eqref{nsp:uni-energy-comparison} to obtain
\eqref{nsp:uni-local-flux} by Cauchy--Schwarz.

\begin{proposition}[Uniform a priori estimates]\label{uni:apriori-proposition}
For the smooth positive solutions under consideration, with
$\sigma(t)=\min\{1,t\}$,
\begin{equation}\label{nsp:uni-apriori-complete}
\begin{aligned}
 &\sup_{0\le t\le S}(\|\rho\|_{L^{\infty}}+\|u\|_{H^1}^2+A_1^2+A_3^2)
       +\int_0^S\|\sqrt\rho\,\dot u\|_{L^{2}}^2\,\mathrm dt\\
 &\quad+\sup_{0<t\le S}\sigma(t)\|\sqrt\rho\,\dot u(t)\|_{L^{2}}^2
       +\int_0^S\sigma(t)\|\nabla\dot u\|_{L^{2}}^2\,\mathrm dt\le C.
\end{aligned}
\end{equation}
The constant is independent of $S$ and the positive lower density bound.
\end{proposition}
\begin{proof}
We now close $R$ uniformly in $S$. The initial-layer
Lemma~\ref{cl:initial}, with upper horizon one, fixes $\tau_0>0$ and $M_0$
independently of $S$. If $S\le\tau_0$, the finite-time estimate already
gives the claim. Otherwise take $\tau=\tau_0/2$ and fix
$M_*>\max\{2,r_\lambda,M_0,2\|s\|_{L^{\infty}}\}$ sufficiently large
for the absorptions below, with $\theta(M_*)>0$.  Set
\[
 r_j=\theta^{-1}(2^j\theta(M_*)),\qquad t_0=\tau,
\]
and let $t_j$ be the first time the running density maximum reaches
$r_j$, considering only attained levels.  Put
$\mathcal D_j=\int_{t_{j-1}}^{t_j}D\,\mathrm dt$.
Then $\sum_j\mathcal D_j\le C$ by \eqref{nsp:uni-budget}.
On a trajectory ending at density $r_j$, let $s_j$ be its last
crossing of $r_{j-1}$ between $t_{j-1}$ and $t_j$.
On $[s_j,t_j]$, $q=\rho^\gamma-s^\gamma\ge r_{j-1}^\gamma/2$.
Since $D_t\theta(\rho)+q=-G$,
\eqref{nsp:uni-local-flux} and Young's inequality give
\begin{equation}\label{nsp:uni-level-cost}
\begin{aligned}
 \frac{\theta(r_j)}2
 &\le C_\tau r_j\log(e+r_j)
             \sqrt{t_j-s_j+\mathcal D_j}
              -\tfrac12r_{j-1}^\gamma(t_j-s_j)\\
 &\le C_\tau r_j\log(e+r_j)\sqrt{\mathcal D_j}
       +C_\tau\frac{r_j^2[\log(e+r_j)]^2}{r_{j-1}^\gamma}.
\end{aligned}
\end{equation}
Call $j$ good if
\[
 \log r_j\le\frac{\gamma+1}{2}\log r_{j-1}.
\]
On a good level the last term is bounded by
\[
 C_\tau r_j^{2/(\gamma+1)}[\log(e+r_j)]^2
       \le\frac{\theta(r_j)}4
\]
for the fixed choice of large $M_*$.
Consequently
\[
 \mathcal D_j\ge c\left(\frac{\theta(r_j)}{r_j\log(e+r_j)}\right)^2
       \ge c[\log(e+r_j)]^{2a-2}\ge\frac c{j+1}.
\]
Here $r_j\le C2^j\theta(M_*)$ follows from the lower growth of
$\theta$, and $2a-2>-1$.
If $j_\nu$ is the $\nu$th bad index, then
\[
 \left(\frac{\gamma+1}{2}\right)^\nu\log M_*
       \le\log r_{j_\nu}\le C(j_\nu+1),\qquad
       \sum_{j\ \mathrm{bad}}\frac1{j+1}<\infty.
\]
Thus every number $L$ of attained levels satisfies
$C\ge\sum_{j\le L}\mathcal D_j\ge c\sum_{j\le L}(j+1)^{-1}-C$.
The divergence of the harmonic sum bounds $L$, and therefore $R$,
independently of $S$.

With $R$ fixed, \eqref{nsp:uni-log-energy} and
\eqref{nsp:uni-energy-comparison} give
\[
 \sup_{t\le S}(A_1^2+A_3^2+\|u\|_{H^1}^2)
       +\int_0^S y\,\mathrm dt\le C.
\]
Multiply \eqref{nsp:uni-H-differential} by
$\sigma(t)=\min\{1,t\}$.  The initial-layer term is bounded because
$\int_0^{\min\{1,S\}}H\le C\int_0^{\min\{1,S\}}(1+y)\le C$.
Gronwall and \eqref{nsp:uni-budget} give
$\sup\sigma H+\int_0^S\sigma D_v\le C$.
The correction estimate also gives
\[
 \int_0^S\sigma\|\nabla Q_n\|_{L^2}^2\,\mathrm dt
 \le C\sup_{t\le S}\bigl[\sigma(t)(1+y(t))\bigr]
                         \int_0^S D\,\mathrm dt\le C.
\]
Since $\dot u=v-Q_n$, this proves \eqref{nsp:uni-apriori-complete}.
\end{proof}

\section{Global strong solutions and large-time behavior}\label{sec:construction}
\subsection{Density and velocity derivatives}

For positive approximations, constants may depend on $T$, $q$,
$\|\rho_0\|_{W^{1,q}}$, the fixed doping norm, and the bound of $\lambda'$ on the density
interval controlled by Proposition~\ref{thm:main}.
Write $A_2=\|\sqrt\rho\,\dot u\|_{L^{2}}$ in this section.

For every $1<r<q/(q-1)$, the low-order estimates and the
density-weighted Poincar\'e inequality for arbitrary vector fields give
\begin{equation}\label{str:force}
\begin{aligned}
 &\int_0^T\|\rho\dot u\|_{L^q}^r\,\mathrm dt
   +\int_0^Tt\bigl(\|\rho\dot u\|_{L^q}^2
                         +\|\dot u\|_{H^1}^2\bigr)\,\mathrm dt
       \le C_{q,r,T},\\
 &\|\nabla F\|_{L^q}+\|\nabla\omega\|_{L^q}
       \le C_q\bigl(\|\rho\dot u\|_{L^q}+\|u\|_{W^{1,q}}
                                      +\|\rho\nabla\Phi\|_{L^q}\bigr).
\end{aligned}
\end{equation}
Choose $2-2/r<\vartheta<2/q$. Weighted Poincar\'e at exponent
$2q(1-\vartheta)/(2-q\vartheta)$ and interpolation with $L^2$ give
\[
 \|\rho\dot u\|_{L^q}
 \le C_{q,r,T}\bigl(A_2+
        A_2^\vartheta\|\nabla\dot u\|_{L^2}^{1-\vartheta}\bigr).
\]
The additive term belongs to $L^r(0,T)$ because $r<2$. For the product,
H\"older's inequality gives
\[
\begin{aligned}
 &\int_0^T A_2^{r\vartheta}
               \|\nabla\dot u\|_{L^2}^{r(1-\vartheta)}\,\mathrm dt\\
 &\quad\le
 \left(\int_0^T A_2^2\,\mathrm dt\right)^{r\vartheta/2}
 \left(\int_0^Tt\|\nabla\dot u\|_{L^2}^2\,\mathrm dt
       \right)^{r(1-\vartheta)/2}
 \left(\int_0^Tt^{-r(1-\vartheta)/(2-r)}\,\mathrm dt
       \right)^{1-r/2}\le C_{q,r,T},
\end{aligned}
\]
since $r(1-\vartheta)/(2-r)<1$. The weighted assertion follows from
the same interpolation, H\"older's inequality, and the weighted acceleration.

The finite-time density bound and the Neumann equation give
$\|\rho\nabla\Phi\|_{L^{q}}\le C_{q,T}$.
At exponent two, \eqref{eng:ellipticlow} with this bounded force gives
$\|F\|_{H^1}+\|\omega\|_{H^1}\le C_T(1+A_2)$.
The div--curl estimate with its lower order term then gives, for each
fixed finite $p$,
\begin{equation}\label{str:uWp}
 \|u\|_{W^{1,p}}\le C_{p,T}(1+A_2).
\end{equation}
In particular \eqref{str:force} gives
$\|\nabla F\|_{L^{q}}+\|\nabla\omega\|_{L^{q}}
 \le C_{q,T}(1+A_2+\|\rho\dot u\|_{L^{q}})$.

We use the bounded-domain logarithmic div--curl estimate
in~\cite[proof of Lemma~5.11, following (5.43)]{FanJiangLi2024}:
\begin{equation}\label{str:BKM}
 \|\nabla u\|_{L^{\infty}}
 \le C_q\bigl(1+\|\diver u\|_{L^{\infty}}+\|\omega\|_{L^{\infty}}
                 +\|u\|_{H^1}\bigr)
                      \log(e+\|u\|_{W^{2,q}}).
\end{equation}
Put $\Xi_i=(2\mu+\lambda(\rho))\partial_i\rho$. Direct substitution
of the divergence relation in the differentiated continuity equation
gives the cancellation
\begin{equation}\label{str:Phi}
 \Dt\Xi_i+(2\mu+\lambda)\partial_i u_j\partial_j\rho
        +\Xi_i\diver u+\rho\partial_iF+\rho\partial_iP=0.
\end{equation}
The $\rho\lambda'\partial_i\rho\diver u$ terms cancel exactly.
Testing \eqref{str:Phi} gives
\[
 (\|\Xi\|_{L^{q}})'
 \le C(1+\|\nabla u\|_{L^{\infty}})\|\Xi\|_{L^{q}}
                  +C_q\|\nabla F\|_{L^{q}}.
\]
The fixed lower $L^2$ flux norms and Gagliardo--Nirenberg imply
\[
 \|F\|_{L^{\infty}}+\|\omega\|_{L^{\infty}}+\|\diver u\|_{L^{\infty}}
 \le C_{q,T}\left(1+A_2^{q/(2(q-1))}
                  +\|\rho\dot u\|_{L^{q}}^{q/(2(q-1))}\right).
\]
Differentiate $\diver u=(F+P)/(2\mu+\lambda)$ once and use the local
bound for $\lambda'$ to obtain
\[
\begin{aligned}
 \|u\|_{W^{2,q}}
 \le C_{q,T}\bigl[&\bigl(1+A_2^{q/[2(q-1)]}
              +\|\rho\dot u\|_{L^{q}}^{q/[2(q-1)]}\bigr)\|\Xi\|_{L^{q}}\\
                 &+1+A_2+\|\rho\dot u\|_{L^{q}}\bigr].
\end{aligned}
\]
Since $q/[2(q-1)]<1$, using
$(1+s^{q/[2(q-1)]})\log(e+s)\le C_q(1+s)$ in \eqref{str:BKM} gives
\[
 \|\nabla u\|_{L^{\infty}}
 \le C_{q,T}(1+A_2+\|\rho\dot u\|_{L^{q}})\log(e+\|\Xi\|_{L^{q}}).
\]
The coefficient is integrable by \eqref{str:force}. Logarithmic
Gronwall proves
\begin{equation}\label{str:densitygradient}
 \sup_{t\le T}\|\rho(t)\|_{W^{1,q}}
            +\int_0^T\|\nabla u\|_{L^{\infty}}\,\mathrm dt\le C_{q,T}.
\end{equation}
The preceding elliptic estimate gives $u\in L^r(0,T;W^{2,q})$
for every $1<r<q/(q-1)$ and $\sqrt t\,u\in L^2(0,T;W^{2,q})$.
At exponent two interpolation absorbs
$\|\nabla\rho\|_{L^{q}}\|\diver u\|_{L^{2q/(q-2)}}$, giving
\[
 \|u\|_{H^2}\le C_{q,T}(1+A_2).
\]
The estimates for $u_t=\dot u-u\cdot\nabla u$ follow from
\[
 \|u\cdot\nabla u\|_{L^{2}}^2\le C_{q,T}(1+\|u\|_{H^2}^2),
 \quad
 \|\nabla(u\cdot\nabla u)\|_{L^{2}}^2
 \le C_{q,T}(1+A_2^2+\|u\|_{W^{2,q}}^2).
\]
Finally, the continuity equation gives $\sup\|\rho_t\|_{L^{2}}\le C_{q,T}$.
The electric regularity and its time continuity follow from
\begin{equation}\label{nsp:potential-regularity}
 \|\nabla\Phi\|_{W^{2,q}}\le C_q\|\rho-b\|_{W^{1,q}},\qquad
 \|\nabla\Phi_t\|_{H^1}\le C\|\rho_t\|_{L^{2}}.
\end{equation}

\subsection{Positive local solutions and continuation}

For the preliminary positive-density construction, let the bulk
coefficient and doping profile be smooth. The final approximation
below removes this auxiliary smoothness. Local constants may depend
on the positive density lower bound, fixed doping norms, and coefficient
derivatives on the attained density interval. All transported densities
have mass $M=\int_\Om b\dd$, so their normalized Neumann
potentials are well defined.

\begin{lemma}[Linearized equation with transported coefficients]
\label{loc:linear}
Let $\rho\in C([0,t_*];H^2)$, $0<c\le\rho\le C$,
$\sup_{t\le t_*}\|\rho(t)\|_{H^2}\le L$,
and $\rho_t+\diver(\rho v)=0$. Assume that $v$ is tangent and
\[
 v\in C H^2\cap L^2H^3,\qquad
 v_t\in L^\infty L^2\cap L^2H^1,\qquad
 \rho_t\in L^\infty H^1,\qquad \|v\|_{L^\infty_{t,x}}\le V.
\]
Define
\[
 \mathcal L_\rho u=-\mu\Delta u
             -\nabla((\mu+\lambda(\rho))\diver u).
\]
For $u_0\in H^2$ satisfying both Navier trace conditions, the problem
\begin{equation}\label{loc:equation}
 \rho(u_t+v\cdot\nabla u)+\mathcal L_\rho u=-\nabla P+\rho\nabla\Phi
\end{equation}
with \eqref{nsp:boundary} has a unique solution, where at each time
$\Delta\Phi=\rho-b$, $\partial_n\Phi=0$, and $\bar\Phi=0$.
Its velocity satisfies
\[
 u\in C H^2\cap L^2H^3,\qquad
 u_t\in L^\infty L^2\cap L^2H^1.
\]
Set $Z=\int_\Om\rho(|u|^2+|u_t|^2)\dd$. Then
\begin{align}
 \|u\|_{H^2}&\le C_{L,V}(1+\|u_t\|_{L^{2}}+\|u\|_{L^{2}}),
                                                   \label{loc:H2}\\
 Z'+c_1(\|\nabla u\|_{L^{2}}^2+\|\nabla u_t\|_{L^{2}}^2)
 &\le C_{L,V}\bigl(1+\|\rho_t\|_{H^1}^2+\|v_t\|_{L^{2}}^{4/3}\bigr)(1+Z),
                                                   \label{loc:energy}\\
 \|u\|_{H^3}&\le C_L\left[\|u_t\|_{H^1}
              +(1+\|v\|_{H^2})\|u\|_{H^2}+1\right].
                                                   \label{loc:H3}
\end{align}
The constants also depend on $c,C$, the geometry and boundary coefficient,
the fixed doping norm, and the coefficient on $[c,C]$.
\end{lemma}

\begin{proof}
For tangent fields the elliptic operator has the variational form
\[
 a_\rho(w,z)=\int_\Om\bigl((2\mu+\lambda(\rho))\diver w\diver z
                         +\mu\curl w\curl z\bigr)\dd
                    +\mu\int_{\partial\Om}\alpha w_\tau z_\tau\ds.
\]
It satisfies
\[
 \|w\|_{H^1}^2\le C\bigl(a_\rho(w,w)+\|w\|_{L^{2}}^2\bigr).
\]

The flux and vorticity equations, $(2\mu+\lambda(\rho))^{-1}\in H^2$,
and the div--curl estimate give
\[
 \|u\|_{H^2}\le C_L(1+\|\rho(u_t+v\cdot\nabla u)\|_{L^{2}}+\|u\|_{H^1}).
\]
Using $\|v\|_{L^{\infty}}\le V$ and interpolating between $L^2$ and $H^2$
proves \eqref{loc:H2}.

Differentiate \eqref{loc:equation} in time. In variational form,
\begin{equation}\label{loc:Uequation}
\begin{aligned}
 \rho(u_{tt}+v\cdot\nabla u_t)+\mathcal L_\rho u_t
 ={}&-\rho_tu_t-\rho_tv\cdot\nabla u-\rho v_t\cdot\nabla u\\
    &+\nabla(\lambda'(\rho)\rho_t\diver u)-\nabla P_t
        +\partial_t(\rho\nabla\Phi).
\end{aligned}
\end{equation}
The added force satisfies, on this positive local coefficient set,
\[
 \|\rho\nabla\Phi\|_{H^1}\le C_L,\qquad
 \|\partial_t(\rho\nabla\Phi)\|_{L^{2}}
 \le \|\rho_t\|_{L^{2}}\|\nabla\Phi\|_{L^{\infty}}
       +\|\rho\|_{L^{\infty}}\|\nabla\Phi_t\|_{L^{2}}
 \le C_L\|\rho_t\|_{L^{2}}.
\]
Here $\Delta\Phi_t=\rho_t$, $\partial_n\Phi_t=0$ and
$\int_\Om\rho_t\dd=0$. Thus its pairing with $u_t$ is bounded by
$C_L(1+\|\rho_t\|_{H^1}^2)(1+Z)$, and it produces no boundary term.
The fixed boundary operator gives the Navier condition for $u_t$
in variational form. Testing by $u_t$ gives the pressure and viscosity
terms
\[
 \int_\Om P_t\diver u_t\dd,\qquad
       -\int_\Om\lambda'(\rho)\rho_t\diver u\diver u_t\dd.
\]
They are controlled by
\[
 \|P_t\|_{L^{2}}\le C_L\|\rho_t\|_{L^{2}},\qquad
 \|\lambda'(\rho)\rho_t\diver u\|_{L^{2}}
            \le C_L\|\rho_t\|_{H^1}\|u\|_{H^2}.
\]
The transport equation cancels the mass derivative and advection.
The remaining terms satisfy
\[
\begin{aligned}
 &\left|\int_\Om\rho_t|u_t|^2\dd\right|
  +\left|\int_\Om\rho_tv\cdot\nabla u\,u_t\dd\right|\\
 &\qquad\le\eps\|\nabla u_t\|_{L^{2}}^2
  +C_{\eps,L,V}(1+\|\rho_t\|_{H^1}^2)
                         (\|u\|_{H^2}^2+\|u_t\|_{L^{2}}^2),\\
 \left|\int_\Om\rho v_t\cdot\nabla u\,u_t\dd\right|
 &\le C\|v_t\|_{L^{2}}\|u\|_{H^2}\|u_t\|_{L^{2}}^{1/2}\|u_t\|_{H^1}^{1/2}\\
 &\le\eps\|\nabla u_t\|_{L^{2}}^2
       +C_{\eps,L,V}(1+\|v_t\|_{L^{2}}^{4/3})
                         (\|u\|_{H^2}^2+\|u_t\|_{L^{2}}^2).
\end{aligned}
\]
Testing \eqref{loc:equation} also by $u$, adding, and using
\eqref{loc:H2} proves \eqref{loc:energy}.

For \eqref{loc:H3}, we have
\[
 \|\rho(u_t+v\cdot\nabla u)\|_{H^1}
       \le C_L\bigl(\|u_t\|_{H^1}+\|v\|_{H^2}\|u\|_{H^2}\bigr).
\]
The Dirichlet equation
\[
 \mu\Delta\omega=\curl\bigl(\rho(u_t+v\cdot\nabla u)-\rho\nabla\Phi\bigr)
\]
with boundary value $-\alpha u_\tau\in H^{3/2}(\partial\Om)$
bounds $\omega$ in $H^2$. The identity
\[
 \nabla F=\rho(u_t+v\cdot\nabla u)-\rho\nabla\Phi-\mu\nabla^\perp\omega
\]
and the mean of
$F=(2\mu+\lambda)\diver u-P$ then bound $F$ in $H^2$.
Use the $H^2$ algebra property in two dimensions to estimate
$(F+P)/(2\mu+\lambda(\rho))$ in $H^2$, and then the higher
div--curl estimate, including its harmonic part, to obtain \eqref{loc:H3}.

Construct $u$ at times $t_n=nh$ by the implicit variational problem
\[
 \rho^n\frac{u^n-u^{n-1}}h+\mathcal L_{\rho^n}u^n
          +\rho^nv^n\cdot\nabla u^n=-\nabla P^n+\rho^n\nabla\Phi^n.
\]
For small $h$ the positive mass term absorbs the negative part of
the advection form. Lax--Milgram and elliptic regularity give
$u^n\in H^2$ with the Navier conditions. Write
$\delta_h f^n=(f^n-f^{n-1})/h$ and $U^n=\delta_hu^n$.
Define the initial $L^2$ value by
\begin{equation}\label{loc:U0}
 U^0=\rho_0^{-1}\left[-\mathcal L_{\rho_0}u_0
          -\rho_0v_0\cdot\nabla u_0-\nabla P_0+\rho_0\nabla\Phi_0\right].
\end{equation}
Subtract the preceding step, using its coefficients as the principal
coefficients, to obtain
\[
\begin{aligned}
 \rho^{n-1}\frac{U^n-U^{n-1}}h
 +\mathcal L_{\rho^{n-1}}U^n
 +\rho^{n-1}v^{n-1}\cdot\nabla U^n
 ={}&-\delta_h\rho^n U^n
       -\delta_h(\rho v)^n\cdot\nabla u^n\\
 &+\nabla(\delta_h\lambda^n\diver u^n)-\nabla\delta_hP^n
       +\delta_h(\rho\nabla\Phi)^n.
\end{aligned}
\]
For $n=1$, use \eqref{loc:U0}; only $U^0\in L^2$ is required.
Test the discrete velocity and difference equations by $u^n$ and
$U^n$. The identity $2(a-b)\cdot a=|a|^2-|b|^2+|a-b|^2$ gives the
energy differences and nonnegative time increments. The remaining
mass and advection terms have coefficients bounded by
$\|\rho_t\|_{L^\infty H^1}$ and are estimated exactly as in
\eqref{loc:energy}; \eqref{loc:H2} controls the previous-step
$H^1$ norm. The mean value formula controls the pressure and viscosity
differences by $\delta_h\rho^n$.
The electric difference obeys
\[
 \|\delta_h(\rho\nabla\Phi)^n\|_{L^{2}}
       \le C_L\|\delta_h\rho^n\|_{L^{2}},
\]
by the common mass and the Neumann estimate. Indeed, for two
bounded densities of the same mass,
\begin{equation}\label{nsp:electric-lipschitz}
\begin{aligned}
 \|\rho_1\nabla\Phi_1-\rho_2\nabla\Phi_2\|_{L^{2}}
 &\le\|\rho_1-\rho_2\|_{L^{2}}\|\nabla\Phi_1\|_{L^{\infty}}
       +\|\rho_2\|_{L^{\infty}}\|\nabla(\Phi_1-\Phi_2)\|_{L^{2}}\\
 &\le C_L\|\rho_1-\rho_2\|_{L^{2}}.
\end{aligned}
\end{equation}
Set $Z_n=\int_\Om\rho^n(|u^n|^2+|U^n|^2)\dd$. With
\[
 C_n=C_{L,V}\bigl(1+\|\rho_t\|_{L^\infty H^1}^2
                              +\|\delta_hv^n\|_{L^{2}}^{4/3}\bigr),
\]
the tested equations yield
\[
 \frac{Z_n-Z_{n-1}}h
    +c_1(\|\nabla u^n\|_{L^{2}}^2+\|\nabla U^n\|_{L^{2}}^2)
       \le C_n(1+Z_n+Z_{n-1}).
\]
Jensen's inequality gives
$h\sum_n\|\delta_hv^n\|_{L^{2}}^{4/3}
\le\int_0^{t_*}\|v_t\|_{L^{2}}^{4/3}\,\mathrm dt$, and the density node terms are
bounded by $2t_*\|\rho_t\|_{L^\infty H^1}^2$.
Since $v_t\in L^\infty L^2$, the $C_n$ are uniformly bounded.
For $h\sup_n C_n$ small, discrete Gronwall and \eqref{loc:H2}
give uniform $L^\infty H^2$, $L^\infty L^2$ bounds for $u^n,U^n$
and an $L^2H^1$ bound for $U^n$; \eqref{loc:H3} then gives the
$L^2H^3$ bound.

Use the affine interpolation $u_h$ for the time derivative and the
backward constant interpolation $\bar u_h$ for the $H^3$ bound.
Since
\[
 \|u_h-\bar u_h\|_{L^2H^1}\le h\|U_h\|_{L^2H^1}\longrightarrow0,
\]
weak compactness and coefficient continuity give a common limit
satisfying the variational equation, with $u_t\in L^2H^1$ and
$u\in L^2H^3$.
The Hilbert triple then gives $u\in C H^2$. Time differences of the
limiting equation, or the same discrete energy limit, justify
\eqref{loc:Uequation}--\eqref{loc:energy}. Linear uniqueness follows
from the weighted $L^2$ energy identity.
\end{proof}

\begin{lemma}[Positive local solution]\label{loc:existence}
For positive $\rho_0\in H^2$ and $u_0\in H^2$ satisfying the Navier
trace conditions and $\int_\Om\rho_0\dd=\int_\Om b\dd$,
a smooth nonnegative bulk coefficient and a fixed smooth doping profile give a unique
local solution
\[
 \rho,u\in C H^2,\qquad u\in L^2H^3,\qquad
 u_t\in L^\infty L^2\cap L^2H^1.
\]
The local time may depend on the positive density lower bound, the
$H^2$ data, and the coefficient on the relevant positive interval.
\end{lemma}

\begin{proof}
For tangent $v$ with $v(0)=u_0$, transport $\rho_0$ along $v$,
solve $\Delta\Phi=\rho-b$ with zero Neumann data and zero mean,
and apply Lemma~\ref{loc:linear}. Transport preserves the charge,
so this defines the map $v\mapsto u$ at every step.

Fix the density bounds and $L>\|\rho_0\|_{H^2}$ using the initial
data, and put $V=\|u_0\|_{L^{\infty}}+1$. For a radius $K$, use the ball
with $v(0)=u_0$, the Navier conditions, and
\begin{equation}\label{loc:ball}
 \|v\|_{L^\infty H^2}+\|v_t\|_{L^\infty L^2}
       +\|v_t\|_{L^2H^1}+\|v\|_{L^2H^3}\le K.
\end{equation}
The transport estimate and Sobolev embedding give
\[
 \|\rho(t)\|_{H^2}\le\|\rho_0\|_{H^2}
             \exp(C\sqrt{t_*}K),\qquad
 \|\rho_t\|_{H^1}\le C_LK.
\]
The flow formula gives the analogous positive lower and upper bounds.
Moreover,
\[
 \|v-u_0\|_{L^{\infty}}
 \le C\|v-u_0\|_{H^1}^{1/2}\|v-u_0\|_{H^2}^{1/2}
 \le Ct_*^{1/4}K^{1/2}(K+\|u_0\|_{H^2})^{1/2}.
\]
After fixing $K$, choose $t_*$ so that the density bounds and
$\|v\|_{L^{\infty}}\le V$ hold for the fixed $L,V$.

The initial value \eqref{loc:U0} is fixed by the data since $v_0=u_0$.
Equation~\eqref{loc:energy} has integrated coefficient at most
$C_{L,V}t_*(1+K^2+K^{4/3})$. Hence, when this is small,
\[
 \|u\|_{L^\infty H^2}+\|u_t\|_{L^\infty L^2}
                 +\|u_t\|_{L^2H^1}\le C_{\rm data},
 \qquad
 \|u\|_{L^2H^3}\le C_{\rm data}(1+\sqrt{t_*}K).
\]
The Navier Hodge semigroup applied to $u_0$ has the norms in
\eqref{loc:ball} bounded by $C\|u_0\|_{H^2}$ for $t_*\le1$.
Choose $K$ larger than four times the sum of this bound and the
displayed data constants, then choose $t_*$ to satisfy the smallness
requirements. The ball is nonempty and the map preserves it.

For two input velocities write $\delta v=v_1-v_2$,
$r=\rho_1-\rho_2$, and $w=u_1-u_2$. The transported density
difference satisfies
\begin{equation}\label{loc:contractrho}
 \sup_{s\le t_*}\|r(s)\|_{L^{2}}^2
       \le C_K t_*\int_0^{t_*}\|\delta v\|_{H^1}^2\,\mathrm dt.
\end{equation}
Indeed, $\nabla\rho_i\in L^\infty L^4$ and
$\int_0^{t_*}\|v_i\|_{W^{1,\infty}}\,\mathrm dt\le C\sqrt{t_*}K$.
The output difference equation is
\[
\begin{aligned}
 \rho_1(w_t+v_1\cdot\nabla w)+\mathcal L_{\rho_1}w
 ={}&-r(u_{2,t}+v_2\cdot\nabla u_2)
                  -\rho_1\delta v\cdot\nabla u_2\\
   &+\nabla((\lambda(\rho_1)-\lambda(\rho_2))\diver u_2)
                         -\nabla(P_1-P_2)
                          +(\rho_1\nabla\Phi_1-\rho_2\nabla\Phi_2).
\end{aligned}
\]
Testing by $w$, use
$u_{2,t}+v_2\cdot\nabla u_2\in L^2H^1$, the local Lipschitz bound
for $\lambda$, the force estimate \eqref{nsp:electric-lipschitz}, and
\[
 \left|\int_\Om\rho_1\delta v\cdot\nabla u_2\,w\dd\right|
 \le C\|\nabla u_2\|_{L^{\infty}}(\|\delta v\|_{L^{2}}^2+\|w\|_{L^{2}}^2).
\]
Positivity of $\rho_1$ and the div--curl estimate give, with
$D(w)=\sup_{t\le t_*}\|w(t)\|_{L^{2}}^2+\int_0^{t_*}\|w\|_{H^1}^2\,\mathrm dt$,
\[
 D(w)\le C_K(\sqrt{t_*}+t_*)D(\delta v).
\]
Here \eqref{loc:contractrho} controls the density sources, and
$\int_0^{t_*}\|\nabla u_2\|_{L^{\infty}}\,\mathrm dt\le C\sqrt{t_*}K$ controls the input
velocity source. For $t_*$ small the map contracts in the distance
$D(v_1-v_2)^{1/2}$. Weak compactness in \eqref{loc:ball} makes the
ball closed in this metric and preserves its initial value and boundary
conditions. Banach iteration gives a solution, and the difference
estimate gives uniqueness. The transport equation and the Hilbert
triple give $\rho,u\in C H^2$.
\end{proof}

To apply Proposition~\ref{thm:main}, note that
\[
 \rho_t\in L^\infty H^1,\quad
 F,\omega\in L^\infty H^1\cap H^1(0,t_*;L^2),\quad
 \dot u\in L^2H^1.
\]
Smooth positive-part truncations and the transport formula justify
the density moments. Time differences and $H^1$ traces give the flux
energy identity. For the acceleration identity, take time differences
of the variational momentum equation and test by the Steklov-regularized
corrected acceleration in \eqref{acc:commutation}. The boundary
terms are interpreted variationally and yield \eqref{acc:boundarytraction};
its factors have the required traces since $u\in L^\infty H^2$ and
the corrected acceleration belongs to $L^2H^1$.
The initial $L^2$ acceleration is finite for each positive
approximation. No uniform bound for it is required, since the
time-weighted initial terms vanish.

For positive smooth data, suppose the maximal $H^2$ time $T^*$ were
finite. The a priori estimates and \eqref{str:densitygradient} hold
uniformly on $[0,S]$, $S<T^*$. The flow formula gives
\[
 \inf_{t<T^*,x}\rho(t,x)
 \ge\inf\rho_0\exp\left(-\int_0^{T^*}\|\diver u\|_{L^{\infty}}\,\mathrm dt\right)>0.
\]
In the remainder of this continuation step constants may depend on
that lower bound and the higher initial norms. The corrected
acceleration inequality, without a time weight, gives
$\sup A_2^2+\int_0^S\|\dot u\|_{H^1}^2\,\mathrm dt\le C_*$.
Thus $\sup\|u\|_{H^2}\le C_*$ and
\[
 \int_0^S(\|\nabla^2F\|_{L^{2}}^2+\|\nabla^2\omega\|_{L^{2}}^2)\,\mathrm dt\le C_*.
\]
The vorticity boundary lift here is controlled by $u$ in $H^2$.
In this estimate the additional source has the bound
\[
 \|\rho\nabla\Phi\|_{H^1}
 \le C_T\bigl(1+\|\nabla\rho\|_{L^{2}}\bigr)\le C_{q,T}.
\]
Thus it requires no second derivative of the density and does not
change the estimate for $\nabla^2F,\nabla^2\omega$.

To control the second density derivatives, we claim that
\[
 (\|\nabla\Xi\|_{L^{2}})'
 \le C_*(1+\|\nabla u\|_{L^{\infty}})(1+\|\nabla\Xi\|_{L^{2}})
        +C_*(\|\nabla^2F\|_{L^{2}}+\|\nabla^2\omega\|_{L^{2}}).
\]
Indeed, $\nabla\Xi=(2\mu+\lambda)\nabla^2\rho
+\lambda'\nabla\rho\otimes\nabla\rho$ and $q>2$ imply
\[
 \|\nabla^2\rho\|_{L^{2}}+\|\nabla\rho\|_{L^{4}}^2
                      \le C_*(1+\|\nabla\Xi\|_{L^{2}}).
\]
For $2<q<4$, the interpolation power of $\|\nabla^2\rho\|_{L^{2}}$
in $\|\nabla\rho\|_{L^{4}}^2$ is $2-q/2<1$. Also,
\[
 \|\nabla\rho\nabla F\|_{L^{2}}
                 \le C_*(1+\|\nabla^2F\|_{L^{2}}^{2/q}),\qquad
 \|\nabla\rho\nabla^2u\|_{L^{2}}
                     \le C_*(1+\|u\|_{H^3}^{2/q}).
\]
The second derivative of $\diver u=(F+P)/(2\mu+\lambda)$ and
the div--curl estimate give
\[
 \|u\|_{H^3}\le C_*\bigl[1+\|\nabla^2F\|_{L^{2}}+\|\nabla^2\omega\|_{L^{2}}
 +(1+\|\nabla u\|_{L^{\infty}})(1+\|\nabla\Xi\|_{L^{2}})
 +\|u\|_{H^3}^{2/q}\bigr].
\]
Absorbing the last term and differentiating \eqref{str:Phi} proves
the claim. Gronwall gives $\sup\|\rho\|_{H^2}\le C_*$.
Now use $\|\nabla u\|_{L^{\infty}}
\le\delta\|u\|_{H^3}+C_\delta\|u\|_{H^2}$ in the preceding
elliptic inequality and absorb to obtain
$\|u\|_{H^3}\le C_*(1+\|\nabla^2F\|_{L^{2}}+\|\nabla^2\omega\|_{L^{2}})$.
Hence $u\in L^2H^3$.
In the twice-differentiated continuity equation, the source terms
$\rho\nabla^2\diver u$, $\nabla\rho\nabla^2u$ and
$\nabla u\nabla^2\rho$ belong to $L^1L^2$.
Transport by the integrably Lipschitz flow therefore gives the
terminal $H^2$ strong density trace. Also
$u_t=\dot u-u\cdot\nabla u\in L^2H^1$; the Hilbert triple gives
the terminal $H^2$ velocity trace. The potential has the terminal trace given by its normalized
Neumann equation. The data retain positive density, common charge,
and the Navier boundary condition. Lemma~\ref{loc:existence} therefore restarts
the local solution from these terminal data.
This contradicts maximality.

\subsection{Passage to vacuum}

For $0<\varepsilon<1$, choose a nonnegative smooth kernel $\chi$ supported in $(-1,1)$ with
integral one, and define
\[
 \lambda_\varepsilon(r)=\int_{-1}^1\chi(z)
                   \lambda(r+\varepsilon(2+z))\,\mathrm dz.
\]
These are smooth and nonnegative, converge locally in $C^1$ on
$[0,\infty)$, and satisfy
\[
 |r\lambda_\varepsilon'(r)|\le C_\lambda(1+\lambda_\varepsilon(r)).
\]
For large $r$ the shifted arguments lie in $[r,2r]$. The logarithmic
derivative condition implies
$\lambda_\varepsilon(r)\ge c\lambda(r)-C$, and hence
\[
 \theta_\varepsilon(r)\ge c\theta(r)-C\log r-C.
\]
The growth condition is preserved with common constants after
increasing a fixed threshold. Local bounds of
$\lambda_\varepsilon'$ are uniform.

Let $\Delta_N$ be the Neumann Laplacian and take
\[
 \rho_{0,\varepsilon}
 =\frac{M}{M+\varepsilon|\Om|}
       \bigl(e^{\varepsilon\Delta_N}\rho_0+\varepsilon\bigr),
 \qquad b_\varepsilon=e^{\varepsilon\Delta_N}b.
\]
These are smooth, have common mass $M$, converge in $W^{1,q}$,
and have uniform $L^\infty$ and $W^{1,q}$ bounds.
The doping satisfies $b_-\le b_\varepsilon\le b_+$.
Obtain $\Phi_{0,\varepsilon}$ from the normalized Neumann equation;
it is not an independent initial datum. For the velocity
use $u_{0,\varepsilon}=e^{-\varepsilon L_\alpha}u_0$, where $L_\alpha$
is the nonnegative Navier Hodge operator associated with the form
\[
 \int_\Om(\diver v\diver w+\curl v\curl w)\dd
                   +\int_{\partial\Om}\alpha v_\tau w_\tau\ds
\]
on tangent $H^1$ fields. The form norm plus the $L^2$ norm is equivalent
to $H^1$, so these approximations converge in $H^1$ and satisfy the
Navier conditions.
The positive construction gives global approximate solutions;
its local constants may depend on $\varepsilon$.
The stationary solutions for $b_\varepsilon$ have uniform positive
bounds and elliptic norms, and converge to $(\rho_s,\Phi_s)$.
The compactness-entry argument is uniform for this family: along a
contradicting sequence $b_\varepsilon$ and its stationary states
converge, and the same limiting hydrostatic uniqueness proof applies.

The common structural constants and local $C^1$ bounds for
$\lambda_\varepsilon$ make \eqref{eq:main}, \eqref{str:force}, and
\eqref{str:densitygradient} uniform on each finite interval, for every
fixed $1<r<q/(q-1)$.
The uniform estimate \eqref{nsp:uni-apriori-complete} also holds
with a constant independent of the interval and the approximation parameter.

Compactness gives, for $0<\eta<1-2/q$,
\[
 \rho_\varepsilon\to\rho\text{ in }C([0,T];C^\eta(\overline\Om)),
 \qquad u_\varepsilon\to u\text{ in }L^2(0,T;H^1).
\]
For the velocity first work on $[\tau,T]$ using the weighted time
derivative and $W^{2,q}$ bounds; the common $L^\infty H^1$ norm
controls the omitted interval as $\tau\downarrow0$.
The momentum is bounded in $L^\infty H^1$ and its time derivative
in $L^\infty H^{-1}$: test the conservative equation against
$H^1_0$ fields, so only the interior divergence of the uniformly
$L^2$ stress and convective flux is needed. Thus
$\rho_\varepsilon u_\varepsilon\to\rho u$ in $C([0,T];L^2)$.
Uniform convergence of $\lambda_\varepsilon(\rho_\varepsilon)$ identifies
the viscous stress. The strong velocity convergence identifies the
convective term.
The Neumann equation and the strong density convergence give
$\nabla\Phi_\varepsilon\to\nabla\Phi$ strongly, for example in
$C([0,T];W^{1,p})$ for every fixed finite $p$.
The common density bound therefore identifies
$\rho_\varepsilon\nabla\Phi_\varepsilon\to\rho\nabla\Phi$
strongly in $C([0,T];L^2)$. The electric boundary condition and
normalization are retained by this same Neumann operator.
 On each $[\tau,T]$, the weak $L^2W^{2,q}$ convergence
and the continuous first-derivative trace map pass
$\curl u_\varepsilon=-\alpha u_{\varepsilon,\tau}$ to the limit.
The normal trace follows already from the $H^1$ convergence.
Weak compactness of derivatives
gives the equations almost everywhere at positive times.

The differentiated transport formula, using
$\nabla u\in L^1L^\infty$ and $\rho\nabla\diver u\in L^1L^q$,
gives $\rho\in C W^{1,q}$. For $\delta>0$,
$\rho u/\sqrt{\rho+\delta}\in C L^2$, and its distance from
$\sqrt\rho\,u$ is at most $\sqrt\delta\|u\|_{L^{2}}$, uniformly in time.
This proves the stated weighted velocity trace. The positive-time
time derivative bounds identify $u_t$ and $\dot u$; lower
semicontinuity followed by $\tau\downarrow0$ retains all weights.
A diagonal sequence over integer time intervals and exponents
$r\uparrow q/(q-1)$ gives one global solution with all these bounds.

\subsection{Uniqueness}

Let $r=\rho_1-\rho_2$ and $w=u_1-u_2$ for two solutions with the
same initial density and momentum. Define
\[
\begin{aligned}
 Y&=\int_\Om\rho_1|w|^2\dd,\\
 D_w&=\int_\Om((2\mu+\lambda(\rho_1))(\diver w)^2+\mu(\curl w)^2)\dd
                      +\mu\int_{\partial\Om}\alpha w_\tau^2\ds.
\end{aligned}
\]
The weighted div--curl inequality gives
$\|w\|_{H^1}^2\le C(D_w+Y)$.
The density difference equation yields
\begin{equation}\label{str:diffdensity}
 \|r(t)\|_{L^{2}}^2\le Ct\int_0^t(D_w+Y)\,\mathrm ds.
\end{equation}
Here $\nabla\rho_1\in L^\infty L^q$ and
$\nabla u_2\in L^1L^\infty$.

On the common bounded density range, $P$ and $\lambda$ are Lipschitz.
Testing the momentum difference by $w$ and using the weighted
div--curl inequality and the Navier condition gives
\begin{equation}\label{str:diffenergy}
 Y'+cD_w
 \le C(1+\|\nabla u_2\|_{L^{\infty}})Y
   +C(1+\|\nabla u_2\|_{L^{\infty}}^2+\|\dot u_2\|_{L^{q}}^2)\|r\|_{L^{2}}^2.
\end{equation}
The viscosity difference contributes at most
$\eps D_w+C_\eps\|\nabla u_2\|_{L^{\infty}}^2\|r\|_{L^{2}}^2$.
The electric difference is
$r\nabla\Phi_2+\rho_1\nabla(\Phi_1-\Phi_2)$.
The densities have equal mass, so
$\|\nabla(\Phi_1-\Phi_2)\|_{L^{2}}\le C\|r\|_{L^{2}}$ and
\[
 \left|\int_\Om
   \bigl(r\nabla\Phi_2+\rho_1\nabla(\Phi_1-\Phi_2)\bigr)\cdot w\dd
                                                       \right|
 \le\eps D_w+C_\eps(\|r\|_{L^{2}}^2+Y).
\]
It therefore adds no time-singular coefficient to
\eqref{str:diffenergy}.

The strong and weighted estimates imply
\[
 \int_0^T\left[1+\|\nabla u_2\|_{L^{\infty}}+
  t(1+\|\nabla u_2\|_{L^{\infty}}^2+\|\dot u_2\|_{L^{q}}^2)\right]\,\mathrm dt<\infty.
\]
The common initial momentum and weighted velocity traces give
$Y(0+)=0$. To check the mixed weight, use
$\|(\sqrt{\rho_1}-\sqrt{\rho_2})u_2\|_{L^{2}}
\le\|\sqrt{\rho_1}-\sqrt{\rho_2}\|_{L^{\infty}}\|u_2\|_{L^{2}}\to0$. Set
$\mathcal Y(t)=\sup_{s\le t}Y(s)+\int_0^tD_w\,\mathrm ds$.
Substitute \eqref{str:diffdensity} into \eqref{str:diffenergy}, integrate
from a positive time, and let it tend to zero. The displayed
integrable coefficient and Gronwall give $\mathcal Y=0$.
Thus $\rho_1=\rho_2$, $u_1=u_2$, and the normalized Neumann
potentials coincide. This completes the existence and uniqueness
assertion. Lower semicontinuity in the uniform a priori estimates,
followed by an increasing sequence of time intervals, gives
\eqref{nsp:uniform-conclusion} and proves Theorem~\ref{nsp:thm-global}.

\subsection{Decay after the uniform density closure}

The conditional estimate \eqref{nsp:conditional-relative-decay},
with the now fixed density bound, gives exponential decay of the full
relative energy. Its energy identity gives
\begin{equation}\label{nsp:relative-energy-tail}
 \mathcal E_s(t)+\|\rho(t)-\rho_s\|_{L^{2}}^2
       +\int_t^\infty A_3^2\,\mathrm ds\le Ce^{-ct}.
\end{equation}
Since the density is bounded,
$\|P-P_s\|_{L^{2}}^2+\|P-P_s\|_{L^{4}}^4\le C\|\rho-\rho_s\|_{L^{2}}^2$.
The comparison of relative-flux energy and dissipation therefore yields
\begin{equation}\label{nsp:relative-dissipation-tail}
 \int_t^\infty\bigl(D+A_1^2\bigr)\,\mathrm ds\le Ce^{-ct}.
\end{equation}
Here $A_1,D$, and the corrected logarithmic energy are as defined in
the preceding section.

Use the corrected energy
\[
 \mathscr E=e_0+A_1^2-2C_s,\qquad
 C_s=\int_\Om(\rho-\rho_s)\nabla\Phi_s\cdot u\dd .
\]
Once the density is bounded,
\[
 |C_s|\le C\|\rho-\rho_s\|_{L^{2}} A_3
       \le\varepsilon A_1^2+C_\varepsilon e^{-ct}.
\]
Choose $a\in[t-1,t]$ with
$A_1^2(a)\le\int_{t-1}^t A_1^2$ and integrate \eqref{nsp:uni-log-differential} from $a$ to $t$.
Equations \eqref{nsp:relative-energy-tail}--%
\eqref{nsp:relative-dissipation-tail} and the preceding bound for $C_s$
give
\begin{equation}\label{nsp:viscous-decay}
 A_1^2(t)+A_3^2(t)+\int_t^\infty
                      \|\sqrt\rho\,\dot u\|_{L^{2}}^2\,\mathrm ds
       \le Ce^{-ct},\qquad t\ge1.
\end{equation}
For the integral, first integrate to a finite endpoint; the raw
$A_1^2$ there is nonnegative and the cross term tends to zero by
the velocity bound and the convergence of the density.

The tangent corrected acceleration has
$z\le2\|\sqrt\rho\,\dot u\|_{L^{2}}^2+CA_3^4$.
The acceleration inequality in the preceding section, its uniform
bound for positive times, and the same unit-interval averaging
argument yield
\[
 z(t)+\int_t^\infty D_v\,\mathrm ds\le Ce^{-ct}.
\]
The geometric correction satisfies
\eqref{nsp:uni-correction-bounds}; its $H^1$ tail is exponentially
small by \eqref{nsp:relative-dissipation-tail}. Hence
\begin{equation}\label{nsp:force-tail}
 \|\sqrt\rho\,\dot u(t)\|_{L^{2}}^2
       +\int_t^\infty\|\dot u\|_{H^1}^2\,\mathrm ds\le Ce^{-ct}.
\end{equation}
The relative electric force obeys
\[
 \|(\rho-\rho_s)\nabla\Phi_s
             +\rho\nabla(\Phi-\Phi_s)\|_{L^{p}}
       \le C_p\|\rho-\rho_s\|_{L^{p}}
\]
after fixing the density bound, with an adjustment of the harmless
finite-exponent constants. The relative elliptic estimates,
\eqref{nsp:viscous-decay}, and \eqref{nsp:force-tail} thus give
exponential $H^1$ decay of $G,\omega$, followed by $W^{1,p}$ decay of
$u$ through the div--curl estimate. Interpolation gives
exponential density decay in every finite $L^p$ space.
The Neumann equation for $\Phi-\Phi_s$ gives its $W^{1,p}$ field decay.
Finally $u_t=\dot u-u\cdot\nabla u$ supplies its weighted estimate.
For every $0<\tau<T$, the time-derivative bound gives
$u_t\in L^2(\tau,T;H^1)$ and hence $u\in C([\tau,T];H^1)$.
Weak lower semicontinuity in $W^{1,p}$ for $p>2$ extends the velocity
decay to every positive time; for $1\le p\le2$, use the $H^1$
continuity directly. Density continuity and the Neumann equation give
the same conclusion for the density and electric field.
This proves \eqref{nsp:decay-conclusion} and
\eqref{nsp:accel-decay-conclusion}, with the pointwise acceleration
assertions understood almost everywhere.

We will also use the following stronger integral form. For every
fixed $2\le p<\infty$,
\begin{equation}\label{nsp:flux-strong-tail}
 \int_t^\infty
   \bigl(\|\rho\dot u\|_{L^{p}}^2+\|G\|_{W^{1,p}}^2
                          +\|\omega\|_{W^{1,p}}^2\bigr)\,\mathrm ds
       \le C_pe^{-c_pt}.
\end{equation}
For $\rho\dot u$ use the bounded density and the weighted Poincar\'e
inequality with \eqref{nsp:force-tail}. The elliptic equations and
the already established relative electric decay give the remaining
terms. Summing over unit intervals, with $p>2$, gives
\begin{equation}\label{nsp:flux-L1-tail}
 \int_t^\infty(\|G\|_{L^{\infty}}+\|\omega\|_{L^{\infty}}+\|u\|_{L^{\infty}})
                                     \,\mathrm ds\le Ce^{-ct}.
\end{equation}
The same integrals are finite near zero by the finite-time
weighted acceleration and flux interpolation.

\section{Convergence in higher Sobolev norms and vacuum}\label{sec:higher-relaxation}
\subsection{Positive density and convergence in higher Sobolev norms}

Write $s=\rho_s$, $P_s=s^\gamma$, $\psi=\Phi-\Phi_s$,
$h(r)=2\mu+\lambda(r)$, and
$G=h(\rho)\diver u-(P-P_s)$.
The fixed stationary solution satisfies $0<c_s\le s\le C_s$ and
has bounded second derivatives. This regularity follows from the
stationary elliptic equation and $b\in W^{1,q}$, $q>2$.
We use the low-order decay established in
\eqref{nsp:decay-conclusion}, \eqref{nsp:force-tail},
\eqref{nsp:flux-strong-tail}, and \eqref{nsp:flux-L1-tail}.
Pointwise derivative assertions are understood almost everywhere.
The flux and velocity integrals in \eqref{nsp:flux-L1-tail} are also
finite near zero by the finite-time estimates. These inputs do not
require a positive density lower bound.

Assume now that $\min_{\overline\Om}\rho_0>0$.
For
\[
 \theta(r)=2\mu\log r+\int_1^r\frac{\lambda(z)}z\,\mathrm dz,
 \qquad W=\theta(\rho)-\theta(s),
\]
the continuity equation gives the exact trajectory equation
\begin{equation}\label{nsp:relative-theta-trajectory}
 D_tW+P-P_s=-G-u\cdot\nabla\theta(s).
\end{equation}
On any portion of a trajectory where $W<0$, one has $\rho<s$ and
therefore
\[
 D_tW\ge-\|G\|_{L^{\infty}}-C\|u\|_{L^{\infty}}.
\]
Integration from the initial time or the last crossing of $W=0$,
using \eqref{nsp:flux-L1-tail}, gives a uniform lower
bound for $W$. Since $\theta(s)$ is bounded and
$\theta(r)\to-\infty$ as $r\downarrow0$, this proves
\begin{equation}\label{nsp:uniform-positive-density}
 \inf_{t\ge0,\,x\in\overline\Om}\rho(t,x)>0.
\end{equation}
Together with the uniform upper bound, this fixes a compact positive
density interval. All coefficient bounds below refer to this interval.
On it,
\[
 \operatorname{sgn}(W)(P-P_s)\ge c|W|,
 \qquad |\rho-s|\le C|W|,
\]
because $P'(r)/\theta'(r)=\gamma r^\gamma/h(r)$ is strictly positive.
Consequently \eqref{nsp:relative-theta-trajectory} yields
\[
 D_t|W|+c|W|\le\|G\|_{L^{\infty}}+C\|u\|_{L^{\infty}}.
\]
Solving this scalar inequality and splitting its forcing convolution
at $t/2$ proves
\begin{equation}\label{nsp:positive-density-sup-decay}
 \|\rho(t)-s\|_{L^{\infty}}\le Ce^{-ct},\qquad t\ge1.
\end{equation}
Indeed the forcing is integrable on the full half-line and has an
exponential $L^1$ tail by \eqref{nsp:flux-L1-tail}.

We next control the relative density gradient. Define
\[
 Z=h(\rho)\nabla\rho-h(s)\nabla s.
\]
With $(\nabla u)_{ij}=\partial_j u_i$, differentiation of continuity
and the relative-flux identity gives
\begin{equation}\label{nsp:relative-density-gradient-equation}
\begin{aligned}
 D_tZ+\frac{\gamma\rho^\gamma}{h(\rho)}Z
 ={}&-(\nabla u)^{\mathsf T}\bigl(Z+h(s)\nabla s\bigr)
       -\bigl(Z+h(s)\nabla s\bigr)\diver u
       -\rho\nabla G\\
 &+\left(\gamma\rho s^{\gamma-1}
              -\frac{\gamma\rho^\gamma h(s)}{h(\rho)}\right)\nabla s
       -u\cdot\nabla\bigl(h(s)\nabla s\bigr).
\end{aligned}
\end{equation}
In particular, the transpose in the deformation term is required
with this gradient convention. The stationary term in parentheses
vanishes at $\rho=s$ and is bounded by $C|\rho-s|$ on the fixed
positive interval.
Since $\gamma\rho^\gamma/h(\rho)\ge c>0$,
testing \eqref{nsp:relative-density-gradient-equation} in $L^q$ gives
\begin{equation}\label{nsp:relative-gradient-damping}
 \frac{\mathrm d}{\mathrm dt}\|Z\|_{L^{q}}+c\|Z\|_{L^{q}}
       \le C\|\nabla u\|_{L^{\infty}}\|Z\|_{L^{q}}
                         +C_q\|\rho\dot u\|_{L^{q}}+Ce^{-ct}.
\end{equation}
Here the terms with the fixed gradient $\nabla s$ are controlled by
$C\|\nabla u\|_{L^{q}}+C\|u\|_{L^{q}}$, which decay by
\eqref{nsp:decay-conclusion}. The relative force
\[
 (\rho-s)\nabla\Phi_s+\rho\nabla\psi
\]
also decays in every fixed finite $L^p$ space. Thus the elliptic
flux and vorticity estimates give
\begin{equation}\label{nsp:positive-flux-gradient}
 \|\nabla G\|_{L^{q}}+\|\nabla\omega\|_{L^{q}}
                 +\|G\|_{L^{\infty}}+\|\omega\|_{L^{\infty}}
       \le C_q\bigl(\|\rho\dot u\|_{L^{q}}+e^{-ct}\bigr).
\end{equation}

It remains to control the Lipschitz coefficient in
\eqref{nsp:relative-gradient-damping}. We have
\[
 \nabla\rho-\nabla s
       =\frac{Z+(h(s)-h(\rho))\nabla s}{h(\rho)},
\]
and consequently
\begin{equation}\label{nsp:relative-gradient-comparison}
 \|\nabla\rho-\nabla s\|_{L^{q}}
       \le C\bigl(\|Z\|_{L^{q}}+\|\rho-s\|_{L^{\infty}}\bigr),
 \qquad \|\nabla\rho\|_{L^{q}}\le C(1+\|Z\|_{L^{q}}).
\end{equation}
Differentiate $\diver u=(G+P-P_s)/h(\rho)$.
Using \eqref{nsp:positive-density-sup-decay},
\eqref{nsp:positive-flux-gradient}, and the div--curl estimate gives
\begin{equation}\label{nsp:positive-second-velocity-estimate}
 \|u\|_{W^{2,q}}
       \le C_q\left[
        (1+\|\rho\dot u\|_{L^{q}})\|Z\|_{L^{q}}
                         +\|\rho\dot u\|_{L^{q}}+e^{-ct}\right].
\end{equation}
The stationary gradient enters only multiplied by the decaying
divergence or by $\rho-s$; it does not produce a constant source.

Apply the logarithmic div--curl inequality to $e^tu$ and scale back.
Its logarithm is bounded using
\eqref{nsp:positive-second-velocity-estimate}. One obtains
\begin{equation}\label{nsp:positive-integrable-lipschitz}
\begin{aligned}
 \|\nabla u\|_{L^{\infty}}
       &\le b_q(t)\bigl[1+\log(e+\|Z\|_{L^{q}})\bigr],\\
 b_q(t)&=C_q(1+t)\left(
                  e^{-ct}+\|\rho\dot u\|_{L^{q}}+\|\rho\dot u\|_{L^{q}}^2
                         \right),\\
 \int_t^\infty b_q(s)\,\mathrm ds&\le C_qe^{-ct},\qquad t\ge1.
\end{aligned}
\end{equation}
The last assertion follows from
\eqref{nsp:flux-strong-tail} by Cauchy--Schwarz on successive
unit intervals. The constants here do not depend on $\|Z\|_{L^{q}}$.
Discard the damping term in \eqref{nsp:relative-gradient-damping}
and divide by $e+\|Z\|_{L^{q}}$. Logarithmic Gronwall, using
\eqref{nsp:positive-integrable-lipschitz}, first gives
$\sup_{t\ge1}\|Z(t)\|_{L^{q}}<\infty$.
It then follows that
\[
 \int_t^\infty\|\nabla u\|_{L^{\infty}}\,\mathrm ds\le Ce^{-ct}.
\]
Returning to the damping inequality and using the exponential
$L^1$ tail of $\|\rho\dot u\|_{L^{q}}$, supplied by
\eqref{nsp:flux-strong-tail}, proves
\[
 \|Z(t)\|_{L^{q}}\le Ce^{-ct}.
\]
Equations \eqref{nsp:positive-density-sup-decay} and
\eqref{nsp:relative-gradient-comparison} therefore give
\begin{equation}\label{nsp:positive-density-gradient-decay}
 \|\rho(t)-s\|_{W^{1,q}}\le Ce^{-ct},\qquad t\ge1.
\end{equation}
In particular it is $\nabla\rho-\nabla s$ that tends to zero.

The second-order $L^2$ div--curl estimate gives
\[
 \|u\|_{H^2}\le Ce^{-ct}
       +C\|\nabla\rho\|_{L^{q}}
                          \|\diver u\|_{L^{2q/(q-2)}}.
\]
The factor $\|\nabla\rho\|_{L^{q}}$ is uniformly bounded. By
\eqref{nsp:decay-conclusion} with $p=2q/(q-2)$, the divergence
factor decays exponentially. Hence $\|u(t)\|_{H^2}\le Ce^{-ct}$.
The positive density bound gives
$\|\dot u\|_{L^{2}}\le C\|\sqrt\rho\dot u\|_{L^{2}}$.
Using $u_t=\dot u-u\cdot\nabla u$ and
$\rho_t=-u\cdot\nabla\rho-\rho\diver u$, we conclude that
\begin{equation}\label{nsp:positive-higher-relaxation}
\begin{aligned}
 &\|\rho(t)-s\|_{W^{1,q}}+\|\rho_t(t)\|_{L^{q}}
                   +\|u(t)\|_{H^2}+\|u_t(t)\|_{L^{2}}\le Ce^{-ct},\\
 &\int_t^\infty\left(
       \|u\|_{W^{2,q}}^2+\|u_t\|_{H^1}^2
                              +\|\nabla u\|_{L^{\infty}}\right)\,\mathrm ds
       \le Ce^{-ct},\qquad t\ge1.
\end{aligned}
\end{equation}
For the second line, use
\eqref{nsp:positive-second-velocity-estimate},
$\|u\cdot\nabla u\|_{H^1}\le C\|u\|_{H^2}^2$, and the
already proved $L^2$ tail of $\dot u$ in $H^1$.
The elliptic equation also gives
\[
 \|\nabla(\Phi-\Phi_s)\|_{W^{2,q}}
                +\|\nabla\Phi_t\|_{W^{1,q}}\le Ce^{-ct}.
\]
The density assertions and the tail assertions hold for every $t$;
the pointwise velocity derivative assertions hold almost everywhere.
The constants in this subsection may depend on the initial
$W^{1,q}$ density norm and on $(\min\rho_0)^{-1}$.

\subsection{Persistent vacuum and growth of density gradients}

Suppose instead that $\rho_0(a_0)=0$ for some
$a_0\in\overline\Om$. On every finite interval the strong solution
has $\nabla u\in L^1L^\infty$ and generates an invertible Lipschitz
flow preserving $\overline\Om$. Along this flow,
\[
 \rho(t,X(t,a))=\rho_0(a)
       \exp\left(-\int_0^t\diver u(\xi,X(\xi,a))\,\mathrm d\xi\right).
\]
Thus $\rho(t,X(t,a_0))=0$, and
\begin{equation}\label{nsp:persistent-vacuum}
 \min_{\overline\Om}\rho(t)=0,\qquad
       \|\rho(t)-s\|_{L^{\infty}}\ge\min_{\overline\Om}s>0,
                                  \qquad t\ge0.
\end{equation}

For $2<r\le q$, the function $\rho-s$ has zero spatial mean.
The mean-zero Gagliardo--Nirenberg inequality gives
\[
 \min_{\overline\Om}s
 \le\|\rho-s\|_{L^{\infty}}
 \le C_r\|\nabla\rho-\nabla s\|_{L^{r}}^{r/[2(r-1)]}
                  \|\rho-s\|_{L^{2}}^{(r-2)/[2(r-1)]}.
\]
The exponential $L^2$ convergence therefore implies
\[
 \|\nabla\rho(t)-\nabla s\|_{L^{r}}\ge C_r^{-1}e^{c_rt}.
\]
Subtracting the fixed norm $\|\nabla s\|_{L^{r}}$ proves exponential growth
of $\|\nabla\rho(t)\|_{L^{r}}$ for all sufficiently large $t$.
To cover the remaining bounded interval, use the fixed positive mass.
Writing $\bar\rho=|\Om|^{-1}\int_\Om s\dd>0$, the persistent vacuum
point and the mean-zero Sobolev inequality give
\[
 \bar\rho\le\|\rho(t)-\bar\rho\|_{L^{\infty}}
                         \le C_r\|\nabla\rho(t)\|_{L^{r}}.
\]
After reducing the positive prefactor, this proves
\begin{equation}\label{nsp:vacuum-gradient-growth}
 \|\nabla\rho(t)\|_{L^{r}}\ge C_r^{-1}e^{c_rt},
                      \qquad t\ge0,\quad 2<r\le q.
\end{equation}
Finally, for $0\le\varepsilon<1$,
\begin{equation}\label{nsp:vacuum-set-decay}
 |\{x\in\Om:\rho(t,x)\le\varepsilon s(x)\}|
       \le\frac{\|\rho(t)-s\|_{L^{2}}^2}
                    {(1-\varepsilon)^2(\min s)^2}
       \le C_\varepsilon e^{-ct}.
\end{equation}

\section*{Acknowledgements and AI disclosure}

The authors acknowledge the use of generative AI tools in developing and checking mathe
matical arguments and in drafting and revising this \textbf{BASIC} version of manuscript.

\end{document}